\documentclass{amsart}
\begin{document}
\allowdisplaybreaks
        \newtheorem{theorem}{{\bf Theorem}}
          \newtheorem{lemma}{Lemma}
      \newtheorem{corollary}[lemma]{Corollary}
     \newtheorem{definition}{Definition}
        \newtheorem{problem}{Problem}
     \newtheorem{conjecture}{Conjecture}
    \newtheorem{proposition}{{\bf Prop.}}
         \newtheorem{remark}{Remark}
\def\con#1{\setbox13\hbox{$#1$}\ifdim\wd13<1.1em\breve{#1}\else{\(#1\)}\breve{\ }\fi}%
\def\halfthinspace{\relax\ifmmode\mskip.5\thinmuskip\relax\else\kern.8888em\fi}%
\let\hts=\halfthinspace%
\def\rp{{\hts;\hts}}%
\def\rs{\mathop\dagger}%
\def\id{{1\kern-.08em\raise1.3ex\hbox{\rm,}\kern.08em}}%
\def\di{{0\kern-.04em\raise1.3ex\hbox{\rm,}\kern.04em}}%
\def\makecs#1#2{\makecsX {#1}#2,.}%
\def\makecsX#1#2#3.{\onecs{#1}{#2}%
\ifx#3,\let\next\eatit\else\let\next\makecsX\fi\next{#1}#3.}%
\def\onecs#1#2{\expandafter\gdef\csname #2\endcsname%
{{\csname #1\endcsname {#2}}}}%
\def\eatit#1#2.{\relax}%
\makecs{}      {abcdefghijklmnopqrstuvwxyzABCDEFGHIJKLMNOPQR TUVWXYZ}%
\def\atsymbol{\char'100}%
\def\ie{{\it i.e.}}%
\def\eg{{\it e.g.}}%
\def\conv#1{\setbox13\hbox{$#1$}\ifdim\wd13<1.5em{#1}^{-1}%
\else{(#1)}^{-1}\fi}%
\def\({\left(}%
\def\){\right)}%
\def\<{\left<}%
\def\>{\right>}%
\def\Id{\mathsf{Id}\hts}%
\def\Di{\mathsf{Di}\hts}%
\def\ex#1{\mathbf\exists_{#1}}%
\def\all#1{\mathbf\forall_{#1}}%
\def\gc#1{\mathfrak{#1}}
\def\etc{{\it etc}}%
\def\RA{\mathsf{RA}}%
\let \ism=\cong%
\def\Cm#1{\setbox13\hbox{$#1$}\ifdim\wd13=0pt{\gc{Cm}}%
\else{\gc{Cm}\({#1}\)}\fi}%
\def\mand{\mathrel{\ \text{and}\ }}
\def\mor{\mathrel{\ \text{or}\ }}
\def\alg#1#2{\text{$\mathsf{#1}_{#2}$}} \def\MP{\mathsf{MP}}%
\def\MP{\text{\it modus ponens}}%
\makeatletter \def\labelenumi{(\@roman\c@enumi)}%
\def\theenumi{(\@roman\c@enumi)}%
\def\labelenumii{(\@alph\c@enumii)}%
\def\theenumii{\@alph\c@enumii}%
\def\labelenumiii{(\@arabic\c@enumiii)}%
\def\theenumiii{(\@arabic\c@enumiii)}%
\def\p@enumiii{\theenumi\theenumii}%
\def\wideA#1{\@mathmeasure\z@\textstyle{#1}\ifdim\wd\z@>\tw@%
  em\mathaccent "055C{#1}%
  \else\mathaccent "0364{#1}\fi}%
\def\wideB#1{\@mathmeasure\z@\textstyle{#1}\ifdim\wd\z@>\tw@%
  em\mathaccent "055E{#1}%
  \else\mathaccent "0366{#1}\fi}%
\makeatother%
\def\ul#1{#1}
\def\ull#1{#1}%
\def\from{\leftarrow}%
\def\rel{{\mathbf{R}}}%
\def\rms{relevant model structure}%
\def\iff{\,\mathrel\leftrightarrow\,}%
\def\ifff{\,\,\text{iff}\,\,\,}%
\def\min#1{\overline{#1}}%
\def\conmin#1{\min{\conv{#1}}}%
\def\CR{{\mathbf{CR}}}%
\def\rmin{{\sim}}%
\def\nec{\square}%
\def\true{\mathbf\t}%
\def\cid{\text{\sf c.i.d.}}%
\def\toE{{\to}\E}%
\def\sent{\mathsf{Sent}}%
\def\var{\mathsf{Pv}}%
\def\fmla{\mathsf{Fmla}}%
\def\sqnt{\mathsf{Sqnt}}%
\def\blank{\phantom{I}}%
\def\DR{\text{\bf{DR2}}}%
\def\bxd#1{$\boxed{\text{\bf#1}}$}%
\def\L{{\mathcal L}}
\def\RR{\text{\sf R}}%
\def\BB{\text{\sf B}}%
\def\pow#1{\mathcal{P}\({#1}\)}%
\def\wid{{\text{\small[\id]}}}%
\def\prop#1{\mathfrak{Pr}(#1)}%
\def\Ap{\A'}%
\def\Bp{\B'}%
\def\Cp{\C'}%
\def\Dp{\D'}%
\def\crc#1#2{(#1\circ#2)}%
\def\proves{\vdash}%
\def\provesFL{\vdash^{\mathbf{4}}}%
\def\provesFl{\vdash_4}%
\def\bnd{\par\smallskip\hrule\smallskip\par}%
\def\bnd{\par\medskip\hrule\medskip\par}%
\def\0{0}%
\def\1{1}%
\def\2{2}%
\def\3{3}%
\def\4{4}%
\def\fm#1#2#3{\setbox13\hbox{$#1$}\ifdim\wd13<1.2em{#1{}_{#2#3}}%
\else{\(#1\){}_{#2#3}}\fi}%
\def\implies{\,{\Rightarrow}\,}
\def\|{{|}}%
\def\qqq{\ \,}%
\def\stand{\hbox{\vrule height12pt depth0pt width0pt}}%
\def\BL{\mathcal{B}}%
\def\FL{{\mathcal{L}_4}}%
\def\EL{{\mathcal{E}_4}}%
\def\oo{\emptyset}%
\def\rationals{\mathbb{Q}}%
\def\four{\{\0,\1,\2,\3\}}%
\def\basis#1{\B_4(#1)}%
\def\Qt{\rationals^2}%
\def\convax#1{\con\Xi(#1)}%
\def\Ka{\mathfrak{K}}%
\def\ra#1{\textup{\rm R$_{#1}$}}%
\let\bp=\cdot%
\def\lref#1{L\eqref{#1}}%
\title[{Tarski's relevance logic; Version~2}]{Tarski's relevance
  logic; Version~2} \author{Roger D.  Maddux}
\address{Department of Mathematics\\ 396 Carver Hall\\
  Iowa State University\\ Ames, Iowa 50011\\ U.S.A.}
\email{maddux\atsymbol iastate.edu} \date\today \subjclass{Primary:
  03G15, 03B47} \keywords{relevance logic, relation algebra}
\begin{abstract}
  Tarski's relevance logic is defined and shown to contain many
  formulas and derived rules of inference. The definition arises from
  Tarski's work on first-order logic restricted to finitely many
  variables. It is a relevance logic because it contains the Basic
  Logic of Routley-Plumwood-Meyer-Brady, has Belnap's variable-sharing
  property, and avoids the paradoxes of implication. It does not
  include several formulas used as axioms in the Anderson-Belnap
  system $\RR$. For example, the Axiom of Contraposition is not in
  Tarski's relevance logic. On the other hand, the Rules of
  Contraposition and Disjunctive Syllogism are derived rules of
  inference in Tarski's relevance logic.  It also contains a formula
  (not previously known or considered as an axiom for any relevance
  logic) that provides a counterexample to a completeness theorem of
  T.\ Kowalski (that the system $\RR$ is complete with respect to the
  class of dense commutative relation algebras).
\end{abstract}
\maketitle
\section{Introduction}
In 1975, Alfred Tarski delivered a pair of lectures on relation
algeras at the University of Campinas. The lectures were videotaped
and transcriptions of them appeared in the book {\bf Alfred Tarski:
  Lectures at UniCamp in 1975} published in 2016. At the end of his
second lecture, Tarski said (p.\ 154),
\begin{quote}
  ``And finally, the last question, if it is so, you could ask me a
  question whether this definition of relation algebra which I have
  suggested and which I have founded --- I suggested it many years ago
  --- is justified in any intrinsic sense.  If we know that these are
  not all equations which are needed to obtain representation
  theorems, this means, to obtain the algebraic expression of
  first-order logic with two-place predicate, if we know that this is
  not an adequate expression of this logic, then why restrict oneself
  to these equations? Why not to add strictly some other equations
  which hold in representable relation algebras or maybe all?''
\end{quote}
Tarski's question arises from the fact that the equations he chose for
his axiomatization of relation algebras are all simple and natural and
occur throughout the nineteenth century literature on the algebra of
logic, such as the works of Peirce~\cite{ MR1505381, MR1403570,
  MR1403566, MR1403567, MR1403568, MR1403569, MR1961048, MR0110632,
  MR1188412, MR1625634} and Schr\"oder~\cite{ MR1509947, MR1510638,
  MR0192997, MR0192998, MR0192999} and yet the choice is clearly
arbitrary. Furthermore, the axioms were shown to be incomplete, hence
insufficient for proving representability, by Lyndon \cite{MR0037278}
in 1950. Back in 1941 Tarski \cite[pp.\,87--88]{Tarski1941} wrote,
\begin{quote}
  ``Is it the case that every sentence of the calculus of relations
  which is true in every domain of individuals is derivable from the
  axioms adopted under the second method? This problem presents some
  difficulties and remains open. I can only say that I am practically
  sure that I can prove with the help of the second method all of the
  hundreds of theorems to be found in Schroder's {\it Algebra und
    Logik der Relative}.''
\end{quote}
The ``second method'' is Tarski's equational axiomatization.  The
problem Tarski posed was to find a true equation that his axioms
\emph{can't} prove.  Lyndon solved Tarski's open problem in his 1950
paper by showing the answer is ``no''. This left only Tarski's rather
practical reason for adopting his axioms: they are good enough to
prove a lot. 

Besides what could be \emph{proved} from his axioms, Tarski was also
concerned from the outset with what could be \emph{expressed} with
equations. This topic had been considered already by Schr\"oder and
L\"owenheim~\cite{ MR1511558, MR1511730, MR1511835, MR1511923,
  MR0001929, MR0018624, MR2374248}.  By the early 1940s Tarski had
proved that the equations of relation algebras have the same
expressive power as first-order logic restricted to three variables.
Tarski took a first-order language with an equality symbol and other
binary relation symbols (but no function symbols or constants),
reduced the usual stock of countably many variables to just three,
added a binary operator $|$ on relation symbols, and included a
definition asserting that the operator produces the relative product
of the relations denoted by the inputs:
\begin{equation}\label{relprod}
  (\A|\B)(\x,\y)\iff\exists\z(\A(\x,\z)\land\B(\z,\y)).
\end{equation}
He included other operators on relation symbols, for union,
complementation, and converse, along with their definitions
\begin{align}
  \label{union}
  (\A\cup\B)(\x,\y)&\iff\A(\x,\y)\lor\B(\x,\y),\\
  \label{comp}
  \min\A(\x,\y)&\iff \neg\A(\x,\y),\\
  \label{conv}
  \conv\A(\x,\y)&\iff\A(\y,\x).
\end{align}
Finally, Tarski introduced a new form of sentence called an equation,
written $\A=\B$, made out of two relation symbols $\A$ and $\B$ and a
new equality symbol, with this definition
\begin{equation*}
  \A=\B\iff\forall\x\forall\y(\A(\x,\y)\iff\B(\x,\y)).
\end{equation*}
In Tarski's definition for $|$, $\z$ is the first variable distinct
from $\x$ and $\y$. Such a variable always exists because Tarski's
language has three variables. To illustrate, the associative law for
relative multiplication is
\begin{equation*}
  (\A|\B)|\C=\A|(\B|\C),
\end{equation*}
and its expansion according to the definition of $|$ is
\begin{equation*}
  \forall\x\forall\y\Big(
  \exists\z\big(\exists\y(\A(\x,\y)\land\B(\y,\z))\land\C(\z,\y)\big)
  \iff
  \exists\z\big(\A(\x,\z)\land\exists\x(\B(\z,\x)\land\C(\x,\y))\big)
  \Big).
\end{equation*}
The burden of parentheses can be reduced by resorting to subscripts.
\begin{equation*}
  \forall\x\forall\y\Big(
  \exists\z\big(\exists\y(\A_{\x,\y}\land\B_{\y,\z})\land\C_{\z,\y}\big)
  \iff
  \exists\z\big(\A_{\x,\z}\land\exists\x(\B_{\z,\x}\land\C_{\x,\y})\big)
  \Big).
\end{equation*}
Tarski observed that every relation-algebraic equation expands to a
formula in first-order logic of binary relations restricted to three
variables, as was just done for the associative law, and then he
proved that every formula is equivalent to such an expansion, \ie,
every formula of 3-variable first-order logic (of binary relations)
can be converted to an equivalent relation-algebraic equation. For
details, consult \cite{Maddux2006, MR920815}.

Naturally, Tarski included the associative law as an axiom for
relation algebras.  With regard to the other axioms, Tarski found that
he could not only express them with three variables, but also prove
them with only three variables.  On the other hand, Tarski's proof of
the associative law used four variables. Could it be proved with only
three variables?  J.\ C.\ C.\ McKinsey had invented an algebra that
satisfies all of Tarski's axioms for relation algebras except the
associative law, thus proving that the associative law is independent
of the other axioms.  Tarski used McKinsey's algebra to prove that the
associative law for relative multiplication cannot be proved in
first-order logic with only three variables.

This is how things stood in 1975, when Tarski asked, ``whether this
definition of relation algebra ...  is justified in any intrinsic
sense''. Tarski had proved that every equation true in all relation
algebras, \ie, every equation that follows from his axioms by the
rules of equational logic (equality is transitive and symmetric, and
equals may be substituted for equals) can be proved in first-order
logic with four variables. Since the associative law is the only axiom
requiring four variables to prove, Tarski asked whether deleting it
would result in an equational theory equivalent to 3-variable logic in
means of proof as well as expression. If not, could the associative be
replaced with a weaker version to yield an equational theory
equivalent to 3-variable logic? 

These problems were included in the draft of the Tarski-Givant book
\cite{MR920815}, which was being written at the time of Tarski's talk.
This book started life as an unpublished manuscript by Tarski from the
early 1940s. Work on the revision was begun in 1971.  It was planned
to become Tarski's contribution to the Proceedings of the Tarski
Symposium~\cite{MR601214, MR0349335}, held in honor of his 70th
birthday, but grew into a project not published until four years after
his death.

Around this time of Tarski's talk it was proved that the answers are
``no'' and ``yes'', \ie, deleting the associative law leaves an axiom
set that is too weak, but a weakened associative law, dubbed the
``semi-associative law'' can replace the associative law to produce an
equational theory that is a precise correlate of first-order logic of
binary relation symbols and only three individual variables---every
sentence of 3-variable logic is equivalent to an equation, and every
provable sentence of 3-variable logic is equivalent to an equation
provable from the weakened axiom set.  (Algebras satisfying this
weakened axiom set are now called semi-associative relation algebras,
but their initial name was ``Tarski algebras''.)  Furthermore, the
equations true in all relation algebras are exactly those that are
equivalent to a statement in 3-variable logic of binary relations and
can be proved with four variables. For details see~\cite{ Maddux1978,
  Maddux1989, Maddux2006, MR920815}.

This last result provides a potential answer to Tarski's question,
``whether this definition of relation algebra ...  is justified in any
intrinsic sense''. The justification of Tarski's axioms would be that
their consequences are the equations that are
\begin{itemize}
\item equivalent to statements in first-order logic of binary
  relations, restricted to three variables, and
\item are provable with four variables.
\end{itemize}
Certainly one can dispute whether this characterization is
``intrinsic'', but any true equation not provable with four variables
must require at least five, and finding such formulas is difficult.
The shortest ones known are quite complicated. It is a safe bet that
no such formula was ever encountered for any other purpose prior to
Lyndon's proof that Tarski's axioms are incomplete.

Tarski never published his proof that the associative law requires
four variables to prove; see \cite[p.\ 65]{MR0124250}.  Henkin
\cite{MR0392564} published such a proof, but for cylindric algebras
rather than relation algebras. The connections between these two
subjects had been studied from the early 1960s by Monk
\cite{Monk1961a, Monk1961}.

A search for Henkin \cite{MR0392564} led to the same volume containing
Routley-Meyer \cite{MR0409114}.  What Routley and Meyer define as a
``relevant model structure'' in that paper was immediately recognized
as nearly the same as the atom structure of an integral relation
relation, but with one property missing and two more added.  The atom
structures of relation algebras with the two additional properties
(density and commutativity) form particularly nice relevant model
structures. They have two other additional properties, one called
``normal'' by Routley and Meyer, the other called ``tagging'' by Dunn;
for more details, see \cite[\S7]{MR2641636}.

Indeed, making use of the database of finite relation algebras
compiled for \cite{Maddux2006}, one can see that out of 4527 integral
relation algebras with five or fewer atoms, all of them are
``normal'', all of them have ``tagging'', 3885 of them are commutative
(satisfy $\x\rp\y=\y\rp\x$), 822 of them are dense (satisfy
$\x\leq\x\rp\x$), and 626 are both commutative and dense. Many of
these 626 relevant model structures are the atom structures of proper
relation algebras. The elements of proper relation algebras are binary
relations, and their operations are the usual set-theoretic operations
on relations: union, intersection, complementation (with respect to
the largest relation), relative multiplication (or composition), and
conversion (forming the converse of a binary relation). This allows
Routley-Meyer semantics to be deciphered into ordinary mathematical
concepts in common use.

Routley and Meyer refer to the objects in a relevant model structure
at first as ``worlds'', but settle on ``set-ups'' (which might be, as
they suggest, ``sets of beliefs''). Other words have been employed on
these objects, such as ``situations'' or ``points''. However, in
relevant model structures arising from proper relation algebras, the
set-ups (or worlds, or situations, or points) are clearly identified;
they are simply binary relations.

The logical connectives considered by Routley and Meyer are
conjunction $\land$, disjunction $\lor$, negation $\rmin$, and
implication $\to$ \cite[p.\ 204, \S1]{MR0409114}.  Every
\emph{valuation} $v$ determines a map that sends each propositional
variable and set-up to a truth value, either $\T$ or $\F$ \cite[p.\
206, \S3]{MR0409114}. A valuation extends to an \emph{interpretation}
$\I$ defined for all formulas and set-ups. An interpretation in turn
determines a map, we call it $\J$, from formulas to sets of set-ups.
Conditions ii and iii \cite[p.\ 206]{MR0409114} defining the extension
show how the connectives are interpreted: conjunction as intersection
and disjunction as union.  The treatment of negation involves (what
has become known as) the Routley star. The Routley star in a relevant
model structure matches up with the unary operation of forming the
converse of the atom in the atom structure of a relation algebra.  In
a proper relation algebra, this is simply the ordinary converse of a
binary relation---the result of turning all the pairs around.
Condition vi \cite[p.\ 206]{MR0409114} shows that negation is to be
treated as the converse of the complement (or, what is the same thing,
the complement of the converse). We call this simply
converse-complementation.  Condition v \cite[p.\,206]{MR0409114} shows
that the binary connective $\circ$, defined by
\begin{equation}\label{D1}
  \tag{D1}  \A\circ\B=\rmin(\A\to\rmin\B)
\end{equation}
in \cite[p.\ 204]{MR0409114} (later called ``fusion''), is interpreted
as relative multiplication in the opposite order, that is,
\begin{equation}\label{fusion}
  \A\circ\B=\B|\A.
\end{equation}
This is a good place to notice a notational coincidence. In the case
where $\A$ and $\B$ are unary functions, \eqref{fusion} shows that
$\circ$ denotes the usual operation of composing these two functions.
The use of $\circ$ for functional composition is common in a wide
range of mathematical literature, including calculus textbooks.
Instead of writing $\<\x,\y\>\in\A$ in case $\A$ is a function, it is
customary to write $\A(\x)=\y$, since there is no other ordered pair
in $\A$ whose first component is $\x$, \ie, $\y$ is uniquely
determined by $\A$ and $\x$.  Composing $\A$ and $\B$ produces a
function denoted $\A\circ\B$, defined by
\begin{equation*}
  (\A\circ\B)(\x)=\A(\B(\x)).
\end{equation*}
This is an abbreviated way of describing the relative product of $\A$
and $\B$ in the opposite order. It says, in more detail, that
$\<\x,\B(\x)\>\in\B$ and $\<\B(\x),\A(\B(\x))\>\in\A$. Combining these
two statements according to \eqref{relprod} yields
$\<\x,\A(\B(\x))\>\in\B|\A$, establishing \eqref{fusion} in case $\A$
and $\B$ are functions. The notational coincidence is that the same
symbol was (inadvertantly, as it turns out) chosen for the same thing.

A discussion of definition \eqref{D1}, incorporating remarks of
Anderson, Nelnap, Dunn, Woodruff, and Meyer, occurs in
\cite[\S27.1.4]{MR0406756}, where the ``memorable and delightful''
properties of $\circ$ are mentioned, including associativity (see
Lemma~\ref{assocfusion} below).  However, they ask
\cite[p.\,345]{MR0406756},
\begin{quote} 
  ``3. How then to interpret $\circ$? We confess puzzlement.

  In some ways $\circ$ looks like conjunction \dots\ 

  But $\circ$ fails to have the property $\A\circ\B\to\A$; so it isn't
  conjunction.''
\end{quote}
Perhaps the proper interpretation of $\circ$ is identified in
\eqref{fusion}. In this context Meyer's remarks seem remarkably
insightful:
\begin{quote}
  ``The term 'fusion' is, I believe, due to Fine, and it is a good
  one; previous tries were 'intensional conjunction', 'relevant
  conjunction', 'consistency', and 'cotenability'. But the first two
  invite confusion with the extensional conjunction '\&', while the
  latter two depend on properties of the negation-of R that have not,
  so far, generalized to related logics. The notion, in one guise or
  another, has been invented or re-invented by Lewis, Nelson,
  Anderson-Belnap, Church, Dunn, Curry, Meredith, Powers, Routley,
  Urquhart, Fine and the author, no doubt among several score others.
  It is to be attributed accordingly to Tarski, on the ground that,
  when it comes to unifying principles, no one is likely to have
  anticipated him. Except, maybe, Peirce.'' \cite[Note 4,
  p.\,85]{MR0419171}
\end{quote}

Condition iv \cite[p.\ 206]{MR0409114} shows that implication should
be interpreted as residuation, defined as an operation on binary
relations by \eqref{comp}, \eqref{conv}, and
\begin{equation}\label{resid}
  \A\to\B=\min{\conv\A|\min\B},
\end{equation}
or, in expanded form
\begin{equation}\label{resid2}
  (\A\to\B)(\x,\y)\iff\forall\z(\A(\z,\x)\to\B(\z,\y)).
\end{equation}
A good example of residuation is the subset relation between sets---it
is the residual of the membership relation with itself.  Formulas of
relevance logic may be interpreted as subsets of a relevant model
structure, \ie, as sets of atoms in the atom structure of a relation
algebra, \ie, as elements of an atomic relation algebra (since the
elements are joins of sets of atoms), or, and this is the most
important case, as binary relations in a proper relation algebra.
This includes an interpretation each connective in any relation
algebra, and in proper relation algebras those interpretations are
disjunction as union, conjunction as intersection,
negation as converse-complementation, and implication as residuation.

What remains is to figure out, from the Routley-Meyer definition of
verification in a relevant model structure, how a formula is verified
in a proper relation algebra.  Routley and Meyer explain,
\begin{quote}
  ``The \emph{real world} plays a distinguished role in our semantical
  postulates.  (Accordingly we call it $0$ rather than $G$; not only
  does the former look better [this is supposed to be, remember, a
  \emph{mathematical} semantics] , but it correctly hints that $0$
  will play the formal role of an identity.) It's necessary to
  distinguish $0$ for the following reason: Logical truth does
  \emph{not} turn out to be \emph{truth} in all set-ups; for the
  strategy which dispatches the paradoxes lies in allowing even
  logical identities to turn out sometimes false.  (What, after all,
  could be better grounds for denying that $q$ entails $p\to p$ than
  to admit that sometimes $q$ is true when, essentially on grounds of
  relevance, $p\to p$ isn't?)

  ``What then is logical truth? Truth in all set-ups, of course,
  \emph{in which all the logical truths are true}!''  \cite[p.\
  202]{MR0409114}

  ``Truth at $0$ is as noted earlier what counts in verifying logical
  truths; accordingly we say simply that $A$ is \emph{verified} on
  $v$, or on the associated $I$, just in case $I(A,0) = T$, and
  otherwise that $A$ is \emph{falsified} on $v$.''  \cite[p.\
  207]{MR0409114}
\end{quote}
In other words, if the map determined by an interpretation sends a
formula $\A$ to a set of set-ups that includes $0$, then that formula
is verified. The distinguished element in the atom structure of an
integral relation algebra is the identity element. Integral relation
algebras are exactly the ones in which the identity element is an
atom. The identity element matches up with the distinguished $0$ of a
relevant model structure. In a proper relation algebra, the identity
element is the identity relation on the underlying set whose pairs
make up the binary relations belonging to the proper relation algebra.

Assign the binary relation symbols of Tarski's extended first-order
logic to binary relations in a proper relation algebra.  According to
the Routley-Meyer definition, a formula is verified under this
assignment if and only if it evaluates (under the interpretation of
its connectives as operations on binary relations) to a relation that
contains the identity relation.  Therefore a formula $\A$ is verified
in a proper relation algebra if and only if
\begin{equation*}	\forall\x(\A(\x,\x))
\end{equation*}
is true under this assignment. What does this mean for an implication?
An implication $\A\to\B$ is verified if and only if
\begin{equation*}
\forall\x((\A\to\B)(\x,\x)),
\end{equation*}
or the equivalent sentence
\begin{equation*}
\forall\x\forall\y(\A(\x,\y)\to\B(\x,\y)),
\end{equation*}
is true. These sentences assert that the relation denoted by $\A$ is
included in the relation denoted by $\B$. The verified implications
are the inclusions between binary relations obtained by interpreting
the conectives as operations on binary relations.

All this is standard operating procedure in the theory of relation
algebras.  Ever since Tarski's and Lyndon's work in the 1950s, it has
been a relevant question to ask for every relation algebra, is it
isomorphic to a proper relation algebra (\ie, representable)?  And if
it is, what does that say about the binary relations in it?

Theorems asserting that relation algebras are representable are among
the most important parts of the subject. Tarski's early QRA Theorem is
a prime example. If Tarski's axioms for relation algebras had turned
out to be complete, then his long and difficult theorem would have
become pointless. Tarski's QRA Theorem (see \cite[VII]{Tarski1953} or
\cite[8.4(iii)]{MR920815} or \cite[Theorem~427]{Maddux2006}) asserts
that if a relation algebra contains a pair of quasi-projections
(elements that behave like projection functions) then it is
representable. The QRA Theorem follows from the main result of the
Tarski-Givant book, called the Main Mapping Theorem for $\L^\times$
and $\L_\n^+$ \cite[4.4(xxxiii)(xxxiv)]{MR920815},
\cite[Theorem~574]{Maddux2006}.  The Main Mapping Theorem says that if
a theory, formalized in first-order logic, proves the existence of a
pair of functions acting sufficiently like projection functions (from
ordered pairs to their components), then that theory can be formalized
as a equational theory in the language of relation algebras.  This
enables Tarski's formalization of set theory without variables
(\cite{Tarski1953a}, \cite[\S4.6]{MR920815}).

To recall the characterization of Tarski's axioms, let $\EL$ be the
equations provable in Tarski's extended system of first-order logic of
equality and other binary relations restricted to four variables.
This class of equations is axiomatized by Tarski's axioms for relation
algebras together with the rules of deduction for equational logic.
These equations contain the entire range of operations used by the
nineteenth century algebraic logicians: union, intersection,
complementation, converse, relative multiplication, and a
distinguished identity relation. We might call $\EL$ ``Tarski's
equational logic'' (for relation algebras).

Applying this characterization with the reduced set of operations
available in relevance logic produces Tarski's relevance logic $\FL$.
By definition, $\FL$ consists of those formulas for which
$\forall\x(\A(\x,\x))$ is provable in first-order logic of binary
relations restricted to four variables.  Unlike $\EL$, the formulas in
$\FL$ contain only the operations corresponding to the connectives of
relevance logic: union, intersection, converse-complementation, and
residuation.  Note that relative multiplication and residuation can be
defined from each other using converse-complementation. On the other
hand, complementation and converse (the Routley star) do \emph{not}
occur in the formulas in $\FL$.

This definition of $\FL$ is precise enough to demonstrate what
formulas are in $\FL$, what derived rules of inference it is closed
under, and what formulas are not in $\FL$. The exact choice of logical
axioms for first-order logic doesn't really matter, as experience has
shown. One can quibble about what 4-variable logic should be. For
example, respelling of bound variables is usually presented as a
consequence of the logical axioms, but its proof requires extra
variables not occurring in a given sentence, and these may not exist
if all four variables already occur in a sentence. Respelling of bound
variables can be excluded or explicitly included, but the result is
the same. For the sake of avoiding such questions and the notational
complexities of quantifiers, a sequent calculus was employed in
\cite{MR722170}, as will be done here. (Another good option are proofs
by natural deduction, restricted to examination of at most four
objects at once.) Some of the rules from \cite{MR722170} can be used
directly (the structural rules and ones for $\land$ and $\lor$), while
new rules are formulated for the connectives $\rmin$ and $\to$.  The
resulting proofs are close in appearance to informing reasoning using
at most four objects.

By \cite[Theorems~2]{MR722170} for $\n=4$ (or $\n=3$), a formula is
provable with four (or three) variables, using the complete set of
rules in \cite{MR722170}, if and only if the corresponding equation is
true in all relation algebras by \cite[Theorems~5]{MR722170} (or
semi-associative relation algebras by \cite[Theorems~4]{MR722170}).
The rules used here are a proper subset of the rules in
\cite{MR722170}, or are the result of the combined application of two
rules from \cite{MR722170}, as is the case for ${\to}|$, $|{\to}$,
$\rmin|$, and $|\rmin$.  Consequently every formula in $\FL$ (or
$\L_3$) corrsponds to an equation true in all relation algebras (or
semi-associative relation algebras). The correspondence is quite
direct in \cite{MR722170}. An inclusion $\A\subseteq\B$ is true in all
relation algebras if and only if the sequent
$\fm\A\0\1\implies\fm\B\0\1$ is provable in 4-variable logic.  The
equivalent condition here is that $\implies\fm{\A\to\B}\0\0$ is
provable in 4-variable logic. These two sequents are interderivable,
corresponding to the fact that one relation is a subset of another if
and only if their residual contains the identity relation:
$\A\subseteq\B$ iff $\Id\subseteq\A\to\B$.
\begin{figure}
  \begin{align*}
    &\begin{aligned}&\\\hline\stand
      \Gamma,\fm\A\i\j&\implies\Delta,\fm\A\i\j
    \end{aligned}\text{\quad\boxed{\text{Axiom}}}\\\\
    &\begin{aligned}
      \Gamma&\implies\Delta,\fm\A\i\j\\
      \fm\A\i\j,\Gamma'&\implies\Delta'\\\hline\stand
      \Gamma,\Gamma'&\implies\Delta,\Delta'
    \end{aligned}\text{\quad\boxed{\text{Cut}}}
    &&\begin{aligned}
      \Gamma&\implies\Delta\\
      \hline\stand\Gamma,\Gamma'&\implies\Delta,\Delta'
      \end{aligned}\text{\quad \boxed{\text{Weakening}}}\\\\
      &\begin{aligned}
        \Gamma,\fm\A\i\j&\implies\Delta\\
        \Gamma',\fm\B\i\j&\implies\Delta'\\
        \hline\stand\Gamma,\Gamma',\fm{\A\lor\B}\i\j&\implies\Delta,
        \Delta'\end{aligned}\quad\boxed{\lor|} &&\begin{aligned}
        \Gamma&\implies\Delta,\fm\A\i\j,\fm\B\i\j\\\hline\stand
        \Gamma&\implies\Delta,\fm{\A\lor\B}\i\j
      \end{aligned}\quad\boxed{|\lor}
      \\\\&\begin{aligned}
        \Gamma,\A_{\i\j},\B_{\i\j}&\implies\Delta\\\hline\stand
        \Gamma,(\A\land\B)_{\i\j}&\implies\Delta
      \end{aligned}\quad\boxed{\land|}
      &&\begin{aligned}
        \Gamma&\implies\Delta,\A_{\i\j}\\
        \Gamma'&\implies\Delta',\B_{\i\j}\\\hline\stand
        \Gamma,\Gamma'&\implies\Delta,\Delta',(\A\land\B)_{\i\j}
      \end{aligned}\quad\boxed{|\land}
      \\\\&\begin{aligned}
        \Gamma&\implies\Delta,\fm\A\i\j\\
        \hline\stand\Gamma,\fm{\rmin\A}\j\i&\implies\Delta
      \end{aligned}\quad\boxed{\rmin|}
      &&\begin{aligned} \Gamma,\fm\A\i\j&\implies\Delta\\\hline\stand
        \Gamma&\implies\Delta,\fm{\rmin\A}\j\i
      \end{aligned}\quad\boxed{|\rmin}
      \\\\&\begin{aligned}
        \Gamma&\implies\Delta,\fm\A\k\i\\
        \Gamma',\fm\B\k\j&\implies\Delta'\\\hline\stand
        \Gamma,\Gamma',\fm{A\to\B}\i\j&\implies\Delta,\Delta'
      \end{aligned}\quad\boxed{{\to}|}
      &&\begin{aligned}
        \Gamma,\fm\A\k\i&\implies\Delta,\fm\B\k\j\\\hline\stand
        \Gamma&\implies\Delta,\fm{\A\to\B}\i\j\text{,\,\,no $\k$}
      \end{aligned}\quad\boxed{|{\to}}
    \end{align*}
    \caption{Rules and Axioms for Tarski's relevance logic}
    \label{fig1}
  \end{figure}
\section{The sequent calculus}
\begin{definition}\label{def1}\ 
  \begin{itemize}
  \item $\var$ is a countable set whose elements
    $\p,\q,\r,\dots\in\var$ are called {\bf propositional variables}.
  \item $\fmla$ is the closure of $\var$ under three binary operations
    $\lor$, $\land$, and $\to$, and one unary operation $\rmin$.
  \item The elements $A,\B,\C,\D,\dots\in\fmla$ are called {\bf
      formulas}.
  \item The four elements of $\four$ are called {\bf objects} or {\bf
      indices}.
  \item An {\bf assertion} $\fm\A\i\j$ is a formula $\A$ together with
    an ordered pair of individual objects $\i,\j\in\four$, added to
    the formula as subscripts.
  \end{itemize}
\end{definition}
Parentheses are omitted according to the convention that the
operations are applied in this order: $\rmin$, $\land$, $\lor$, and
finally $\to$.  An assertion $\fm\A\i\j$ should be read as if it said
$\<\i,\j\>\in\A$, that is, $\A$ is a relation that holds between
objects $\i$ and $\j$. In first-order logic an assertion might more
commonly be written $\A(\i,\j)$, as was done earlier. The subscript
style of writing an assertion was common in nineteenth century
algebraic logic, and it reduces the burden of parentheses.
\begin{definition}\label{def2}\ 
  \begin{itemize}
  \item A {\bf sequent} $\Gamma\implies\Delta$ is an ordered pair
    $\<\Gamma,\Delta\>$ of sets of assertions $\Gamma$ and $\Delta$.
  \item The sequent $\Gamma\implies\Delta$ is an {\bf Axiom} if
    $\Gamma\cap\Delta\neq\emptyset$.
  \item A {\bf 4-proof} is a finite sequence of sequents in which
    every sequent is either an Axiom or follows from one or two
    previous sequents by one of the rules of inference shown in
    Figure~\ref{fig1}: {\rm Cut}, {\rm Weakening}, {$\lor|$},
    {$|\lor$}, {$\land|$}, {$|\land$}, {$\rmin|$}, {$|\rmin$},
    {${\to}|$}, and, if $k$ does not appear in $\Gamma\cup\Delta$,
    {$|{\to}$}.
  \end{itemize}
\end{definition}
The restriction to finite proofs in Definition~\ref{def2} is motivated
by the fact that a rule can ``do nothing''. For example, every sequent
follows from itself by Weakening (take $\Gamma'=\Delta'=\emptyset$),
Without the restriction, the infinite $\mathbb\Z$-indexed sequence in
which every sequent is $\fm\A\0\1\implies\fm\B\0\1$ would be a
``proof'' of $\fm\A\0\1\implies\fm\B\0\1$.  Abbreviations used in the
notation for sequents is standard. For example,
\begin{equation*}	\Delta,\Gamma,\A\implies\Delta',\Gamma',\B,\C
\end{equation*}
is short for
\begin{equation*}	\Delta\cup\Gamma\cup\{\A\}\implies
\Delta'\cup\Gamma'\cup\{\B,\C\}.
\end{equation*}
A sequent $\Gamma\implies\Delta$ should be read, ``If all the
assertions in $\Gamma$ are true, then one of the assertions in
$\Delta$ is true.'' For example, the sequent
$\fm\A\i\j\implies\fm\B\i\j$ should be read, ``If $\<\i,\j\>\in\A$
then $\<\i,\j\>\in\B$''. Under this reading, together with the
intended interpretation of the connectives as set-theoretical
operations, it is easy to see why all the rules in Figure~\ref{fig1}
are correct. In particular, the rule $|{\to}$ requires that
$\k\neq\i,\j$ and $\k$ does not occur as a subscript in any assertion
in $\Gamma$ or $\Delta$, as indicated by the notation ``no $\k$''.
The reason for this is the universal quantifier in the definition of
residuation, and is reflected in one of the common logical validities
used in axiomatizations of first-order logic, namely
$\forall\x(\varphi\to\psi)\to(\varphi\to\forall\x\psi)$, where it is
required that $\x$ does not occur free in $\varphi$.  In proofs that a
formula belongs to $\FL$, the notation ``no $\k$'' will accompany
every application of rule {$|{\to}$}, explicitly identifying the
universally quantified object.
\begin{definition}\label{def3}\ 
\begin{itemize}
\item A {\bf 4-proof of the sequent} $\Gamma\implies\Delta$ is a
  4-proof in which $\Gamma\implies\Delta$ appears.  We write
  \begin{equation*}
    \provesFL\Gamma\implies\Delta
  \end{equation*}
  just in case $\Gamma\implies\Delta$ has a 4-proof.
\item A {\bf 4-proof of the formula} $\A$ is a 4-proof of the sequent
  $\implies\fm\A\0\0$.
\item $\FL$ is the set of formulas that have 4-proofs:
  \begin{equation*}
    \FL=\{\A:\,\,\provesFL\,\,\implies\fm\A\0\0\}.
  \end{equation*}
\end{itemize}
\end{definition}
\begin{table}
\begin{align*}
  \text{Lemma}&			&&\text{Objects}&\text{Formula}\\
  \lref{t6}&			&&\{\0\}    	&\A&\lor\rmin\A\\
  \lref{A1.}&			&&\{\0,\1\}		&\A&\to\A\\
  \lref{A2.}&			&&\{\0,\1\}		&\A\land\B&\to\A\\
  \lref{A3.}&			&&\{\0,\1\}		&\A\land\B&\to\B\\
  \lref{A5.}&			&&\{\0,\1\}		&\A&\to\A\lor\B\\
  \lref{A6.}&			&&\{\0,\1\}		&\B&\to\A\lor\B\\
  \lref{comm1}&		&&\{\0,\1\}		&\B\lor\A&\to\A\lor\B\\
  \lref{comm2}&		&&\{\0,\1\}		&\B\land\A&\to\A\land\B\\
  \lref{assoc1}&		&&\{\0,\1\}		&(\A\land\B)\land\C&\to\A\land(\B\land\C)\\
  \lref{assoc2}&		&&\{\0,\1\}		&(\A\lor\B)\lor\C&\to\A\lor(\B\lor\C)\\
  \lref{A8.}&			&&\{\0,\1\}		&\A\land(\B\lor\C)&\to(\A\land\B)\lor(\A\land\C)\\
  \lref{*T9.}&			&&\{\0,\1\}		&(\A\to\rmin\C)\land(\B\to\C)&\to\rmin(\A\land\B)\\
  \lref{T10.}&			&&\{\0,\1\}		&(\A\to\rmin\B)\land(\rmin\A\to\rmin\C)&\to\rmin\B\lor\rmin\C\\
  \lref{A9.}&			&&\{\0,\1\}		&\rmin\rmin\A&\to\A\\
  \lref{t3}&			&&\{\0,\1\}		&\A&\to\rmin\rmin\A\\
  \lref{T2.}&			&&\{\0,\1\}		&\rmin(\A\lor\B)&\to\rmin\A\land\rmin\B\\
  \lref{t4}&			&&\{\0,\1\}		&\rmin(\A\land\B)&\to\rmin\A\lor\rmin\B\\
  \lref{t5}&			&&\{\0,\1\}		&\rmin\A\land\rmin\B&\to\rmin(\A\lor\B)\\
  \lref{t5a}&			&&\{\0,\1\}		&\rmin\A\lor\rmin\B&\to\rmin(\A\land\B)\\
  \lref{T11.}&			&&\{\0,\1\}		&((\A\to\A)\to\B)&\to\B\\
\end{align*}\caption{Formulas in Tarski's relevance logic, provable with 1 or 2 objects}\label{table1}
\end{table}
\begin{table}
\begin{align*}
  \text{Lemma}&		&&\text{Objects}&\text{Formula}\\
  \lref{A4.}&		&&\{\0,\1,\2\}	&(\A\to\B)\land(\A\to\C)&\to(\A\to\B\land\C)\\
  \lref{A7.}&		&&\{\0,\1,\2\}	&(\A\to\C)\land(\B\to\C)&\to(\A\lor\B\to\C)\\
  \lref{t11}&		&&\{\0,\1,\2\}	&(\A\to\B)\land(\C\to\D)&\to(\A\land\C\to\B\land\D)\\
  \lref{T6.}&		&&\{\0,\1,\2\}	&(\A\to\B)\land(\C\to\D)&\to(\A\lor\C\to\B\lor\D)\\
  \lref{T8.}&		&&\{\0,\1,\2\}	&(\A\to\B)\lor(\C\to\D)&\to(\A\land\C\to\B\lor\D)\\
  \lref{t7}&		&&\{\0,\1,\2\}	&\A&\to(\rmin\B\to\rmin(\A\to\B))\\
  \lref{t9}&		&&\{\0,\1,\2\}	&\A&\to(\B\to\rmin(\A\to\rmin\B))\\
  &			&&		&\A&\to(\B\to\A\circ\B)\\
  \lref{t13}&		&&\{\0,\1,\2\}	&\A&\to((\rmin\B\to\rmin\A)\to\B)\\
  \lref{t14}&		&&\{\0,\1,\2\}	&\A&\to((\B\to\rmin\A)\to\rmin\B)\\
  \lref{T12.}&		&&\{\0,\1,\2\}	&\rmin((\A\to\B)\to\rmin\A)&\to\B\\
  &			&&		&(\A\to\B)\circ\A&\to\B\\
  \lref{*T15.}&	&&\{\0,\1,\2\}	&\rmin\A&\to((\B\to\A)\to\rmin\B))\\
  \lref{reflection}&	&&\{\0,\1,\2\}	&\rmin(\A\to\rmin\B)\land\C&\to\rmin((\A\land\rmin\D)\to\rmin\B)\\
  &			&&		&&\quad{}\lor\rmin(\A\to\rmin(\B\land\rmin(\D\to\rmin\C)))\\
  &			&&		&(\A\circ\B)\land\C&\to((\A\land\rmin\D)\circ\B)\lor(\A\circ(\B\land(\D\circ\C)))\\
  \lref{t??}&		&&\{\0,\1,\2\}	&(\A\to\B)\land\rmin\(\C\to\rmin\D\)&\to\rmin(\C\land\B\to\rmin\D)\\
  &			&&		&&\quad{}\lor\rmin(\C\to\rmin(\D\land\rmin\A))\\
  &			&&		&(\A\to\B)\land\(\C\circ\D\)&\to ((\C\land\B)\circ\D)\lor(\C\circ(\D\land\rmin\A))\\
  \lref{prefixingA}&	&&\{\0,\1,\2,\3\}&(\A\to\B)&\to((\C\to\A)\to(\C\to\B))\\
  \lref{t10}&		&&\{\0,\1,\2,\3\}&(\B\to(\C\to\A))&\to(\rmin(\B\to\rmin\C)\to\A)\\
  &			&&		&(\B\to(\C\to\A))&\to((\B\circ\C)\to\A)\\
  \lref{T19.}&		&&\{\0,\1,\2,\3\}&(\rmin(\A\to\rmin\B)\to\C)&\to(\A\to(\B\to\C))\\
  &			&&		&((\A\circ\B)\to\C)&\to(\A\to(\B\to\C))\\
  \lref{t?}&		&&\{\0,\1,\2,\3\}&(\A\to\B)&\to(\rmin(\A\to\C)\to\rmin(\B\to\C))\\
  &			&&		&(\A\to\B)&\to((\A\circ\D)\to(\B\circ\D))\\
  \lref{assocfusion}& &&\{\0,\1,\2,\3\}
  &(\A\circ\B)\circ\C&\to\A\circ(\B\circ\C)
\end{align*}\caption{Formulas in Tarski's relevance logic, provable with 3 or 4 objects}\label{table2}
\end{table}
\begin{table}
  \begin{align*}
  \text{Lemma}			&&&\text{Objects}&\text{Rule}\\
    \lref{adjunction}		&&&\{\0\}		&\A,\,\,\B&\proves\A\land\B\\
    \lref{modusponens}		&&&\{\0\}		&\A\to\B,\,\A&\proves\B\\
    \lref{disjunctivesyllogism}	&&&\{\0\}		&\A\lor\B,\,\,\rmin\A&\proves\B\\
    \lref{transitivity}		&&&\{\0,\1\}		&\A\to\B,\,\,\B\to\C&\proves\A\to\C\\
    \lref{contraposition}	&&&\{\0,\1\}		&\A\to\B&\proves\rmin\B\to\rmin\A\\
    \lref{contraposition.2}	&&&\{\0,\1\}		&\A\to\rmin\B&\proves\B\to\rmin\A\\
    \lref{cut}			&&&\{\0,\1\}		&\A\land\B\to\C,\,\,\B\to\C\lor\A&\proves\B\to\C\\
    \lref{E-rule}		&&&\{\0,\1\}		&\A&\proves(\A\to\B)\to\B\\
    \lref{suffixing}		&&&\{\0,\1,\2\}		&\A\to\B&\proves(\B\to\C)\to(\A\to\C)\\
    \lref{cycling}		&&&\{\0,\1,\2\}		&\A\to(\B\to\C)&\proves\B\to(\rmin\C\to\rmin\A)\\
    \lref{prefixingR}		&&&\{\0,\1,\2,\3\}	&\A\to\B&\proves(\C\to\A)\to(\C\to\B)\\
    \lref{R3}			&&&\{\0,\1,\2,\3\}	&\A\to\B,\,\,\C\to\D&\proves(\B\to\C)\to(\A\to\D)\\
    \lref{monotonicfusion}	&&&\{\0,\1,\2,\3\}	&\A\to\B,\,\,\C\to\D&\proves(\A\circ\C)\to(\B\circ\D)\\
    & &&
    &\A\to\B,\,\,\C\to\D&\proves\rmin(\A\to\rmin\C)\to\rmin(\B\to\rmin\D)
  \end{align*}\caption{Some derived rules of inference in Tarski's relvance logic}\label{table3}
\end{table}
\section{Tables of formulas and rules}
Tables~\ref{table1} and ~\ref{table2} show more than three dozen
formulas in $\FL$.  Each entry begins with the number in parentheses,
preceded by ``L'', of the lemma in which that formula is shown to have
a 4-proof. For example, the proof of Lemma~\ref{t6} is a 4-proof of
formula \lref{t6}.  In a 4-proof, every sequent is either an Axiom or
follows from the one or two sequents immediately preceding it
according to the rule mentioned to the right. Line numbers in 4-proofs
are included whenever the justifying sequents are not the previous one
or two. The second entry in Tables~\ref{table1} and~\ref{table2} is a
list of the objects that are actually used in the 4-proof of the
formula. This provides a rough classification of the formulas into
those belonging to what we might call $\L_1$, $\L_2$, $\L_3$, and
$\FL$, depending on the number of objects needed for their 4-proofs.

All of the formulas in $\FL$ make assertions about binary relations
that are universally true. As was observed earlier, the verification
of a formula of the form $\A\to\B$ in every proper relation algebra
confirms that $\A\subseteq\B$, no matter how the propositional
variables in $\A$ and $\B$ are interpreted as binary relations.  For
example, formula \lref{A1.} asserts the universal truth that for
every binary relation $\A$, $\A\subseteq\A$, while \lref{A2.} asserts
that for all binary relations $\A$ and $\B$, $\A\cap\B\subseteq\A$, as
one would expect if the interpretation of $\land$ is intersection.

Table~\ref{table3} shows more than a dozen derived rules of inference
in $\FL$.  Each entry begins with ``L'' and the number in parentheses
of the lemma in which the rule is shown to have a 4-proof under the
assumption that the inputs to the rule have 4-proofs.  For example,
the proof of Lemma~\ref{adjunction} shows how to assemble 4-proofs of
$\A$ and $\B$ into a 4-proof of $\A\land\B$.  The second entry in
Table~\ref{table3} is a list of the objects needed for this assembly.
Of course the 4-proofs of $\A$ and $\B$ may use all available objects,
but if not, then the second entry in the table shows what additional
objects, if any, might be required.  Again, this provides a rough
classification of the rules into those belonging to $\L_1$, $\L_2$,
$\L_3$, and $\FL$.

\section{Comparisons with other systems}
Tables~\ref{table1}, \ref{table2}, and \ref{table3} help locate
Tarski's relevance logic in the pantheon of relevance logics. There is
a very great contrast here between Tarski's relevance logic and the
usual world of relevance logics.  Indeed, the situation is well
expressed by the following quotations.
\begin{quote}
  ``Old friends of our project will be surprised to find that we were
  forced to split the book into two volumes -- in order, of course, to
  avoid weighing the reader down either literally or financially --
  when we finally realized that the universe of relevance logics had
  expanded unnoticed overnight.'' \cite[p.\,xxiii]{MR0406756}
\end{quote}
\begin{quote}
  ``This book mentions or discusses so many different systems (Meyer
  claims the count exceeds that of the number of ships in Iliad II)
  that we have been driven \dots\ to try to devise a reasonably
  rational nomenclature.'' \cite[p.\,xxv]{MR0406756}
\end{quote}
\begin{quote}
  ``Additional axiom schemes drawn from the following lists may be
  added to basic system $\sf\B$ \dots\ singly or in combination to
  yield a wealth of stronger systems:--''
  \cite[p.\,288]{Routleyetal1982}
\end{quote}
\begin{quote}
  ``The following postulates are added \dots, singly or in
  combination, to provide modellings for the wealth of further systems
  of sentential logics introduced in the previous section.''
  \cite[p.\,300]{Routleyetal1982}
\end{quote}
\begin{quote}
  ``In this chapter we first present algebraic analyses for an
  important and extensive class of affixing systems: the class
  comprises not only a great many relevant logics including all the
  more standard systems but also all the usual irrelevant logics and
  some unusual ones as well'' \cite[p.\,72]{Brady2003}
\end{quote}
Tarski's relevance logic $\FL$ does not have this sort of variation.
There is no list of formulas and rules from which to choose ``singly
or in combination''.  The only available parameter is the number of
variables used to prove any particular formula or deductive rule
expressing a property of binary relations.  The most interesting cases
are when the number of variables is 1, 2, 3, or 4.  The logic $\L_1$
already has the Law of the Excluded Middle, and among its rules are
Adjunction, {\it modus ponens}, and Disjunctive Syllogism.  The logic
$\L_2$ picks up all the formulas in Table~\ref{table1} (many of which
are part of various systems of Basic Logic), plus some more rules from
Table~\ref{table3}, such as the Rules of Transitivity, Contraposition,
and Cut.  The logic $\L_3$ adds to this list the Rule of Suffixing,
for example, along with some key formulas governing conjunction,
disjunction, and fusion.  However, the associative law for fusion is
missing from $\L_3$, along with those axioms (such as Suffixing) and
rules (such as Prefixing) from Tables~\ref{table3} and~\ref{table4}
whose sequent proofs require four objects. (These omissions can be
proved by examining semi-associative relation algebras that are not
associative, hence not relation algebras, which fail to satisfy the
the appropriate rules and equations).  The logic $\L_5$, however, is
(or, at least, has been) well beyond the consideration of even the
most ardent inventors of systems.

Even $\FL$ misses standard axioms used in various relevance logics.
Such axioms can be added, perhaps yielding a ``wealth of systems''.
$\FL$ is a ``naturally occurring'' system. The motivation for studying
$\FL$ certainly involves relevance logics. But $\FL$ arises from
entirely different considerations.  Indeed, Tarski's relevance logic
$\FL$ may satisfy van Benthem's~\cite{vanBenthem1984} suggestion that
\begin{quote}
  ``\dots, the Routley semantics still has to prove its mettle. On the
  realistic side, its model structures ought to admit of, if not a
  natural linguistic anchoring, then at least one mathematical
  `standard example', providing some food for independent
  reflection.''
\end{quote}
Perhaps Tarski's relevance logic should be considered as a
``standard mathematical example.''
\section{Basic logic}
Tarski's relevance logic contains the Basic Logic $\sf\B$ of
\cite{Brady2003} and \cite{Routleyetal1982}.  The axioms of Basic
Logic in \cite[pp.\,287--8]{Routleyetal1982} are A1--A9, and its rules
are R1--R5, with R3$'$ as an alternative to rules R3 and R4.  In
Tables~\ref{table1} and~\ref{table2}, A1 is \lref{A1.}, A2 is
\lref{A2.}, A3 is \lref{A3.}, A4 is \lref{A4.}, A5 is \lref{A5.}, A6
is \lref{A6.}, A7 is \lref{A7.}, A8 is \lref{A8.}, and A9 is
\lref{A9.}.  The rules of Basic Logic in
\cite[pp.\,287--8]{Routleyetal1982} are derived rules of inference in
Tarski's relevance logic.  In Table~\ref{table3}, R1 is
\lref{modusponens}, R2 is \lref{adjunction}, R3$'$ is \lref{R3}, R3 is
\lref{suffixing}, R4 is \lref{prefixingR}, and R5 is
\lref{contraposition.2}.  The axioms of Basic Logic in
\cite[pp.\,192--3]{Brady2003} are A1--A9, the same as axioms A1--A9 of
\cite[pp.\,287--8]{Routleyetal1982}.  The rules of Basic Logic in
\cite[p.\,193]{Brady2003} are R1--R4, where R1 is \lref{modusponens},
R2 is \lref{adjunction}, R3 is \lref{R3}, and R4 is
\lref{contraposition.2}.  Rule R5 in \cite[p.\,192--3]{Brady2003} is
part of systems $\sf\E$ and $\sf\E\W$; R5 is \lref{E-rule}.  Axiom A13
of system $\sf\T\W$ in \cite[p.\,193]{Brady2003} is \lref{prefixingA}.
Axiom A17 of systems $\sf\D\K$ and $\sf\T\K$ in
\cite[p.\,193]{Brady2003} is \lref{t6}.
\section{Properties of binary relations}
All the formulas and rules of inference in Tarki's relevance logic are
true for arbitrary binary relations. They are verified in all proper
relation algebras. More generally, they hold in every algebra of the
form 
\begin{align*}
  \Ka&=\<\K,\cup,\cap,\to,\rmin\>,
\end{align*} 
where $\K$ is a set of binary relations on some set $\U$, and $\K$ is
closed under union $\cup$, intersection $\cap$, residuation $\to$, and
converse-complementation $\rmin$. This means that if the propositional
variables in a formula $\A$ in Tarski's relevance logic are assigned
to binary relations in $\K$, then the binary relation assigned to $\A$
will contain the identity relation on $\U$. Conversely, any formula
that holds in every such algebra $\Ka$ will be part of Tarski's
relevance logic if it can be proved by looking at no more than four
objects at a time.  Thus, formulas not belonging to Tarski's relevance
logic are of two kinds.  They are either valid for all binary
relations but require more than five objects to prove, or else they
postulate properties of binary relations that do not hold in general.

Here are some examples of formulas expressing special properties of
binary relations; for details see \cite{KramerMaddux2019}. The axiom
of contraposition,
\begin{equation}
  \label{contra}
  (\A\to\rmin\B)\to(\B\to\rmin\A)
\end{equation}
holds in $\Ka$ if and only if the relations in $\K$ commute with each
other under relative multiplication, \ie, fusion is commutative.  The
same applies to the axiom of permutation,
\begin{equation}
  \label{perm}
  (\A\to(\B\to\C))\to(\B\to(\A\to\C)).
\end{equation}
Commutativity of $\Ka$ is enough to insure that the suffixing and
modus ponens axioms
\begin{align}
  \label{suff}
  (\A\to\B)&\to((\B\to\C)\to(\A\to\C))\\
  \label{mp}
  \A&\to((\A\to\B)\to\B)
\end{align}
hold in $\Ka$, but neither of them is equivalent to assuming $\Ka$ is
commutative.  The contraction axiom and the reductio axiom
\begin{align}
  \label{contr}
  (A\to(\A\to\B))&\to(\A\to\B)\\
  \label{reduc}
  (\A\to\rmin\A)&\to\rmin\A
\end{align}
are each equivalent to assuming every relation in $\K$ is dense.  The
$\RR$-mingle axiom
\begin{equation}\label{ming}
  \A\to(\A\to\A)
\end{equation}
holds in $\Ka$ if and only if every relation in $\K$ is transitive. 
\section{System $\RR$}
An axiom set for the Anderson-Belnap system $\RR$ of relevant
implication is presented by Routley-Meyer \cite[p.\,204]{MR0409114}.
It contains axioms A1--A13 along with axioms A14 and A15
\cite[p.\,224]{MR0409114} when fusion $\circ$ is included as primitive
rather than defined as in \eqref{D1}.  Eleven of these fifteen axioms
occur in Tarski's relevance logic.  In Tables~\ref{table1}
and~\ref{table2}, A1 is \lref{A1.}, A5 is \lref{A2.}, A6 is
\lref{A3.}, A7 is \lref{A4.}, A8 is \lref{A5.}, A9 is \lref{A6.}, A10
is \lref{A7.}, A11 is \lref{A8.}, A13 is \lref{A9.}, A14 is \lref{t9},
and A15 is \lref{t10}.  The remaining four axioms of $\RR$ (A2, A3,
A4, A12) do not occur in Tarski's relevance logic: A2 is \eqref{mp},
A3 is \eqref{suff}, A4 is \eqref{contr}, and A12 is \eqref{contra}.
The rules for $\RR$ are \lref{adjunction} and \lref{modusponens}, both
part of Tarski's relevance logic.  If Tarski's relevance logic is
extended by adding axioms A1--A15, then all formulas of the logic
$\RR$ of Anderson-Belnap~\cite{MR0406756} become provable; see
\cite[Corollary 5.2(i)]{MR2641636}.

Adding axioms to $\FL$ may be done by supplementing the rules in
Figure~\ref{fig1}. For example, to add \eqref{contra}, one may include
the rule
\begin{equation*}
\begin{aligned}&\\\hline\stand
  \Gamma,\fm{\A\to\rmin\B}\i\j&\implies\Delta,\fm{\B\to\rmin\A}\i\j
\end{aligned}\text{\quad\boxed{\text{Contraposition}}}\\\\
\end{equation*}
\section{System $\RR$-mingle}
If Tarski's relevance logic is extended by adding the axioms
\eqref{contra}, \eqref{suff}, \eqref{mp}, \eqref{contr}, and
\eqref{ming}, the result is the Dunn-McCall system $\RR$-mingle.
$\RR$-mingle contains every formula valid for transitive dense
commutative binary relations, no matter how many objects are needed
for its proof.  This is just a restatement of
\cite[Theorem~6.2]{MR2641636}.  (See \cite{KramerMaddux2019} for
another proof.)  In more detail, $\A\to\B$ is a theorem of
$\RR$-mingle if and only if the inclusion $\A\subseteq\B$ is true
whenever all its propositional variables are interpreted as relations
in a set $\K$ of dense transitive binary relations, where $\K$ is
closed under union, intersection, residation, and
converse-complementation (the interpretations of the connectives in
$\A\to\B$) and $\K$ is commutative under relative multiplication.  The
underlying reason is that, as Meyer proved \cite[Corollaries 3.1, 3.5,
p.\,413--4]{MR0406756}, the theorems of $\RR$-mingle are the formulas
valid in all Sugihara matrices, and all Sugihara matrices are
representable as sets of transitive dense binary relations,
commutative under relative multiplication \cite{KramerMaddux2019,
  MR2641636}.  An informal mnemonic for this result might be
\begin{equation*}
\text{$\RR$-mingle = $\L_\infty$ + all relations are dense,
  transitive, and commute under $\circ$.}
\end{equation*}
This completeness result involves binary relations and their natural
operations.  $\RR$-mingle \emph{is} the set of laws (expressible with
$\cap$, $\cup$, $\to$, $\rmin$) that hold for all transitive, dense,
commutative binary relations.  Perhaps this provides another standard
mathematical example of a relevance logic, as van Benthem suggested,
although sometimes $\RR$-mingle is not regarded as a true relevance
logic because of Meyer's result \cite[RM84,\,p.\,417]{MR0406756} that
$\RR$-mingle has only the weak variable sharing property that if
$\A\to\B$ is a theorem of $\RR$-mingle then either $\A$ and $\B$ share
a propositional variable or $\rmin\A$ and $\B$ are both theorems of
$\RR$-mingle.
\section{System $\L_5$}
Given Lyndon's initial result and the non-finite axiomatizability
results that followed, starting with Monk's \cite{Monk1964} proof that
the equational theory of representable relation algebras is not
finitely based, it was easy to suspect that formulas must exist that
are valid for all binary relations but are not in Tarski's relevance
logic because they require more than five objects to prove; see
\cite{Maddux2007}, \cite[(Q1), p.\,52]{MR2641636}.  Indeed,
Mikul\'as~\cite{zbMATH05544425} proved such a non-finite
axiomatizability result for relevance logic. The formulas involved are
complicated and not generally considered as potential axioms for
relevance logics. Two formulas that are not theorems of $\RR$ are
given in \cite[Theorem~8.2]{MR2641636}. Here is the shorter one.
(Because of the associativity of $\circ$, one set of parentheses has
been omitted from the final term.)
\begin{gather*}
  \Bigg(\Big( (\A_{34}\circ\A_{23})\land\A_{24}\Big)
  \circ\Big((\A_{12}\circ\A_{01})\land\A_{02}\Big)\Bigg)\land\A_{04}\\
  \to\Bigg(\Big((\A_{34}\circ\A_{23})\land\A_{24}\Big)\circ\Big(\big(\A_{12}\circ[\A_{01}
  \land\rmin\A_{01}]\big)\land\A_{02}\Big)\Bigg)\land\A_{04}\\
  {}\lor\,\,\Bigg(\Big(\big([\A_{34}\land\rmin\A_{34}]\circ\A_{23}\big)\land\A_{24}\Big)
  \circ\Big((\A_{12}\circ\A_{01})\land\A_{02}\Big)\Bigg)\land\A_{04}\\
  {}\lor\,\,\A_{34}\circ\Bigg((\A_{23}\circ\A_{12})\land\Big(\big((\A_{23}\circ\A_{02})
  \land(\A_{43}\circ\A_{04})\big)\circ\A_{10}\Big)\\
  {}\land\Big(\A_{43}\circ\big((\A_{04}\circ\A_{10})\land(\A_{24}\circ\A_{12})\big)\Big)\Bigg)\circ\A_{01}
\end{gather*}
The subscripts on the variables indicate which objects should appear
as subscripts in assertions based on that formula. For example, the
assertion $\fm{\A_{24}}\2\4$ will appear in a properly constructed
5-proof.
\section{Relevant model structures}
Relevant model structures~\cite[\S2]{MR0409114} provide sound and
complete semantics for system $\RR$. We will use them to show various
formulas are not in $\RR$ or not in $\FL$.
\begin{definition}\label{def4}
  A {\bf relevant model structure} $\gc\K=\<\K,\R,{}^*,0\>$ consists
  of a non-empty set $\K$, a ternary relation $\R\subseteq\K^3$, a
  unary operation ${}^*:\K\to\K$, and a distinguished element
  $0\in\K$, such that postulates \thetag{p1}--\thetag{p6} hold for all
  $\a,\b,\c\in\K$, where
\begin{align*}
  \tag{d1}&\R^2\a\b\c\d\ifff\ex\x(\R\a\b\x\text{ and }\R\x\c\d),\\
  \tag{d2}&\R^2\a(\b\c)\d\ifff\ex\x(\R\b\c\x\text{ and }\R\a\x\d).
\end{align*}
\begin{align*}
  \tag{p1}&\R0\a\a&&\text{\rm($0$-reflexivity)}\\
  \tag{p2}&\R\a\a\a&&\text{\rm(density)}\\
  \tag{p3}&\R^2\a\b\c\d\implies\R^2\a\c\d\b\\
  \tag{p4}&\R^20\a\b\c\implies\R\a\b\c&&\text{\rm($0$-cancellation)}\\
  \tag{p5}&\R\a\b\c\implies\R\a\c^*\b^*\\
  \tag{p6}&\a^*{}^*=a &&\text{\rm(involution)}
\end{align*}
\end{definition}
By \cite[Theorem7.1]{MR2641636}, $\gc\K=\<\K,\R,{}^*,0\>$ is a
relevant model structure if and only if it satisfies \thetag{p1},
\thetag{p2}, \thetag{p3$'$}, \thetag{p4}, \thetag{p5$'$}, \thetag{p6},
and \thetag{comm}, where
\begin{align*}
  \tag{comm}&\R\a\b\c\implies\R\b\a\c&&\text{(commutativity)}\\
  \tag{p3$'$}&\R^2\a\b\c\d\implies\R^2\a(\b\c)\d&&\text{(associativity)}\\
  \tag{p5$'$}&\R\a\b\c\implies\R\c^*\a\b^*&&\text{(rotation)}
\end{align*}
\begin{definition}\label{def5}
  Let $\gc\K=\<\K,\R,{}^*,0\>$ be a relevant model structure.  A {\bf
    valuation} in $\gc\K$ is a function
  $\nu\colon\var\times\K\to\{\T,\F\}$ such that, for all $\a,\b\in\K$
  and $\p\in\var$, if $\R0\a\b$ and $\nu(\p,\a)=\T$ then
  $\nu(\p,\b)=\T$. $\I$ is the {\bf interpretation associated with
    $\nu$} if $\I\colon\fmla\times\K\to\{\T,\F\}$, and for all
  $\A,\B\in\fmla$ and $\c\in\K$,
  \begin{enumerate}
  \item\label{i} $\I(\p,\c)=\nu(\p,\c),$
  \item\label{ii} $\I(\A\land\B,\c)=\T$ iff $\I(\A,\c)=\T$ and
    $\I(\B,\c)=\T,$
  \item\label{iii} $\I(\A\lor\B,\c)=\T$ iff $\I(\A,\c)=\T$ or
    $\I(\B,\c)=\T,$
  \item\label{iv} $\I(\A\to\B,\c)=\T$ iff for all $\a,\b$, if
    $\R\c\a\b$ and $\I(\A,\a)=\T$ then $\I(\B,\b)=\T$,
  \item\label{v} $\I(\A\circ\B,\c)=\T$ iff for some $\a,\b$,
    $\R\a\b\c$, $\I(\A,\a)=\T$, and $\I(\B,\b)=\T$,
  \item\label{vi} $\I(\rmin\A,\c)=\T$ iff $\I(\A,\c^*)=\F.$
  \end{enumerate}
A formula $\A$ is {\bf true} on a valuation $\nu$, or on the
associated $\I$, at $\c\in\K$ if $\I(\A,\c)=\T$, and {\bf false} on
$\nu$ at $\c$ if $\I(\A,\c)=\F$.  A formula $\A$ is {\bf verified} on
$\nu$, or on the associated $\I$, if $\I(\A,0)=\T$, otherwise {\bf
  falsified}.  A formula $\A$ is {\bf valid} in $\gc\K$ if $\A$ is
verified on every valuation in $\gc\K$, and {\bf$\RR$-valid} if $\A$ is
valid in every relevant model structure, otherwise {\bf$\RR$-invalid}.
\end{definition}
Condition \ref{v} follows from \ref{iv} and \ref{vi} when definition
\eqref{D1} is used instead of taking $\circ$ as primitive; see
\cite[footnote~10, \p.\,206]{MR0409114}.  By
\cite[Theorem~2]{MR0409114}, all theorems of $\RR$ are $\RR$-valid,
and by \cite[Theorem~3]{MR0409114}, all $\RR$-valid formulas are
theorems of $\RR$.  A relevant model structure
$\gc\K=\<\K,\R,{}^*,0\>$ is {\bf normal} if
$0^*=0$~\cite[p.\,218]{MR0409114}.  By \cite[Theorem~4]{MR0409114}, a
formula $\A$ is a theorem of $\RR$ if and only if $\A$ is valid in
every normal relevant model structure.
\begin{definition}\label{def6}
  Given a relevant model structure $\gc\K=\<\K,\R,{}^*,0\>$, define
  operations $\circ$, $\to$, ${}^*$, and $\rmin$ on subsets
  $\X,\Y\subseteq\K$ by
  \begin{align}
    \label{circ}
    \X\circ\Y&=\{\z:\text{$\R\x\y\z$ for some $\x\in\X$ and $\y\in\Y$}\},\\
    \label{to}
    \X\to\Y&=\{\z:\text{if $\R\z\x\y$ and $\x\in\X$ then $\y\in\Y$}\},\\
    \label{star}
    \X^*&=\{\z^*:\z\in\X\},\\
    \label{rmin}
    \rmin\X&=\K\setminus\X^*.
  \end{align}
  For any valuation $\nu$ in $\gc\K$ with associated interpretation
  $\I$, let
  \begin{equation*}
    \J_\nu(\A)=\J(\A)=\{\c:\I(\A,\c)=\T\}
  \end{equation*}
  for every formula $\A$.
\end{definition}
These operations (and their notation) are designed for the following
consequences of Definition~\ref{def5}. For all formulas $\A$ and $\B$,
\begin{align*}
  \J(\A\land\B)&=\J(\A)\cap\J(\B),\\
  \J(\A\lor\B)&=\J(\A)\cup\J(\B),\\
  \J(\A\to\B)&=\J(\A)\to\J(\B),\\
  \J(\A\circ\B)&=\J(\A)\circ\J(\B),\\
  \J(\rmin\A) &=\rmin\J(\A),
\end{align*}
and $\A$ is valid in $\gc\K$ if $0\in\J(\A)$ for every valuation on
$\gc\K$.  Some useful observations to make at this point are, for all
$\X,\Y\subseteq\K$,
\begin{itemize}
\item $\X\circ\emptyset=\emptyset\circ\X=\emptyset$,
\item $\X\circ(\Y\cup\Z)=\X\circ\Y\,\,\cup\,\,\X\circ\Z$,
\item $(\Y\cup\Z)\circ\X=\Y\circ\X\,\,\cup\,\,\Z\circ\X$,
\item $\X\to\Y=\rmin(\X\circ\rmin\Y)$.
\end{itemize}
Every relevant model structure $\gc\K=\<\K,\R,{}^*,0\>$ has an
associated algebra, called its ``complex algebra''. The elements of
the complex algebra of $\gc\K$ are all the subsets of $\K$, and the
operations of the complex algebra are $\cup$, $\cap$, and the
operations $\circ$, $\to$, and $\rmin$ from Definition~\ref{def6}.

\begin{table}
\begin{align*}\gc\K_1&=
\begin{array}{|l|llll|}\hline
  \circ	&\{0\}		&\{a\}		&\{b\}		&\{b^*\}\\\hline
\{0\}	&\{0\}		&\{a\}		&\{b\}		&\{b^*\}\\
\{a\}	&\{a\}		&\{0,a,b\}	&\{b,b^*\}	&\{a,b,b^*\}\\
\{b\}	&\{b\}		&\{b,b^*\}	&\{a,b,b^*\}	&\{0,a,b,b^*\}\\
\{b^*\}	&\{b^*\}	&\{a,b,b^*\}	&\{0,a,b,b^*\}	&\{a,b,b^*\}\\\hline
\end{array}\\
\gc\K_2&=
\begin{array}{|l|llll|}\hline
  \circ	&\{0\}		&\{a\}		&\{b\}		&\{b^*\}\\\hline
\{0\}	&\{0\}		&\{a\}		&\{b\}		&\{b^*\}\\
\{a\}	&\{a\}		&\{0,a,b,b^*\}	&\{a,b,b^*\}	&\{a,b,b^*\}\\
\{b\}	&\{b\}		&\{a,b,b^*\}	&\{a,b,b^*\}	&\{0,a,b,b^*\}\\
\{b^*\}	&\{b^*\}	&\{a,b,b^*\}	&\{0,a,b,b^*\}	&\{a,b^*\}\\\hline
\end{array}\\
\gc\K_3&=
\begin{array}{|l|llll|}\hline
  \circ	&\{0\}		&\{a\}		&\{b\}		&\{b^*\}\\\hline
\{0\}	&\{0\}		&\{a\}		&\{b\}		&\{b^*\}\\
\{a\}	&\{a\}		&\{0,a,b,b^*\}	&\{a,b,b^*\}	&\{a,b,b^*\}\\
\{b\}	&\{b\}		&\{a,b,b^*\}	&\{a,b,b^*\}	&\{0,a,b,b^*\}\\
\{b^*\}	&\{b^*\}	&\{a,b,b^*\}	&\{0,a,b,b^*\}	&\{a,b,b^*\}\\\hline
\end{array}
\end{align*}
\caption{Three normal relevant model structures.}\label{table4}
\end{table}
\section{A formula in $\L_3$ but not $\RR$}
Lemma~\ref{reflection} shows that \lref{reflection} in
Table~\ref{table2} is in Tarski's relevance logic $\FL$.  In fact,
this formula is already part of $\L_3$.  However, \lref{reflection} is
not a theorem of $\RR$. Here we present two normal relevant model
structures that invalidate an instance of \lref{reflection}, namely
$\A\to\B$, where $\p,\q,\r,\s\in\var$,
\begin{align*}
  \A&=(\p\circ\q)\land\r,\\
  \B&=((\p\land\rmin\s)\circ\q)\lor(\p\circ(\q\land(\s\circ\r))).
\end{align*}

Let $\K=\{0,\a,\b,\b^*\}$, where $|\K|=4$, $0^*=0$, $\a^*=\a$, and
${}^*$ interchanges $\b$ and $\b^*$, as suggested by the notation.
Two relevant model structures on $\K$, $\gc\K_1$ and $\gc\K_2$, are
obtained by using two ternary relations $\R$ on $\K$.  The ternary
relation for $\gc\K_1$ has 34 triples, while the ternary relation for
$\gc\K_2$ has 36 triples.  Table~\ref{table4} lists the
$\circ$-products of all singleton subsets of $\K$ in both structures.
The products for larger sets can be computed by using the distributive
laws listed above. The triples can be read from the tables. For
example, $\<\a,\a,\b^*\>$ is a triple in the ternary relation of
$\gc\K_2$ but not $\gc\K_1$ because
$\b^*\in\{\0,\a,\b,\b^*\}=\{\a\}\circ\{\a\}$ in the table for
$\gc\K_2$, while $\b^*\notin\{\0,\a,\b\}=\{\a\}\circ\{\a\}$ in the
table for $\gc\K_1$.

Neither $\gc\K_1$ nor $\gc\K_2$ is the atom structure of a relation
algebra. Their ternary relations fail to have the property, possessed
by all atom structures of relation algebras, that
$\R\x\y\z\iff\R\z\y^*\x$.  In $\gc\K_1$, the triples $\<\a,\a,\b\>$,
$\<\a,\b^*,\a\>$, and $\<\b^*,\a,\a\>$ are present, but
$\<\b,\a,\a\>$, $\<\a,\b,\a\>$, and $\<\a,\a,\b^*\>$ are missing.  In
$\gc\K_2$, $\<\b,\b,\b^*\>$ is present but $\<\b^*,\b^*,\b\>$ is
missing. Adding the missing triples to either structure produces
$\gc\K_3$ in Table~\ref{table4}.

Now we proceed to use $\gc\K_1$ and $\gc\K_2$ to show that
\lref{reflection} in Table~\ref{table2} is not a theorem of $\RR$.
For $\gc\K_1$, choose valuation $\nu$ so that $\J(\p)=\J(\s)=\{\a\}$,
that is,
\begin{align*}
  \nu(\p,\a)&=\nu(\s,\a)=\T,\\
  \nu(\p,0)&=\nu(\s,0)=\nu(\p,\b)=\nu(\p,\b^*)=\nu(\s,\b)=\nu(\s,\b^*)=\F.
\end{align*}
Then $\J(\rmin\s)=\{0,\b,\b^*\}$ so 
\begin{align*}
  \J(\p\land\rmin\s)&=\{\a\}\cap\{0,\b,\b^*\}=\emptyset,
\end{align*}
hence 
\begin{align*}
  \J((\p\land\rmin\s)\circ\q)&=\J(\p\land\rmin\s)\circ\J(\q)=
  \emptyset\circ\J(\q)=\emptyset,
\end{align*}
regardless of the
action of $\nu$ on $\q$. Let $\J(\q)=\{\a\}$. Then 
\begin{equation*}
  \J(\p\circ\q)=\J(\p)\circ\J(\q)=\{\a\}\circ\{\a\}=\{0,\a,\b\}.
\end{equation*}
Next, let $\J(\r)=\{\b\}$. Then
\begin{align*}
  \J(\s\circ\r)&=\J(\s)\circ\J(\r)=\{\a\}\circ\{\b\}
  =\{\b,\b^*\},\\
  \J(\q\land(\s\circ\r))&=\J(\q)\cap\,\,\J(\s\circ\r)=
  \{\a\}\cap\{\b,\b^*\}=\emptyset.
\end{align*}
This last equation gives us
\begin{equation*}
  \J(\p\circ(\q\land(\s\circ\r)))=
  \J(\p)\circ\J(\q\land(\s\circ\r))=\J(\p)\circ\emptyset=\emptyset,
\end{equation*}
regardless of our choice for $\J(\p)$, and this, together with
$\J((\p\land\rmin\s)\circ\q)=\emptyset$, gives us
\begin{equation*}
  \J(\B)=\J((\p\land\rmin\s)\circ\q)\cup\J(\p\circ(\q\land(\s\circ\r)))
  =\emptyset\cup\emptyset=\emptyset.
\end{equation*}
However, we also have
\begin{align*}
\J(\A)&=\J((\p\circ\q)\land\r)\\
&=\J(\p\circ\q)\cap\J(\r)\\
&=\{0,\a,\b\}\cap\{\b\}\\
&=\{\b\}
\end{align*}
Now, by definition, $\A\to\B$ is verified if $\I(\A\to\B,0)=\T$. This
means that for all $\x,\y\in\K$, if $\R0\x\y$ and $\I(\A,\x)=\T$ then
$\I(\B,\y)=\T$.  However, from the table we have
$\{0\}\circ\{\b\}=\{\b\}$, which tells us that $\R0\b\b$ by the
definition of the operation $\circ$, and $\I(\A,\b)=\T$ since
$\J(\A)=\{\b\}$, so we ought to have $\I(\B,\b)=\T$ if $\A\to\B$ were
verified, but we don't, because $\J(\B)=\emptyset$. By the
Routley-Meyer completeness results mentioned earlier, we conclude that
$\A\to\B$ is not a theorem of $\RR$.

For $\gc\K_2$, we repeat all these steps, put with different values.
Choose $\nu$ so that $\J(\p)=\{\b\}$ and $\J(\s)=\{\b^*\}$.  Then
$\J(\rmin\s)=\{0,\a,\b^*\}$,
\begin{align*}
  \J(\p\land\rmin\s)&=\{\b\}\cap\{0,\a,\b^*\}=\emptyset,
\end{align*}
hence
\begin{align*}
  \J((\p\land\rmin\s)\circ\q)&=\J(\p\land\rmin\s)\circ\J(\q)=
  \emptyset\circ\J(\q)=\emptyset,
\end{align*}
regardless of the action of $\nu$ on $\q$. Let $\J(\q)=\{\b\}$.  Then
\begin{equation*}
  \J(\p\circ\q)=\J(\p)\circ\J(\q)=\{\b\}\circ\{\b\}=\{\a,\b,\b^*\}.
\end{equation*}
Let $\J(\r)=\{\b^*\}$.  Then
\begin{align*}
  \J(\s\circ\r)&=\J(\s)\circ\J(\r)=\{\b^*\}\circ\{\b^*\}=\{\a,\b^*\},\\
  \J(\q\land(\s\circ\r))&=\J(\q)\cap\J(\s\circ\r)
  =\{\b\}\cap\{\a,\b^*\}=\emptyset,
\end{align*}
hence
\begin{equation*}
  \J(\p\circ(\q\land(\s\circ\r)))=
  \J(\p)\circ\J(\q\land(\s\circ\r))=\J(\p)\circ\emptyset=\emptyset,
\end{equation*}
regardless of our choice for $\J(\p)$. Together with
$\J((\p\land\rmin\s)\circ\q)=\emptyset$, this gives us
\begin{equation*}
  \J(\B)=\J((\p\land\rmin\s)\circ\q)\cup\J(\p\circ(\q\land(\s\circ\r)))
  =\emptyset\cup\emptyset=\emptyset.
\end{equation*}
We also have
\begin{align*}
\J(\A)&=\J((\p\circ\q)\land\r)\\
&=\J(\p\circ\q)\cap\J(\r)\\
&=\{\a,\b,\b^*\}\cap\{\b\}\\
&=\{\b\}
\end{align*}
From the table we have $\{0\}\circ\{\b\}=\{\b\}$, hence $\R0\b\b$, and
$\I(\A,\b)=\T$ since $\J(\A)=\{\b\}$. We ought to get $\I(\B,\b)=\T$
if $\A\to\B$ were verified, but we don't since $\J(\B)=\emptyset$.  By
the Routley-Meyer completeness results, $\A\to\B$ is not a theorem of
$\RR$.

By the way, all three relevant model structures $\gc\K_1$, $\gc\K_2$,
and $\gc\K_3$ validate formulas \eqref{contra}, \eqref{perm},
\eqref{suff}, \eqref{mp}, \eqref{contr}, and \eqref{reduc}, but
\eqref{ming} is invalidated in many ways. For example, choosing $\nu$
so that $\J(\p)=\{\a\}$ yields the same calculations in all three
structures:
\begin{align*}
  \J(\p\to\p)
  &= \J(\p)\to\J(\p)\\
  &=\rmin(\{\a\}\circ\rmin\{\a\})\\
  &=\rmin(\{\a\}\circ\{0,\b,\b^*\})\\
  &=\rmin\{\a,\b,\b^*\}\\
  &=\{0\}\\
  \J(\p\to(\p\to\p))
  &= \rmin(\J(\p)\circ\rmin\J(\p\to\p))\\
  &= \rmin(\{\a\}\circ\rmin\{0\})\\
  &= \rmin(\{\a\}\circ\{\a,\b,\b^*\})\\
  &= \rmin\{0,\a,\b,\b^*\}\\
  &= \emptyset
\end{align*}
\section{A normal relevant model structure on a 21-element group}
The two normal relevant model structures $\gc\K_1$ and $\gc\K_2$ are
the simplest of several found by Prover9/Mace4 \cite{prover9-mace4}.
Although neither of them is the atom structure of a relation algebra,
they turn out to be very nearly the same proper relation algebra. If
three more triples are added to $\gc\K_1$, namely $\<\b,\a,\a\>$,
$\<\a,\b,\a\>$, and $\<\a,\a,\b^*\>$, or if one more triple is added
to $\gc\K_2$, namely $\<\b^*,\b^*,\b\>$, they both become the normal
relevant model structure $\gc\K_3$, also shown in Table~\ref{table4}.
\begin{table}
\begin{equation*}
\begin{array}{|ll|}\hline\stand
  \a=\{\f,\f^2,\g,\g^2,\g^5,\g^6\}&\\
  \b=&\b^*=\\
  \{\f\g,\f^2\g,\f\g^2,\f\g^3,\g^4,\f\g^4,\f\g^6\}
  &\{\g^3,\f^2\g^2,\f^2\g^3,\f^2\g^4,\f\g^5,\f^2\g^5,\f^2\g^6\}\\
  \{\f\g,\g^3,\f\g^3,\f\g^4,\f\g^5,\f\g^6,\f^2\g^6\}
  &\{\f^2\g,\f\g^2,\f^2\g^2,\g^4,\f^2\g^3,\f^2\g^4,\f^2\g^5\}\\
  \{\f\g,\g^3,\f\g^3,\f^2\g^4,\f\g^5,\f^2\g^5,\f^2\g^6\}
  &\{\f^2\g,\f\g^2,\f^2\g^2,\g^4,\f^2\g^3,\f\g^4,\f\g^6\}\\
  \{\f\g,\f^2\g^2,\g^4,\f\g^4,\f^2\g^4,\f\g^5,\f^2\g^6\}
  &\{\f^2\g,\f\g^2,\g^3,\f\g^3,\f^2\g^3,\f^2\g^5,\f\g^6\}\\
  \{\f\g,\f\g^3,\g^4,\f^2\g^4,\f\g^5,\f^2\g^5,\f^2\g^6\}
  &\{\f^2\g,\f\g^2,\g^3,\f^2\g^2,\f^2\g^3,\f\g^4,\f\g^6\}\\
  \hline\stand
  \a=\{\f,\f^2,\f\g,\f\g^2,\f^2\g^3,\f^2\g^6\}&\\
  \b=&\b^*=\\
  \{\g,\f^2\g,\g^3,\f^2\g^2,\f\g^4,\g^5,\f^2\g^4\}
  &\{\g^2,\f\g^3,\g^4,\f\g^5,\g^6,\f^2\g^5,\f\g^6\}\\
  \{\g,\f^2\g,\g^3,\f^2\g^2,\f\g^4,\g^5,\f\g^6\}
  &\{\g^2,\f\g^3,\g^4,\f^2\g^4,\f\g^5,\g^6,\f^2\g^5\}\\
  \{\g,\f^2\g,\g^3,\f\g^3,\f\g^4,\g^5,\f\g^6\}
  &\{\g^2,\f^2\g^2,\g^4,\f^2\g^4,\f\g^5,\g^6,\f^2\g^5\}\\\hline
\end{array}
\end{equation*}
\caption{Representations on $\mathfrak{K}_3$ on a 21-element group.}
\label{table5}
\end{table}
$\gc\K_3$ coincides with the relation algebra called \alg{37}{37} in
\cite{Maddux2006}.  This relation algebra is representable, as was
first shown by Stephen D.\ Comer ~\cite{Comer1979}, and is actually
isomorphic to a proper relation algebra whose base set is a 21-element
group, first shown by Peter Jipsen.  

Let $\G$ be the group generated by $\f$ and $\g$, subject to the
relations $\f^3=\g^7=1$ and $\g\f=\f\g^2$, where $1$ is the identity
element of $\G$.  Alternatively, let $\G$ be the group generated by
the permutations of $\{1,\cdots,21\}$ defined by
\begin{align*}
  \f&=(3,6,12)(5,8,14)(7,10,16)(9,18,15)(11,20,17)(13,21,19),\\
  \g&=(2,20,17,14,11,8,5)(4,16,7,19,10,21,13).
\end{align*}
Up to isomorphism, there are 8 ways (found by GAP~\cite{GAP4}) to
obtain $\gc\K_3$ from $\G$.  First let $0=\{1\}$ be the singleton
containing the identity element of $\G$.  Next, choose one of the 8
partitions listed in Table~\ref{table5} of the 20 non-identity
elements into 3 sets $\a$, $\b$, and $\b^*$.  Then $0$ and $\a$ are
closed under the formation of inverses in $\G$, $\b$ is the set of
inverses of elements in $\b^*$, and {\it vice versa},
\begin{align*}
  0&=\{1^{-1}\}&\a&=\{\h^{-1}\colon\h\in\a\},
  &\b^*&=\{\h^{-1}\colon\h\in\b\},&\b&=\{\h^{-1}\colon\h\in\b^*\}.
\end{align*}
Here the Routley star ${}^*$ is the operation of forming all the
inverses of the elements in a subset of $\G$.  For any
$\x,\y,\z\in\K=\{0,\a,\b,\b^*\}$, let the ternary relation $\R$ hold
on the triple $\<\x,\y,\z\>$ just in case $\z$ is included in the set
products of elements from $\x$ and $\y$, that is,
\begin{equation*}
  \R\x\y\z\iff\x\y\supseteq\z
  \iff\z\subseteq\{\h\k\colon\h\in\x,\,\k\in\y\},
\end{equation*}
where $\h\k$ is the product in $\G$ of the two group elements
$\h,\k\in\G$ and $\x\y$ is the set of products in $\G$ of pairs of
elements, one from $\x$ and one from $\y$, in that order.  This
completes the construction of $\gc\K_3$ from $\G$. (Infinitely many
other groups can be used in a similar way; $\G$ is just the smallest
one.) Every choice of partition from Table~\ref{table5} produces that
same table for the operation $\circ$ in $\gc\K_3$, as defined in
\eqref{circ} and shown in Table~\ref{table4}.

To show that this relevant model structure $\gc\K_3$ is isomorphic to
a proper relation algebra we use the right regular representation of
the group $\G$, as is done in the proof of the Cayley representation
theorem (every group is isomorphic to a group of permutations).  First
we recall that subsets of a group were once called ``complexes'', and
the set of subsets of the group $\G$ forms an algebra called its
``complex algebra'', whose operations are union, intersection,
complementation with respect to $\G$, multiplication of complexes as
defined above, and the operation of forming all the inverses of
elements in a subset of $\G$, here denoted by the Routley star ${}^*$.
The complex algebra also has, as a distinguished element, the
singleton consisting of just the identity element of the group.  For
every $\x\subseteq\G$, define the binary relation $\sigma(\x)$ on $G$
by
\begin{equation*}
  \sigma(\x)=\{\<\k,\k\h\>:\k\in\G,\h\in\x\}\subseteq\G\times\G.
\end{equation*}
Then $\sigma$ is an injective homomorphism from the complex algebra of
$\G$ into the proper relation algebra of all binary relations on
$\G$, in the sense that, for all $\x,\y\subseteq\G$, recalling
definitions \eqref{relprod} and \eqref{conv}, we have
\begin{align*}
  \sigma(\x\cup\y)&=\sigma(\x)\cup\sigma(\y),\\
  \sigma(\x\cap\y)&=\sigma(\x)\cap\sigma(\y),\\
  \sigma(\G\setminus\x)&=(\G\times\G)\setminus\sigma(\x),\\
  \sigma(\x\y)&=\sigma(\x)|\sigma(\y),\\
  \sigma(\x^*)&=\sigma(\x)^{-1},\\
  \sigma(\{1\})&=\{\<\h,\h\>\colon\h\in\G\}.
\end{align*}
If $\h\in\G$ then $\sigma(\{\h\})$ is the permutation used in the
proof of Cayley's theorem.  The right regular representation has a
property required by representations of relation algebras: the
permutations associated with $\{\h\}$ and $\{\k\}$ must be disjoint
(as sets) whenever $\h\neq\k$, simply because
$\{\h\}\cap\{\k\}=\emptyset$ and this fact must be reflected in any
representation.  Applying $\sigma$ to the elements of $\gc\K_3$
produces four binary relations on the 21-element set $\G$,
\begin{align*}
  \sigma(0)&=\{\<\h,\h\>\colon\h\in\G\},&
  \A&=\sigma(\a),&\B&=\sigma(\b),&\conv\B&=\sigma(\b^*).
\end{align*}
These four relations form a partition of $\G\times\G$. One of them is
the identity relation on $\G$, and the converse of any one of them is
either itself or another one of them. In particular,
$\conv{\sigma(0)}=\sigma(0)$ and $\conv\A=\A$, while $\B$ and
$\conv\B$ are converses of each other.  The table for $\gc\K_3$ yields
these conclusions about the relative products of these relations:
$\sigma(0)|\x=\x|\sigma(0)=\x$ for all
$\x\in\{\sigma(0),\A,\B,\conv\B\}$,
\begin{equation*}
  \A|\A=\B|\conv\B=\conv\B|\B=\G\times\G,
\end{equation*} 
and all other products are equal to $(\G\times\G)\setminus\sigma(0)$.
Any four relations with these properties gives us yet another
representation of the relation algebra \alg{37}{37}, alias $\gc\K_3$.
\section{Counterexample to a theorem of Kowalski}
Kowalski \cite[Theorem 8.1]{MR3289545} proved, ``The relevant logic
$\RR$ is sound and complete with respect to square-increasing,
commutative, integral relation algebras.'' \lref{reflection} is a
counterexample to this theorem.  It is not a theorem of $\RR$ (because
it is invalid the the relevant model structures $\gc\K_1$ and
$\gc\K_2$) and yet holds in all relation algebras (including the
square-increasing, commutative, integral ones).  In fact,
Lemma~\ref{reflection} shows that it is in $\L_3$ and is true in all
semi-associative relation algebras (because it is provable with only
three objects).

\cite[Theorem 8.1]{MR3289545} is obtained as an immediate consequence
of \cite[Theorem 7.1]{MR3289545}, that ``Every normal De Morgan monoid
is embeddable as a bare [no constants] De Morgan monoid into a
square-increasing, commutative, integral relation algebra.'' The
complex algebras of $\gc\K_1$ and $\gc\K_2$ are counterexamples to
Theorem~7.1.

Part of the proof of Theorem~7.1 reads, ``By definition of
$\varepsilon$, it is an embedding of the lattice reduct of [the normal
De Morgan monoid] $\M$ into the lattice reduct of [the relation
algebra] $\U_\M$; in particular, $\varepsilon$ is injective.
Lemma~7.1 shows that the multiplication, implication and De Morgan
negation are preserved as well, \dots''. The proof of Theorem~7.1 uses
Lemma~5.4 to show half of the preservation of multiplication, namely
$\varepsilon(\a\b)\subseteq\varepsilon(\a)\circ\varepsilon(\b)$.  

The difficulty arises in the proof of Lemma~5.4(1) at this point:
``Since $\M$ is a distributive lattice, $\R$ is a prime filter and
$\R'=\overline{\R}$, proving (1).''  At this stage in the proof,
$\<\R,\R'\>$ is known to be a maximal disjoint pair in $S$ (hence $\R$
a proper filter and $\R'$ a proper ideal disjoint from $\R$, and they
satisfy two additional technical conditions).  If the desired
conclusion that $\overline{R}=R'$ were to fail, there would be some
element $\x\notin\R\cup\R'$.  A desired contradiction could then be
attained by showing that $\x$ could be added to $\R$ or to $\R'$, \ie,
either the filter generated by $\x$ and $\R$ is disjoint from the
ideal $\R'$, or else the ideal generated by $\R'$ and $x$ is disjoint
from $\R$. The distributivity of $\M$ insures that one of these two
possibilities happens. For example, if the filter $R^+$ generated by
$\x$ and $\R$ is disjoint from the ideal $\R'$, but the ideal
generated by $\R'$ and $\x$ is not disjoint from $\R$, then
$\<\R^+,\R'\>$ would be a strictly larger disjoint filter-ideal pair.
This would yield a contradiction if $\<\R,\R'\>$ were maximal among
\emph{all} disjoint filter-ideal pairs, but it is only known to be
maximal in $S$ (and subject to the technical conditions). The goal
would be to show that $\<\R^+,\R'\>$ is actually a strictly larger
pair in $S$ (that this pair also satisfies the technical conditions).
The difficulty in achieving this goal would be revealed by a more
detailed examination of this situation for the complex algebras of
$\gc\K_1$ and $\gc\K_2$.
\section{Deriving \lref{reflection} from Tarski's axioms}
\lref{reflection} in Table \ref{table2} is shown to be in $\L_3$ by
Lemma~\ref{reflection}. It is the translation into relevance logic
notation of the following equation, which is true in all relation
algebras.
\begin{align}\label{refleq}
  \x\rp\y\cdot\z
  &\leq(\x\cdot\min{\con\w})\rp\y+\x\rp(\y\cdot\w\rp\z).
\end{align}
The equation \eqref{refleq} can therefore be derived from Tarski's ten
axioms for relation algebras~\cite[8,2(i)]{MR920815} (treated as
algebras of the form $\<\U,+,\min\blank,\,\rp\,,\con\ ,\id\>$).
Tarski's axioms are
\begin{align*}
  \tag{\ra1}\x+\y&=\y+\x,\\
  \tag{\ra2}\x+(\y+\z)&=(\x+\y)+\z,\\
  \tag{\ra3}\min{\min\x+\min\y}+\min{\min\x+\y}&=\x,\\
  \tag{\ra4}\x\rp(\y\rp\z)&=(\x\rp\y)\rp\z,\\
  \tag{\ra5}(\x+\y)\rp\z&=\x\rp\z+\y\rp\z,\\
  \tag{\ra6}\x\rp\id&=\x,\\
  \tag{\ra7}\con{\con\x}&=\x,\\
  \tag{\ra8}\con{\x+\y}&=\con\x+\con\y,\\
  \tag{\ra9}\con{\x\rp\y}&=\con\y\rp\con\x,\\
  \tag{\ra{10}}\con\x\rp\min{\x\rp\y}+\min\y&=\min\y.
\end{align*}
These are the axioms about which Tarski asked ``whether this
definition of relation algebra \dots\ is justified in any intrinsic
sense.'' 

The first three axioms are a set of postulates for Boolean algebras
(treated as algebras of the form $\<\U,+,\min\blank\>$) due to E.\ V.\
Huntington~\cite{MR1501702, MR1501684, MR1501729}.  Proving all the
usual equations true in Boolean algebras from the Huntington axioms is
an interesting and challenging homework problem.  One must first prove
$\x+\min\x=\y+\min\y$ in order to define the maximum element $1$ by
$1=\x+\min\x$. (See \cite{MR1392462} for a solution.)  We will prove
the following purely relation-algebraic facts directly from Tarski's
axioms.  In the following proofs, any step that requires only Boolean
algebra is marked ``BA''.
\begin{align}
  \label{ra1}
  \x\leq\y&\to\x\rp\z\leq\y\rp\z\\
  \label{ra2}
  \z\rp(\x+\y)&=\z\rp\x+\z\rp\y\\
  \label{ra3}
  \x\leq\y&\to\z\rp\x\leq\z\rp\y\\
  \label{ra4}
  \con1&=1\\
  \label{ra5}
  \min{\con\x}&=\con{\min\x}\\
  \label{ra6}
  \con{\x\cdot\y}&=\con\x\cdot\con\y\\
  \label{ra7}
  \x\rp\y\bp\z&\leq\x\rp(\y\bp\con\x\rp\z)
\end{align}
Proof of \eqref{ra1}:
\begin{align*}
  \x\leq\y&\iff\x+\y=\y&&\text{BA}\\
  &\to(\x+\y)\rp\z=\y\rp\z\\
  &\iff\x\rp\z+\y\rp\z=\y\rp\z&&\thetag{\ra5}\\
  &\iff\x\rp\z\leq\y\rp\z&&\text{BA}
\end{align*}
Proof of \eqref{ra2}:
\begin{align*}
  \z\rp(\x+\y)
  &=\con{\con{\z\rp(\x+\y)}}&&\thetag{\ra7}\\
  &=\con{\con{\x+\y}\rp\con\z}		&&\thetag{\ra9}\\
  &=\con{(\con\x+\con\y)\rp\con\z}	&&\thetag{\ra8}\\
  &=\con{\con\x\rp\con\z+\con\y\rp\con\z}&&\thetag{\ra5}\\
  &=\con{\con{\z\rp\x}+\con{\z\rp\y}}	&&\thetag{\ra9}\\
  &=\con{\con{\z\rp\x+\z\rp\y}}		&&\thetag{\ra8}\\
  &=\z\rp\x+\z\rp\y &&\thetag{\ra7}
\end{align*}
The proof of \eqref{ra3} is like the proof of \eqref{ra2}, but turned
around in the obvious way. Proof of \eqref{ra4}:
\begin{align*}
 1&=	1+\con1			&&\text{BA}\\
 &=	\con{\con1}+\con1	&&\thetag{\ra7}\\
 &=	\con{\con1+1}		&&\thetag{\ra8}\\
 &=	\con1			&&\text{BA}
\end{align*}
For \eqref{ra5}, first note that, for any $\y$, the following
statements are equivalent.
\begin{align*}
  &\con{\min\x}\leq\y\\
  &\con{\min\x}+\y=\y&&\text{BA}\\
  &\min\x+\con\y=\con\y&&\thetag{\ra7}, \thetag{\ra8}\\
  &\x+\con\y=1&&\text{BA}\\
  &\con\x+\y=1&&\thetag{\ra7}, \thetag{\ra8}, \eqref{ra4}\\
  &\min{\con\x}\leq\y&&\text{BA}
\end{align*}
We need only two instances of these equivalences. When $\y$ is either
$\con{\min\x}$ or $\min{\con\x}$, we deduce that
$\min{\con\x}\leq\con{\min\x}$ and $\con{\min\x}\leq\min{\con\x}$,
respectively, hence \eqref{ra5} holds. Proof of \eqref{ra6}:
\begin{align*}
  \con{\x\cdot\y}
  &=\con{\min{\min\x+\min\y}}		&&\text{BA}\\
  &=\min{\con{\min\x+\min\y}}		&&\eqref{ra5}\\
  &=\min{\con{\min\x}+\con{\min\y}}	&&\thetag{\ra8}\\
  &=\min{\min{\con\x}+\min{\con\y}}	&&\eqref{ra5}\\
  &=\con\x\cdot\con\y &&\text{BA}
\end{align*}
Proof of \eqref{ra7}:
\begin{align*}
  \x\rp\y
  &=\x\rp(\y\bp(\con\x\rp\z+\min{\con\x\rp\z}))&&\text{BA}\\
  &=\x\rp(\y\bp\con\x\rp\z+\y\bp\min{\con\x\rp\z})&&\text{BA}\\
  &=\x\rp(\y\bp\con\x\rp\z)+\x\rp(\y\bp\min{\con\x\rp\z})&&\eqref{ra2}\\
  &\leq\x\rp(\y\bp\con\x\rp\z)+\x\rp(\min{\con\x\rp\z})&&\eqref{ra3}\\
  &\leq\x\rp(\y\bp\con\x\rp\z)+\min\z&&\thetag{\ra{10}}
  \intertext{From the previous equation we get} \x\rp\y\bp\z
  &\leq(\x\rp(\y\bp\con\x\rp\z)+\min\z)\bp\z&&\text{BA}\\
  &=\x\rp(\y\bp\con\x\rp\z)\bp\z+\min\z\bp\z&&\text{BA}\\
  &=\x\rp(\y\bp\con\x\rp\z)\bp\z+0&&\text{BA}\\
  &=\x\rp(\y\bp\con\x\rp\z)\bp\z&&\text{BA}
\end{align*}
Associativity is not needed in any form for the proof of
\eqref{refleq}.  Consequently \eqref{refleq} holds in all
non-associative relation algebras (the class of algebras obtained by
dropping \thetag{\ra4} from the list of axioms).  Here is a direct
equational proof of \eqref{refleq}.
\begin{align*}
  \x\rp\y\cdot\z &=(\x\cdot(\min{\con\w}+\con\w))\rp\y\cdot\z
  &&\text{BA}\\
  &=(\x\cdot\min{\con\w}+\x\cdot\con\w)\rp\y\cdot\z
  &&\text{BA}\\
  &=((\x\cdot\min{\con\w})\rp\y+(\x\cdot\con\w)\rp\y)\cdot\z
  &&\thetag{\ra5}\\
  &=(\x\cdot\min{\con\w})\rp\y\cdot\z+(\x\cdot\con\w)\rp\y\cdot\z
  &&\text{BA}\\
  &\leq(\x\cdot\min{\con\w})\rp\y+(\x\cdot\con\w)\rp\y\cdot\z
  &&\text{BA}\\
  &\leq(\x\cdot\min{\con\w})\rp\y
  +(\x\cdot\con\w)\rp(\y\cdot\con{\x\cdot\con\w}\rp\z)
  &&\eqref{ra7}\\
  &=(\x\cdot\min{\con\w})\rp\y
  +(\x\cdot\con\w)\rp(\y\cdot(\con\x\cdot\w)\rp\z)
  &&\eqref{ra6}, \thetag{\ra7}\\
  &\leq(\x\cdot\min{\con\w})\rp\y+\x\rp(\y\cdot\w\rp\z) &&\eqref{ra1},
  \eqref{ra3}
\end{align*}
\section{Variable-sharing}
Tarski's relevance logic has the variable-sharing property, even if
extended beyond $\RR$ by adding axioms insuring commutativity and
density.  Belnap's~\cite{MR0141590} original proof of this fact for
the logic $\sf\E$ of Anderson-Belnap~\cite{MR0115902} applies with no
changes. Belnap's construction and proof are presented in this
section.  Belnap gave matrices for $\land$, $\lor$, $\to$, $\rmin$,
and two defined unary connectives, $\N(\A)=(\A\to\A)\to\A$ and
$\M(\A)=\rmin(\N(\rmin\A))$. 

From the matrices for $\land$ and $\lor$ it is apparent that the eight
values appearing in them, namely $-3$, $-2$, $-1$, $-0$, $+0$, $+1$,
$+2$, and $+3$ (the last four are the designated values), form a
lattice isomorphic to the lattice of subsets of the 3-element set
$\{-1,+0,-2\}$, with $+3$ at the top and $-3$ at the bottom, if
$\land$ and $\lor$ are interpreted as intersection and union.  This
observation does not occur in~\cite{MR0141590}, but in subsequent
literature they are usually protrayed this way; see, for example,
\cite[pp.\,198,\,252]{MR0406756}, \cite[p.\,178]{Routleyetal1982},
and~\cite[p.\,102]{Brady2003}.

What took nearly half a century after their introduction in 1960 was
the realization in \cite{Maddux2007} that Belnap's matrices define a
proper relation algebra; see also \cite{KramerMaddux2019, MR2641636}.
This proper relation algebra was known to Lyndon~\cite{MR0037278} in
1950, and became well known in the 1980s under the name ``The Point
Algebra'', because it describes the ways two points on the real line
can be related to each other; the three atomic relations between two
real numbers are $\x<\y$, $\x=\y$, and $\x>\y$.  The joins of pairs of
these relations are $\leq$, $\geq$, and $\neq$.

Two formulas $\A,\B\in\sent$ are said to {\bf share a variable} if
some propositional variable $\p\in\var$ occurs in both $\A$ and $\B$.
\begin{table}
\begin{align*}\gc\K_4&=
\begin{array}{|l|lll|}\hline
  \circ	&\{0\}		&\{a\}		&\{a^*\}\\\hline
\{0\}	&\{0\}		&\{a\}		&\{a^*\}\\
\{a\}	&\{a\}		&\{a\}		&\{0,a,a^*\}\\
\{a^*\}	&\{a^*\}	&\{0,a,a^*\}	&\{a^*\}\\\hline
\end{array}
\end{align*}
\caption{Belnap's normal relevant model structure.}\label{table6}
\end{table}
To show $\A$ and $\B$ share a propositional variable whenever
$\A\to\B\in\FL$, we use the normal relevant model structure $\gc\K_4$
shown in Table~\ref{table6}.  Choose a valuation $\nu$ so that
\begin{align*}
  \J(\p)&=\begin{cases}
    \{\a\}&\text{ if $\p$ occurs in $\A$,}\\
    \{\a^*\}&\text{ if $\p$ does not occur in $\A$.}
  \end{cases}
\end{align*}
One key feature of $\gc\K_4$ is that $\{\{\a\},\{0,\a\}\}$ and
$\{\{\a^*\},\{0,\a^*\}\}$ are both closed under $\cup$, $\cap$, $\to$,
$\circ$, and $\rmin$. This is obvious for $\cup$ and $\cap$, clear for
$\circ$ from Table~\ref{table6}, easy to check for $\rmin$, and
therefore is also true for $\to$.  The other key feature is that
$\X\to\Y=\emptyset$ whenever $\X\in\{\{\a\},\{0,\a\}\}$ and
$\Y\in\{\{\a^*\},\{0,\a^*\}\}.$ For this we provide two sample
computations.
\begin{align*}
  \{\a\}\to\{\a^*\}
  &=\rmin(\{\a\}\circ\rmin\{\a^*\})=\rmin(\{\a\}\circ\{0,\a^*\})\\
  &=\rmin\{0,\a,\a^*\}=\K\setminus\{0,\a,\a^*\}=\emptyset,\\
  \{0,\a\}\to\{0,\a^*\}
  &=\rmin(\{0,\a\}\circ\rmin\{0,\a^*\})=\rmin(\{0,\a\}\circ\{\a^*\})\\
  &=\rmin\{0,\a,\a^*\}=\K\setminus\{0,\a,\a^*\}=\emptyset.
\end{align*}
By the choice of $\nu$, the closure of $\{\{\a\},\{0,\a\}\}$ give us
\begin{equation*}
  \J(\A)\in\{\{\a\},\{0,\a\}\}.
\end{equation*}
Suppose that $\B$ is a formula whose propositional variables do not
occur in $\A$. Then, by the choice of $\nu$ and the closure of
$\{\{\a^*\},\{0,\a^*\}\}$,
\begin{equation*}
  \J(\B)\in\{\{\a^*\},\{0,\a^*\}\}.
\end{equation*}
By the second key feature, we conclude that
$\J(\A\to\B)=\J(\A)\to\J(\B)=\emptyset$. Since $0$ is not in
$\emptyset$, $\A\to\B$ is not valid in $\gc\K_4$. The contrapositive
of what we have just proved is that if $\A\to\B$ is valid in
$\gc\K_4$, then $\A$ and $\B$ must share a variable.\endproof
\section{Representing Belnap's normal relevant model structure}
A representation of $\gc\K_4$ as the atom structure of a proper
relation algebra can be obtained as follows.  Let $\rationals$ be the
set of rational numbers. Let
\begin{align*}
\sigma(\a)&=\{\<\x,\y\>\colon\x<\y,\,\x,\y\in\rationals\},\\
\sigma(\a^*)&=\{\<\x,\y\>\colon\x>\y,\,\x,\y\in\rationals\},\\
\sigma(0)&=\{\<\x,\y\>\colon\x=\y,\,\x,\y\in\rationals\}.
\end{align*}
Extend $\sigma$ to all subsets of $\K=\{0,\a,\a^*\}$, by sending each
subset of $\K$ to the union of the images of its elements under
$\sigma$. For example,
\begin{align*}
  \sigma(\{\a\})&=\sigma(\a),\\
  \sigma(\{\a,\a^*\})&=\{\<\x,\y\>\colon\x\neq\y,\,\x,\y\in\rationals\},\\
  \sigma(\{0,\a\})&=\{\<\x,\y\>\colon\x\leq\y,\,\x,\y\in\rationals\}.
\end{align*}
Thus $\sigma$ maps the complex algebra of $\gc\K_4$ onto the proper
relation algebra whose universe consists of the eight binary relations
on the rationals usually denoted in a more colloquial notation as $=$,
$\neq$, $<$, $>$, $\leq$, $\geq$, $\emptyset$, and
$\rationals\times\rationals$.

\begin{table}
\begin{align*}
\gc\K_5=
&\begin{array}{|c|cccc|}\hline
\circ
	&\{0\}		&\{\a\}			&\{\b\}		&\{\b^*\}\\\hline
\{0\}	&\{0\}		&\{\a\}			&\{\b\}		&\{\b^*\}\\
\{\a\}	&\{\a\}		&\{0,\a,\b,\b^*\}	&\{\a,\b\}	&\{\a\}\\
\{\b\}	&\{\b\}		&\{\a\}			&\{\b\}		&\{0,\a,\b,\b^*\}\\
\{\b^*\}&\{\b^*\}	&\{\a,\b^*\}		&\{0,\b,\b^*\}	&\{\b^*\}\\\hline
\end{array}
\end{align*}
\caption{A non-commutative normal relevant model structure}\label{table7}
\end{table}
\section{Four axioms of $\RR$ not in $\L_4$}
Table \ref{table7} shows the atom structure of a non-commutative
proper relation algebra called $\alg{13}{37}$ in \cite{Maddux2006}.
It satisfies conditions \thetag{p1}, \thetag{p2}, \thetag{p4}, and
\thetag{p6} in Definition \ref{def4}, plus \thetag{p3$'$} and
\thetag{p5$'$}, and therefore has all the required properties to be a
relevant model structure except \thetag{comm}. It is normal since
$0^*=0$. It is therefore called a ``non-commutative normal relevant
model structure''. 

Although condition \thetag{p1} is called ``$0$-reflexivity'', it
insures that the proper relation algebra $\alg{13}{37}$ is dense, \ie,
satisfies $\x\leq\x^2$, where $\x^2=\x\rp\x$. Condition \thetag{p1}
should therefore be called ``density'', but the term
``square-increasing'' is commonly used instead because it describes
the shape of the equation that defines density.  Since $\gc\K_5$
satisfies \thetag{p1}, it validates the formulas that assert density
for all relations, namely the contraction axiom \eqref{contr} and the
reductio axiom \eqref{reduc}. On the other hand, since it is not
commutative, the axioms depending on that assumption are invalid in
$\gc\K_5$, namely contraposition \eqref{contra}, permutation
\eqref{perm}, suffixing \eqref{suff}, and modus ponens \eqref{mp}.
These formulas are invalidated in many ways, but in rather few ways if
the valuations are restricted so the propositional variables are
mapped to singletons and the formulas are mapped to the empty set.
Here is a complete list of such valuations (calculated with
GAP~\cite{GAP4}).
\begin{itemize}\item\eqref{contra} is invalid in $\gc\K_5$ because
  \begin{equation*}
    \J((\p\to\rmin\q)\to(\q\to\rmin\p))=\emptyset
  \end{equation*}
  whenever $\nu$ is chosen so that one of these three sets of
  equations holds:
\begin{align*}
  \J(\p)&=\{\a\}&\J(\q)&=\{\b\}\\
  \J(\p)&=\{\b\}&\J(\q)&=\{\b^*\}\\
  \J(\p)&=\{\b^*\}&\J(\q)&=\{\a\}\\
\end{align*}
\item \eqref{perm} is invalid in $\gc\K_5$ because
  \begin{equation*}
    \J((\p\to(\q\to\r))\to(\q\to(\p\to\r)))=\emptyset
  \end{equation*}
  whenever $\nu$ is chosen so that one of these two sets of equations
  holds:
\begin{align*}
  \J(\p)&=\{\a\}&\J(\q)&=\{\b\}&\J(\r)&=\{\a\}\\
  \J(\p)&=\{\b^*\}&\J(\q)&=\{\a\}&\J(\r)&=\{\a\}\\
\end{align*}
\item \eqref{suff} is invalid in $\gc\K_5$ because
  \begin{equation*}
    \J((\p\to\q)\to((\q\to\r)\to(\p\to\r)))=\emptyset
  \end{equation*}
  whenever $\nu$ is chosen so that one of these four sets of equations
  holds:
\begin{align*}
  \J(\p)&=\{0\}&\J(\q)&=\{\a\}&\J(\r)&=\{\a\}\\
  \J(\p)&=\{0\}&\J(\q)&=\{\b\}&\J(\r)&=\{\a\}\\
  \J(\p)&=\{\b\}&\J(\q)&=\{\a\}&\J(\r)&=\{\a\}\\
  \J(\p)&=\{\b\}&\J(\q)&=\{\b\}&\J(\r)&=\{\a\}
\end{align*}
\item \eqref{mp} is invalid in $\gc\K_5$ because
  \begin{equation*}
    \J(\p\to((\p\to\q)\to\q))=\emptyset
  \end{equation*}
  whenever $\nu$ is chosen so that one of these two sets of equations
  holds:
\begin{align*}
  \J(\p)&=\{\a\}&\J(\q)&=\{\a\}\\
  \J(\p)&=\{\b\}&\J(\q)&=\{\a\}
\end{align*}
\end{itemize}
\section{Representing $\mathfrak{K}_5$ as a proper relation algebra}
As with Belnap's normal relevant model structure, there is a
representation of $\gc\K_5$ as the atom structure of a proper relation
algebra.  Again, $\rationals$ is the set of rational numbers.  Let
$\U$ be the set of finite sequences of one or more rational numbers,
in which the first is arbitrary and all others are positive.  
Define a binary relation $\B\subseteq\U\times\U$ as follows.

Think of each element of $\U$ as representating a location, from which
it is possible to either travel some positive distance in ``the same
direction'', or to ``branch off'' and travel some positive distance in
``the new direction''. If $\s$ is the new point at which one arrives
by moving as described, then the pair $\<\r,\s\>$ is in the relation
$\B$. Finally, $\B$ is the transitive closure of the set of all pairs
obtained from this description.

More formally, an ordered pair $\<\r,\s\>$ of sequences $\r,\s\in\U$
is in $\B_0$ if and only if $\r\neq\s$ and either $\s$ can be obtained
from $\r$ by adding a nonnegative rational to the last entry of $\r$
(travel in the same direction by that amount) or appending a positive
rational number to the end of $\r$ (travel in the new direction by
that amount). Let $\B$ be the transitive closure of $\B_0$. Since $\B$
is a partial ordering, we will symbolize it with ``$<$'' in these
examples:
\begin{equation*}
  \<-8\><\<0\><\<1\><\<1,2,3\><\<1,2,4\><\<1,2,4,5\><\<1,2,4,5,6\><\dots
\end{equation*}
Let
\begin{align*}
  \sigma(\b^*)&=\B,\\
  \sigma(\b)&=\conv\B,\\
  \sigma(\a)&=\U\times\U\setminus(\B\cup\conv\B),\\
  \sigma(0)&=\{\<\x,\y\>\colon\x=\y,\,\x,\y\in\U\},
\end{align*}
and extend $\sigma$ to all subsets of $\K=\{0,\a,\b,\b^*\}$ by sending
each subset to the union of the images of its elements under $\sigma$.
\section{Axiomatizing classical relevant logic}
In \cite[p.\,183]{MR0363789a}, Meyer and Routley define a CR* model
structure $\gc\K=\<\K,\R,{}^*,0\>$.  Their definition is the same as
that of a normal relevant model structure except that conditions
\thetag{p1}, \thetag{p4}, and $0^*=0$ are replaced by $\R0\a\b\iff
\a=\b$, from which the three conditions can be derived (using the
remaining conditions \thetag{p2}, \thetag{p3}, \thetag{p5}, and
\thetag{p6}). Hence every CR* model structure is a normal relevant
model structure (but not conversely).

Their language contains connectives $\to$, $\land$, $\neg$, and ${}^*$
\cite[p.\,184]{MR0363789a}, while $\lor$ is recovered by the
definition $\A\lor\B=\neg(\neg\land\neg\B)$ \cite[d5.,
p.\,187]{MR0363789a} and $\rmin$ is defined by $\rmin\A=\neg(\A^*)$
\cite[d4., p.\,186]{MR0363789a}. The notions of valuation and
interpretation in Definition~\ref{def5} are suitably altered by
retaining the conditions pertaining to the connectives $\to$ and
$\land$, adding the conditions $\I(\neg\A,\c)=\T$ iff $\I(\A,\c)=\F$
and $\I(\A^*,\c)=\T$ iff $\I(\A,\c^*)=\T$, and deriving the conditions
for $\lor$, $\circ$, and $\rmin$ through their definitions.  Their
system CR* of classical relevant logic is defined as the set of
formulas valid in all CR* model structures \cite[(9),
p.\,185]{MR0363789a}.  In section III they
\begin{quote}
  ``\dots\ show that the system CR${}^*$, characterized so that its
  set of theorems is exactly the CR${}^*$ valid formulas, exactly
  contains the system $\RR$ of relevant implication on the definition
  of $\rmin$ by d4.'' \cite[p.\,187]{MR0363789a}
\end{quote}
This means that a formula $\A$, written in the language of the
connectives $\land$, $\lor$, $\to$, $\circ$, and $\rmin$, is a theorem
of $\RR$ if and only if the same formula, but with the connectives
$\lor$, $\circ$, and $\rmin$ defined in terms of $\land$, $\to$,
$\neg$, and ${}^*$, is valid in all CR${}^*$ model structures. Their
concluding remarks concern axiomatization.
\begin{quote}
  ``In conclusion, it will be noted that we have neglected to
  axiomatize CR${}^*$. The reason isn't that it's unaxiomatizable or
  anything like that; indeed, we presume that just putting together
  the axiomatization of CR in \cite{MR0363789} and of $\RR$ in
  \cite{MR0409114} or \cite{MR0406756} one would have an
  axiomatization of CR${}^*$, near enough, reversing d4 by then
  defining A${}^*$ as $\neg\rmin\A$. Frankly, however, we can't at
  this point stomach yet another completeness proof on ground that we
  have been over so often before; any readers that have stuck with us
  through the series of papers that began with \cite{MR0409114} feel
  as we do, no doubt, letting the semantic characterization of
  CR${}^*$ above suffice.  But the case is now pretty strong that
  $\neg$ was just left out of Anderson-Belnap formulations of their
  logics, and evidence is building that the entire project of relevant
  logic is unified and simplified when the \emph{semantic} $\neg$,
  with a different function from the \emph{deduction-theoretic} ---
  that has been present from the start, is added. This paper is part
  of that evidence.''
\end{quote}
Meyer and Routley \cite[p.\,53]{MR0363789} axiomatize the system
$\RR^+$ with axioms A1--A11, A14, and A15, and rules R1 and R2 from
\cite[p.\,204]{MR0409114}. These axioms and rules are the ones that do
not mention negation. To combine these with the axioms and rules of
\cite{MR0409114}, as they suggest, would seem to do nothing more than
restore axioms A12 and A13 that involve negation.  All these axioms
and rules are recounted in \cite[pp.340-1]{MR0406756}. Since they
require defining $\A^*$ as $\neg\rmin\A$, they may intend that the
axioms involving negation appear twice, once with $\rmin$, and once
with $\neg$. This is also suggested by Meyer~\cite{MR0419171}.  It
would have been interesting if Meyer and Routley had attempted a more
explicit axiomatization, for once the language includes the full range
of connectives $\land$, $\lor$, $\to$, $\circ$, $\rmin$, $\neg$, and
${}^*$, either primitively or by definition, the opportunity exists to
axiomatize classical relevant logic with Tarski's axioms.  We might
describe classical relevant logic as the system obtained from Tarski's
ten axioms \ra1--\ra{10}, suitable renotated using these translations:
\begin{align*}
  \A\lor\B&=\A+\B,\\
  \A\land\B&=\A\bp\B,\\
  \B\circ\A&=\A\rp\B,\\
  \A\to\B&=\min{\conv\A\rp\min\B},\\
  \neg\A&=\min\A,\\
  \rmin\A&=\con{\min\A},\\
  \A^*&=\conv\A.
\end{align*}
Be that as it may, they (and their readers, perhaps) feel that the
semantic characterization of CR* suffices. Certainly that is all we
need to observe that \lref{reflection} is not a theorem of CR${}^*$,
because it is invalidated in the CR* model structures $\gc\K_1$ and
$\gc\K_2$.  Indeed, all five normal relevant model structures
$\gc\K_1$--$\gc\K_5$ used in this paper satisfy the condition
$\R0\a\b\iff\a=\b$, and are therefore CR* model structures.
\section{Formulas in $\mathcal{L}_4$}
The 38 lemmas presented next establish membership in Tarski's
relevance logic of the formulas in Tables~\ref{table1} and
\ref{table2}.
\begin{lemma}\label{t6}
  $\A\lor\rmin\A$
\end{lemma}
\proof
\begin{align*}
  1.\quad&\fm\A\0\0\implies\fm\A\0\0 &&\text{Axiom}\\
  2.\quad&\implies\fm\A\0\0,\fm{\rmin\A}\0\0 &&\text{$|\rmin$}\\
  3.\quad&\implies\fm{\A\lor\rmin\A}\0\0 &&\text{$|\lor$}
\end{align*}
\endproof

\begin{lemma}\label{A1.}
  $\A\to\A$
\end{lemma}
\proof
  \begin{align*} 1.\quad&\fm\A\1\0\implies\fm\A\1\0
    &&\text{Axiom}\\
    2.\quad&\implies\fm{\A\to\A}\0\0 &&\text{$|{\to}$,~no~$\1$}
  \end{align*}
  \endproof

\begin{lemma}\label{A2.}
  $\A\land\B\to\A$
\end{lemma}
\proof
  \begin{align*}
    1.\quad&\fm\A\1\0,\fm\B\1\0\implies\fm\A\1\0
    &&\text{Axiom)}\\
    2.\quad&\fm{\A\land\B}\1\0\implies\fm\A\1\0
    &&\text{$\land|$}\\
    3.\quad&\implies\fm{(\A\land\B)\to\A}\0\0 &&\text{$|{\to}$, no
      $\1$} \end{align*}
  \endproof

\begin{lemma}\label{A3.}
$\A\land\B\to\B$
\end{lemma}
\proof
  \begin{align*}
    1.\quad&\fm\A\1\0,\fm\B\1\0\implies\fm\B\1\0
    &&\text{Axiom}\\
    2.\quad&\fm{\A\land\B}\1\0\implies\fm\B\1\0
    &&\text{1, $\land|$}\\
    3.\quad&\implies\fm{(\A\land\B)\to\B}\0\0 &&\text{2, $|{\to}$, no
      $\1$} \end{align*}
  \endproof

\begin{lemma}\label{A5.}
$\A\to\A\lor\B$
\end{lemma}
\proof
  \begin{align*}
    1.\quad&\fm\A\1\0\implies\fm\A\1\0,\fm\B\1\0	&&\text{Axiom}\\
    2.\quad&\fm\A\1\0\implies\fm{\A\lor\B}\1\0		&&\text{$|\lor$}\\
    3.\quad&\implies\fm{\A\to\A\lor\B}\0\0 &&\text{${\to}|$,~no~$\1$}
  \end{align*}
  \endproof

\begin{lemma}\label{A6.}
$\B\to\A\lor\B$
\end{lemma}
\proof
  \begin{align*}
    1.\quad&\fm\B\1\0\implies\fm\A\1\0,\fm\B\1\0
    &&\text{Axiom}\\
    2.\quad&\fm\B\1\0\implies\fm{\A\lor\B}\1\0
    &&\text{$|\lor$}\\
    3.\quad&\implies\fm{\B\to\A\lor\B}\0\0 &&\text{${\to}|$,~no~$\1$}
  \end{align*}
  \endproof

\begin{lemma}\label{comm1}
  $\B\lor\A\to\A\lor\B$
\end{lemma}
\proof
\begin{align*}
  1.\quad&\fm\A\1\0\implies\fm\A\1\0,\fm\B\1\0  &&\text{Axiom}\\
  2.\quad&\fm\A\1\0\implies\fm{\A\lor\B}1\0  &&\text{$|\lor$}\\
  3.\quad&\fm\B\1\0\implies\fm\A\1\0,\fm\B\1\0  &&\text{Axiom}\\
  4.\quad&\fm\B\1\0\implies\fm{\A\lor\B}1\0  &&\text{$|\lor$}\\
  5.\quad&\fm{\B\lor\A}\1\0\implies\fm{\A\lor\B}1\0  &&\text{2, 4, $\lor|$}\\
  6.\quad&\implies\fm{\B\lor\A\to\A\lor\B}\0\0 &&\text{$|{\to}$,
    no $\1$}
\end{align*}
\endproof
\begin{lemma}\label{comm2}
  $\B\land\A\to\A\land\B$
\end{lemma}
\proof
\begin{align*}
  1.\quad&\fm\A\1\0\implies\fm\A\1\0  &&\text{Axiom}\\
  2.\quad&\fm\B\1\0\implies\fm\B\1\0  &&\text{Axiom}\\
  3.\quad&\fm\A\1\0,\,\fm\B\1\0\implies\fm{\A\land\B}1\0  &&\text{$|\land$}\\
  4.\quad&\fm{\B\land\A}\1\0\implies\fm{\A\land\B}1\0  &&\text{$\land|$}\\
  5.\quad&\implies\fm{\B\land\A\to\A\land\B}\0\0  &&\text{$|{\to}$, no $\1$}
\end{align*}
\endproof
\begin{lemma}\label{assoc1}
  $(\A\land\B)\land\C\to\A\land(\B\land\C)$
\end{lemma}
\proof
  \begin{align*}
    1.\quad&\fm\A\1\0,\,\fm\B\1\0,\,\fm\C\1\0\implies\fm\B\1\0
    &&\text{Axiom}\\
    2.\quad&\fm\A\1\0,\,\fm\B\1\0,\,\fm\C\1\0\implies\fm\C\1\0
    &&\text{Axiom}\\
    3.\quad&\fm\A\1\0,\,\fm\B\1\0,\,\fm\C\1\0\implies\fm{\B\land\C}\1\0
    &&\text{$|\land$}\\
    4.\quad&\fm\A\1\0,\,\fm\B\1\0,\,\fm\C\1\0\implies\fm\A\1\0
    &&\text{Axiom}\\
    5.\quad&\fm\A\1\0,\,\fm\B\1\0,\,\fm\C\1\0\implies\fm{A\land(\B\land\C)}\1\0
    &&\text{$|\land$}\\
    6.\quad&\fm{\A\land\B}\1\0,\,\fm\C\1\0\implies\fm{A\land(\B\land\C)}\1\0
    &&\text{$\land|$}\\
    7.\quad&\fm{(\A\land\B)\land\C}\1\0\implies\fm{\A\land(\B\land\C)}\1\0
    &&\text{$\land|$}\\
    8.\quad&\implies\fm{(\A\land\B)\land\C\to\A\land(\B\land\C)}\0\0
    &&\text{$|{\to}$, no $\1$}
  \end{align*}
\endproof

\begin{lemma}\label{assoc2}
  $(\A\lor\B)\lor\C\to\A\lor(\B\lor\C)$
\end{lemma}
\proof
  \begin{align*}
    1.\quad&\fm\A\1\0\implies\fm\A\1\0,\,\fm\B\1\0,\,\fm\C\1\0
    &&\text{Axiom}\\
    2.\quad&\fm\B\1\0\implies\fm\A\1\0,\,\fm\B\1\0,\,\fm\C\1\0
    &&\text{Axiom}\\
    3.\quad&\fm{\A\lor\B}\1\0\implies\fm\A\1\0,\,\fm\B\1\0,\,\fm\C\1\0
    &&\text{$|\lor$}\\
    4.\quad&\fm\C\1\0\implies\fm\A\1\0,\,\fm\B\1\0,\,\fm\C\1\0
    &&\text{Axiom}\\
    5.\quad&\fm{(\A\lor\B)\lor\C}\1\0\implies\fm\A\1\0,\,\fm\B\1\0,\,\fm\C\1\0
    &&\text{$\lor|$}\\
    6.\quad&\fm{(\A\lor\B)\lor\C}\1\0\implies\fm\A\1\0,\,\fm{\B\lor\C}\1\0
    &&\text{$|\lor$}\\
    7.\quad&\fm{(\A\lor\B)\lor\C}\1\0\implies\fm{\A\lor(\B\lor\C)}\1\0
    &&\text{$|\lor$}\\
    8.\quad&\implies\fm{(\A\lor\B)\lor\C\to\A\lor(\B\lor\C)}\0\0
    &&\text{$|{\to}$, no $\1$}
  \end{align*}
\endproof

\begin{lemma}\label{A8.}
  $\A\land(\B\lor\C)\to(\A\land\B)\lor(\A\land\C)$
\end{lemma}
\proof
  \begin{align*}
    1.\quad&\fm\A\1\0\implies\fm\A\1\0&&\text{Axiom}\\
    2.\quad&\fm\B\1\0\implies\fm\B\1\0&&\text{Axiom}\\
    3.\quad&\fm\C\1\0\implies\fm\C\1\0&&\text{Axiom}\\
    4.\quad&\fm\A\1\0,\fm\B\1\0\implies\fm{\A\land\B}\1\0
    &&\text{1, 2, $|\land$}\\
    5.\quad&\fm\A\1\0,\fm\C\1\0\implies\fm{\A\land\C}\1\0
    &&\text{1, 3, $|\land$}\\
    6.\quad&\fm\A\1\0,\fm{\B\lor\C}\1\0\implies\fm{\A\land\B}\1\0,\fm{\A\land\C}\1\0
    &&\text{$\lor|$}\\
    7.\quad&\fm{\A\land(\B\lor\C)}\1\0\implies\fm{\A\land\B}\1\0,\fm{\A\land\C}\1\0
    &&\text{$\land|$}\\
    8.\quad&\fm{\A\land(\B\lor\C)}\1\0\implies\fm{(\A\land\B)\lor(\A\land\C)}\1\0
    &&\text{$|\lor$}\\
    9.\quad&\implies\fm{(\A\land(\B\lor\C))\to((\A\land\B)\lor(\A\land\C))}\0\0
    &&\text{$|{\to}$,~no~$\1$}
  \end{align*}
  \endproof

\begin{lemma}\label{*T9.}
$(\A\to\rmin\C)\land(\B\to\C)\to\rmin(\A\land\B)$
\end{lemma}
\proof
\begin{align*}
  1.\quad&\fm\A\0\1\implies\fm\A\0\1
  &&\text{Axiom}\\
  2.\quad&\fm\B\0\1\implies\fm\B\0\1
  &&\text{Axiom}\\
  3.\quad&\fm\C\0\0\implies\fm\C\0\0
  &&\text{Axiom}\\
  4.\quad&\fm{\rmin\C}\0\0,\,\fm\C\0\0\implies
  &&\text{$\rmin|$}\\
  5.\quad&\fm{\A\to\rmin\C}\1\0,\,\fm\C\0\0,\,\fm\A\0\1\implies
  &&\text{1, 4, ${\to}|$}\\
  6.\quad&\fm{\A\to\rmin\C}\1\0,\,\fm{\B\to\C}\1\0,\,\fm\A\0\1,\,\fm\B\0\1\implies
  &&\text{2, 5, ${\to}|$}\\
  7.\quad&\fm{\A\to\rmin\C}\1\0,\,\fm{\B\to\C}\1\0,\,\fm{\A\land\B}\0\1\implies
  &&\text{$\land|$}\\
  8.\quad&\fm{\A\to\rmin\C}\1\0,\,\fm{\B\to\C}\1\0\implies\fm{\rmin(\A\land\B)}\1\0
  &&\text{$|\rmin$}\\
  9.\quad&\fm{(\A\to\rmin\C)\land(\B\to\C)}\1\0\implies\fm{\rmin(\A\land\B)}\1\0
  &&\text{$\land|$}\\
  10.\quad&\implies\fm{((\A\to\rmin\C)\land(\B\to\C))\to\rmin(\A\land\B)}\0\0
  &&\text{$|{\to}$,~no~$\1$}
\end{align*}
\endproof

\begin{lemma}\label{T10.}
  $(\A\to\rmin\B)\land(\rmin\A\to\rmin\C)\to(\rmin\B\lor\rmin\C)$
\end{lemma}
\proof
\begin{align*}
  1.\quad&\fm\A\1\1\implies\fm\A\1\1
  &&\text{Axiom}\\
  2.\quad&\implies\fm\A\1\1,\fm{\rmin\A}\1\1
  &&\text{$|\rmin$}\\
  3.\quad&\fm\C\0\1\implies\fm\C\0\1
  &&\text{Axiom}\\
  4.\quad&\fm{\rmin\C}\1\0,\fm\C\0\1\implies
  &&\text{$\rmin|$}\\
  5.\quad&\fm{\rmin\A\to\rmin\C}\1\0,\fm\C\0\1\implies\fm\A\1\1
  &&\text{2, 4, ${\to}|$}\\
  6.\quad&\fm\B\0\1\implies\fm\B\0\1
  &&\text{Axiom}\\
  7.\quad&\fm{\rmin\B}\1\0,\fm\B\0\1\implies
  &&\text{$\rmin|$}\\
  8.\quad&\fm{\A\to\rmin\B}\1\0,\fm{\rmin\A\to\rmin\C}\1\0,\fm\B\0\1,\fm\C\0\1\implies
  &&\text{5, 7, ${\to}|$}\\
  9.\quad&\fm{\A\to\rmin\B}\1\0,\fm{\rmin\A\to\rmin\C}\1\0,\fm\B\1\0\implies\fm{\rmin\C}\1\0
  &&\text{$|\rmin$}\\
  10.\quad&\fm{\A\to\rmin\B}\1\0,\fm{\rmin\A\to\rmin\C}\1\0\implies\fm{\rmin\B}\1\0,\fm{\rmin\C}\1\0
  &&\text{$|\rmin$}\\
  11.\quad&\fm{\A\to\rmin\B}\1\0,\fm{\rmin\A\to\rmin\C}\1\0\implies\fm{\rmin\B\lor\rmin\C}\1\0
  &&\text{$|\lor$}\\
  12.\quad&\fm{(\A\to\rmin\B)\land(\rmin\A\to\rmin\C)}\1\0\implies\fm{\rmin\B\lor\rmin\C}\1\0
  &&\text{$\land|$}\\
  13.\quad&\implies\fm{((\A\to\rmin\B)\land(\rmin\A\to\rmin\C))\to(\rmin\B\lor\rmin\C)}\0\0
  &&\text{$|{\to}$,~no~$\1$}
\end{align*}
\endproof

\begin{lemma}\label{A9.}
  $\rmin\rmin\A\to\A$
\end{lemma}
\proof
\begin{align*}
  1.\quad&\fm\A\1\0\implies\fm\A\1\0&&\text{Axiom}\\
  2.\quad&\implies\fm\A\1\0,\fm{\rmin\A}\0\1&&\text{$|\rmin$}\\
  3.\quad&\fm{\rmin\rmin\A}\1\0\implies\fm\A\1\0&&\text{$\rmin|$}\\
  4.\quad&\implies\fm{\rmin\rmin\A\to\A}\0\0 &&\text{$|{\to}$, no
    $\1$}
\end{align*}
\endproof
\begin{lemma}\label{t3}
  $\A\to\rmin\rmin\A$
\end{lemma}
\proof
\begin{align*}
  1.\quad&\fm\A\1\0\implies\fm\A\1\0&&\text{Axiom}\\
  2.\quad&\fm\A\1\0,\,\fm{\rmin\A}\0\1\implies&&\text{$\rmin|$}\\
  3.\quad&\fm\A\1\0\implies\fm{\rmin\rmin\A}\1\0&&\text{$|\rmin$}\\
  4.\quad&\implies\fm{\A\to\rmin\rmin\A}\0\0&&\text{$|{\to}$,~no~$\1$}
\end{align*}
\endproof

\begin{lemma}\label{T2.}
$\rmin(\A\lor\B)\to(\rmin\A\land\rmin\B)$
\end{lemma}
\proof
\begin{align*}
  1.\quad&\fm\A\0\1\implies\fm\A\0\1
  &&\text{Axiom}\\
  2.\quad&\implies\fm{\rmin\A}\1\0,\fm\A\0\1
  &&\text{$|\rmin$}\\
  3.\quad&\fm\B\0\1\implies\fm\B\0\1
  &&\text{Axiom}\\
  4.\quad&\implies\fm{\rmin\B}\1\0,\fm\B\0\1
  &&\text{$|\rmin$}\\
  5.\quad&\implies\fm{\rmin\A\land\rmin\B}\1\0,\fm\A\0\1,\fm\B\0\1
  &&\text{2, 4, $|\land$}\\
  6.\quad&\implies\fm{\rmin\A\land\rmin\B}\1\0,\fm{\A\lor\B}\0\1
  &&\text{$|\lor$}\\
  7.\quad&\fm{\rmin(\A\lor\B)}\1\0\implies\fm{\rmin\A\land\rmin\B}\1\0
  &&\text{$\rmin|$}\\
  8.\quad&\implies\fm{\rmin(\A\lor\B)\to(\rmin\A\land\rmin\B)}\0\0
  &&\text{$|{\to}$,~no~$\1$}
\end{align*}
\endproof

\begin{lemma}\label{t4}
  $\rmin(\A\land\B)\to(\rmin\A\lor\rmin\B)$
\end{lemma}
\proof
\begin{align*}
  1.\quad&\fm\A\0\1\implies\fm\A\0\1&&\text{Axiom}\\
  2.\quad&\fm\B\0\1\implies\fm\B\0\1&&\text{Axiom}\\
  3.\quad&\fm\A\0\1,\,\fm\B\0\1\implies\fm{\A\land\B}\0\1
  &&\text{$|\land$}\\
  4.\quad&\fm\A\0\1,\,\fm\B\0\1,\,\fm{\rmin(\A\land\B)}\1\0\implies
  &&\text{$\rmin|$}\\
  5.\quad&\fm\A\0\1,\,\fm{\rmin(\A\land\B)}\1\0\implies\fm{\rmin\B}\1\0
  &&\text{$|\rmin$}\\
  6.\quad&\fm{\rmin(\A\land\B)}\1\0\implies\fm{\rmin\A}\1\0,\,\fm{\rmin\B}\1\0
  &&\text{$|\rmin$}\\
  7.\quad&\fm{\rmin(\A\land\B)}\1\0\implies\fm{\rmin\A\lor\rmin\B}\1\0
  &&\text{$|\lor$}\\
  8.\quad&\implies\fm{\rmin(\A\land\B)\to(\rmin\A\lor\rmin\B)}\0\0
  &&\text{$|{\to}$,~no~$\1$}\\
\end{align*}
\endproof
\begin{lemma}\label{t5}
  $(\rmin\A\land\rmin\B)\to\rmin(\A\lor\B)$
\end{lemma}
\proof
\begin{align*}
  1.\quad&\fm\A\0\1\implies\fm\A\0\1&&\text{Axiom}\\
  2.\quad&\fm\B\0\1\implies\fm\B\0\1&&\text{Axiom}\\
  3.\quad&\fm{\A\lor\B}\0\1\implies\fm\A\0\1,\,\fm\B\0\1
  &&\text{$\lor|$}\\
  4.\quad&\fm{\rmin\B}\1\0,\,\fm{\A\lor\B}\0\1\implies\fm\A\0\1
  &&\text{$\rmin|$}\\
  5.\quad&\fm{\rmin\A}\1\0,\,\fm{\rmin\B}\1\0,\,\fm{\A\lor\B}\0\1\implies
  &&\text{$\rmin|$}\\
  6.\quad&\fm{\rmin\A\land\rmin\B}\1\0,\,\fm{\A\lor\B}\0\1\implies
  &&\text{$\land|$}\\
  7.\quad&\fm{\rmin\A\land\rmin\B}\1\0\implies\fm{\rmin(\A\lor\B)}\1\0
  &&\text{$|\rmin$}\\
  8.\quad&\implies\fm{(\rmin\A\land\rmin\B)\to\rmin(\A\lor\B)}\0\0
  &&\text{$|{\to}$,~no~$\1$}
\end{align*}
\endproof

\begin{lemma}\label{t5a}
  $(\rmin\A\lor\rmin\B)\to\rmin(\A\land\B)$
\end{lemma}
\proof
\begin{align*}
  1.\quad&\fm\A\0\1\implies\fm\A\0\1&&\text{Axiom}\\
  3.\quad&\fm{\rmin\A}\1\0,\,\fm\A\0\1\implies
  &&\text{$\rmin|$}\\
  2.\quad&\fm\B\0\1\implies\fm\B\0\1&&\text{Axiom}\\
  4.\quad&\fm{\rmin\B}\1\0,\,\fm\B\0\1\implies
  &&\text{$\rmin|$}\\
  5.\quad&\fm{\rmin\A\lor\rmin\B}\1\0,\,\fm\A\0\1,\,\fm\B\0\1\implies
  &&\text{2, 4, $\lor|$}\\
  6.\quad&\fm{\rmin\A\lor\rmin\B}\1\0,\,\fm{\A\land\B}\0\1\implies
  &&\text{$\land|$}\\
  7.\quad&\fm{\rmin\A\lor\rmin\B}\1\0\implies\fm{\rmin(\A\land\B)}\1\0
  &&\text{$|\rmin$}\\
  8.\quad&\implies\fm{(\rmin\A\lor\rmin\B)\to\rmin(\A\land\B)}\0\0
  &&\text{$|{\to}$,~no~$\1$}
\end{align*}
\endproof

\begin{lemma}\label{T11.}
  $((\A\to\A)\to\B)\to\B$
\end{lemma}
\proof
\begin{align*}
  1.\quad&\fm\A\0\1\implies\fm\A\0\1
  &&\text{Axiom}\\
  2.\quad&\implies\fm{\A\to\A}\1\1
  &&\text{$|{\to}$,~no~$\0$}\\
  3.\quad&\fm\B\1\0\implies\fm\B\1\0
  &&\text{Axiom}\\
  4.\quad&\fm{(\A\to\A)\to\B}\1\0\implies\fm\B\1\0
  &&\text{${\to}|$}\\
  5.\quad&\implies\fm{((\A\to\A)\to\B)\to\B}\0\0
  &&\text{$|{\to}$,~no~$\1$}
\end{align*}
\endproof

\begin{lemma}\label{A4.}
$(\A\to\B)\land(\A\to\C)\to(\A\to(\B\land\C))$
\end{lemma}
\proof
  \begin{align*} 1.\quad&\fm\A\2\1\implies\fm\A\2\1
    &&\text{Axiom}\\
    2.\quad&\fm\B\2\0\implies\fm\B\2\0
    &&\text{Axiom}\\
    3.\quad&\fm{\A\to\B}\1\0,\fm\A\2\1\implies\fm\B\2\0
    &&\text{1, 2, ${\to}|$}\\
    4.\quad&\fm\C\2\0\implies\fm\C\2\0
    &&\text{Axiom}\\
    5.\quad&\fm{\A\to\C}\1\0,\fm\A\2\1\implies\fm\C\2\0
    &&\text{1, 4, ${\to}|$}\\
    6.\quad&\fm{\A\to\B}\1\0,\fm{\A\to\C}\1\0,\fm\A\2\1\implies\fm{\B\land\C}\2\0
    &&\text{3, 5, $|\land$}\\
    7.\quad&\fm{(\A\to\B)\land(\A\to\C)}\1\0,\fm\A\2\1\implies\fm{\B\land\C}\2\0
    &&\text{$\land|$}\\
    8.\quad&\fm{(\A\to\B)\land(\A\to\C)}\1\0\implies\fm{\A\to(\B\land\C)}\1\0
    &&\text{$|{\to}$,~no~$\2$}\\
    9.\quad&\implies\fm{((\A\to\B)\land(\A\to\C))\to(\A\to(\B\land\C))}\0\0
    &&\text{$|{\to}$,~no~$\1$}
  \end{align*}
  \endproof

\begin{lemma}\label{A7.}
$(\A\to\C)\land(\B\to\C)\to((\A\lor\B)\to\C)$
\end{lemma}
\proof
  \begin{align*} 1.\quad&\fm\A\2\1\implies\fm\A\2\1
    &&\text{Axiom}\\
    2.\quad&\fm\C\2\0\implies\fm\C\2\0
    &&\text{Axiom}\\
    3.\quad&\fm{\A\to\C}\1\0,\fm\A\2\1\implies\fm\C\2\0
    &&\text{${\to}|$}\\
    4.\quad&\fm\B\2\1\implies\fm\B\2\1
    &&\text{Axiom}\\
    5.\quad&\fm\C\2\0\implies\fm\C\2\0
    &&\text{Axiom}\\
    6.\quad&\fm{\B\to\C}\1\0,\fm\B\2\1\implies\fm\C\2\0
    &&\text{${\to}|$}\\
    7.\quad&\fm{\A\to\C}\1\0,\fm{\B\to\C}\1\0,\fm{\A\lor\B}\2\1\implies\fm\C\2\0
    &&\text{3, 6, $\lor|$}\\
    8.\quad&\fm{(\A\to\C)\land(\B\to\C)}\1\0,\fm{\A\lor\B}\2\1\implies\fm\C\2\0
    &&\text{$\land|$}\\
    9.\quad&\fm{(\A\to\C)\land(\B\to\C)}\1\0\implies\fm{(\A\lor\B)\to\C}\1\0
    &&\text{${\to}|$,~no~$\2$}\\
    10.\quad&\implies\fm{((\A\to\C)\land(\B\to\C))\to((\A\lor\B)\to\C)}\0\0
    &&\text{${\to}|$,~no~$\1$}
  \end{align*}
  \endproof

\begin{lemma}\label{t11}
  $(\A\to\B)\land(\C\to\D)\to(\A\land\C\to\B\land\D)$
\end{lemma}
\proof
  \begin{align*}
    1.\quad&\fm\A\2\1,\fm\C\2\1\implies\fm\C\2\1&&\text{Axiom}\\
    2.\quad&\fm{\A\land\C}\2\1\implies\fm\C\2\1&&\text{$\land|$}\\
    3.\quad&\fm\B\2\0,\fm\D\2\0\implies\fm\B\2\0&&\text{Axiom}\\
    4.\quad&\fm\B\2\0,\fm\D\2\0\implies\fm\D\2\0&&\text{Axiom}\\
    5.\quad&\fm\B\2\0,\fm\D\2\0\implies\fm{\B\land\D}\2\0&&\text{$|\land$}\\
    6.\quad&\fm\B\2\0,\fm{\C\to\D}\1\0,\fm{\A\land\C}\2\1
    \implies\fm{\B\land\D}\2\0&&\text{2, 5, ${\to}|$}\\
    7.\quad&\fm\A\2\1,\fm\C\2\1\implies\fm\A\2\1&&\text{Axiom}\\
    8.\quad&\fm{\A\land\C}\2\1\implies\fm\A\2\1&&\text{$\land|$}\\
    9.\quad&\fm{\A\to\B}\1\0,\fm{\C\to\D}\1\0,\fm{\A\land\C}\2\1
    \implies\fm{\B\land\D}\2\0&&\text{6, 8, ${\to}|$}\\
    10.\quad&\fm{\A\to\B}\1\0,\fm{\C\to\D}\1\0
    \implies\fm{\A\land\C\to\B\land\D}\1\0&&\text{$|{\to}$,~no~$\2$}\\
    11.\quad&\fm{(\A\to\B)\land(\C\to\D)}\1\0
    \implies\fm{\A\land\C\to\B\land\D}\1\0&&\text{$\land|$}\\
    12.\quad&\implies\fm{(\A\to\B)\land(\C\to\D)\to(\A\land\C\to\B\land\D)}\0\0
    &&\text{$|{\to}$,~no~$\1$}
  \end{align*}
  \endproof

\begin{lemma}\label{T6.}
  $(\A\to\B)\land(\C\to\D)\to(\A\lor\C\to\B\lor\D)$
\end{lemma}
\proof
\begin{align*}
  1.\quad&\fm\A\2\1\implies\fm\A\2\1
  &&\text{Axiom}\\
  2.\quad&\fm\B\2\0\implies\fm\B\2\0,\fm\D\2\0
  &&\text{Axiom}\\
  3.\quad&\fm\B\2\0\implies\fm{\B\lor\D}\2\0
  &&\text{$|\lor$}\\
  4.\quad&\fm{\A\to\B}\1\0,\fm\A\2\1\implies\fm{\B\lor\D}\2\0
  &&\text{1, 3, ${\to}|$}\\
  5.\quad&\fm\C\2\1\implies\fm\C\2\1
  &&\text{Axiom}\\
  6.\quad&\fm\D\2\0\implies\fm\B\2\0,\fm\D\2\0
  &&\text{Axiom}\\
  7.\quad&\fm\D\2\0\implies\fm{\B\lor\D}\2\0
  &&\text{$|\lor$}\\
  8.\quad&\fm{\C\to\D}\1\0,\fm\C\2\1\implies\fm{\B\lor\D}\2\0
  &&\text{5, 7, ${\to}|$}\\
  9.\quad&\fm{\A\to\B}\1\0,\fm{\C\to\D}\1\0,\fm{\A\lor\C}\2\1\implies\fm{\B\lor\D}\2\0
  &&\text{4, 8, $\lor|$}\\
  10.\quad&\fm{\A\to\B}\1\0,\fm{\C\to\D}\1\0\implies\fm{\A\lor\C\to\B\lor\D}\1\0
  &&\text{$|{\to}$,~no~$\2$}\\
  11.\quad&\fm{(\A\to\B)\land(\C\to\D)}\1\0\implies\fm{\A\lor\C\to\B\lor\D}\1\0
  &&\text{$\land|$}\\
  12.\quad&\implies\fm{(\A\to\B)\land(\C\to\D)\to(\A\lor\C\to\B\lor\D)}\0\0
  &&\text{$|{\to}$,~no~$\1$}
\end{align*}
\endproof

\begin{lemma}\label{T8.}
$(\A\to\B)\lor(\C\to\D)\to((\A\land\C)\to(\B\lor\D))$
\end{lemma}
\proof
\begin{align*}
  1.\quad&\fm\A\2\1\implies\fm\A2\1
  &&\text{Axiom}\\
  2.\quad&\fm\B\2\0,\fm\C\2\1\implies\fm\B\2\0,\fm\D\2\0
  &&\text{Axiom}\\
  3.\quad&\fm{\A\to\B}\1\0,\fm\A\2\1,\fm\C\2\1\implies\fm\B\2\0,\fm\D\2\0
  &&\text{${\to}|$}\\
  4.\quad&\fm{\A\to\B}\1\0,\fm\A\2\1,\fm\C\2\1\implies\fm{\B\lor\D}\2\0
  &&\text{$|\lor$}\\
  5.\quad&\fm{\A\to\B}\1\0,\fm{\A\land\C}\2\1\implies\fm{\B\lor\D}\2\0
  &&\text{$\land|$}\\
  6.\quad&\fm{\A\to\B}\1\0\implies\fm{(\A\land\C)\to(\B\lor\D)}\1\0
  &&\text{no $\2$}\\
  7.\quad&\fm\C\2\1\implies\fm\C\2\1
  &&\text{Axiom}\\
  8.\quad&\fm\D\2\0,\fm\A\2\1\implies\fm\B\2\0,\fm\D\2\0
  &&\text{Axiom}\\
  9.\quad&\fm{\C\to\D}\1\0,\fm\A\2\1,\fm\C\2\1\implies\fm\B\2\0,\fm\D\2\0
  &&\text{${\to}|$}\\
  10.\quad&\fm{\C\to\D}\1\0,\fm\A\2\1,\fm\C\2\1\implies\fm{\B\lor\D}\2\0
  &&\text{$|\lor$}\\
  11.\quad&\fm{\C\to\D}\1\0,\fm{\A\land\C}\2\1\implies\fm{\B\lor\D}\2\0
  &&\text{$\land|$}\\
  12.\quad&\fm{\C\to\D}\1\0\implies\fm{(\A\land\C)\to(\B\lor\D)}\1\0
  &&\text{no $\2$}\\
  13.\quad&\fm{(\A\to\B)\lor(\C\to\D)}\1\0\implies\fm{(\A\land\C)\to(\B\lor\D)}\1\0
  &&\text{6, 12, $\lor|$}\\
  14.\quad&\implies\fm{((\A\to\B)\lor(\C\to\D))\to((\A\land\C)\to(\B\lor\D))}\0\0
  &&\text{$|{\to}$,~no~$\1$}
\end{align*}
\endproof

\begin{lemma}\label{t7}
$\A\to(\rmin\B\to\rmin(\A\to\B))$
\end{lemma}
\proof
\begin{align*}
  1.\quad&\fm\A\1\0\implies\fm\A\1\0&&\text{Axiom}\\
  2.\quad&\fm\B\1\2\implies\fm\B\1\2&&\text{Axiom}\\
  3.\quad&\fm\A\1\0,\,\fm{\A\to\B}\0\2\implies\fm\B\1\2
  &&\text{${\to}|$}\\
  4.\quad&\fm\A\1\0\implies\fm{\rmin(\A\to\B)}\2\0,\,\fm\B\1\2
  &&\text{$|\rmin$}\\
  5.\quad&\fm\A\1\0,\,\fm{\rmin\B}\2\1\implies\fm{\rmin(\A\to\B)}\2\0
  &&\text{$\rmin|$}\\
  6.\quad&\fm\A\1\0\implies\fm{\rmin\B\to\rmin(\A\to\B)}\1\0
  &&\text{$|{\to}$,~no~$\2$}\\
  7.\quad&\implies\fm{\A\to(\rmin\B\to\rmin(\A\to\B))}\0\0
  &&\text{$|{\to}$,~no~$\1$}
\end{align*}
\endproof

\begin{lemma}\label{t9}
  $\A\to(\B\to\rmin(\A\to\rmin\B))$
\end{lemma}
\proof
  \begin{align*}
    1.\quad&\fm\A\1\0\implies\fm\A\1\0&&\text{Axiom}\\
    2.\quad&\fm\B\2\1\implies\fm\B\2\1&&\text{Axiom}\\
    3.\quad&\fm\B\2\1,\,\fm{\rmin\B}\1\2\implies&&\text{$\rmin|$}\\
    4.\quad&\fm\A\1\0,\,\fm\B\2\1,\,\fm{\A\to\rmin\B}\0\2\implies
    &&\text{1, 3, ${\to}|$}\\
    5.\quad&\fm\A\1\0,\,\fm\B\2\1\implies\fm{\rmin(\A\to\rmin\B)}\2\0
    &&\text{$|\rmin$}\\
    6.\quad&\fm\A\1\0\implies\fm{\B\to\rmin(\A\to\rmin\B)}\1\0
    &&\text{$|{\to}$,~no~$\2$}\\
    7.\quad&\implies\fm{\A\to(\B\to\rmin(\A\to\rmin\B))}\0\0
    &&\text{$|{\to}$,~no~$\1$}
  \end{align*}
  \endproof

\begin{lemma}\label{t13}
  $\A\to((\rmin\B\to\rmin\A)\to\B)$
\end{lemma}
\proof
  \begin{align*}
    1.\quad&\fm\B\2\0\implies\fm\B\2\0&&\text{Axiom}\\
    2.\quad&\implies\fm\B\2\0,\,\fm{\rmin\B}\0\2&&\text{$|\rmin$}\\
    3.\quad&\fm\A\1\0\implies\fm\A\1\0&&\text{Axiom}\\
    4.\quad&\fm\A\1\0,\,\fm{\rmin\A}\0\1\implies&&\text{$\rmin|$}\\
    5.\quad&\fm\A\1\0,\,\fm{\rmin\B\to\rmin\A}\2\1\implies\fm\B\2\0
    &&\text{2, 4, ${\to}|$}\\
    6.\quad&\fm\A\1\0\implies\fm{(\rmin\B\to\rmin\A)\to\B}\1\0
    &&\text{$|{\to}$,~no~$\2$}\\
    7.\quad&\implies\fm{\A\to((\rmin\B\to\rmin\A)\to\B)}\0\0
    &&\text{$|{\to}$,~no~$\1$}
  \end{align*}
\endproof

\begin{lemma}\label{t14}
  $\A\to((\B\to\rmin\A)\to\rmin\B)$
\end{lemma}
\proof
  \begin{align*}
    1.\quad&\fm\A\1\0\implies\fm\A\1\0&&\text{Axiom}\\
    2.\quad&\fm\A\1\0,\,\fm{\rmin\A}\0\1\implies&&\text{$\rmin|$}\\
    3.\quad&\fm\B\0\2\implies\fm\B\0\2&&\text{Axiom}\\
    4.\quad&\fm\A\1\0,\,\fm{\B\to\rmin\A}\2\1,\,\fm\B\0\2\implies
    &&\text{${\to}|$}\\
    5.\quad&\fm\A\1\0,\,\fm{\B\to\rmin\A}\2\1\implies\fm{\rmin\B}\2\0
    &&\text{$|\rmin$}\\
    6.\quad&\fm\A\1\0\implies\fm{(\B\to\rmin\A)\to\rmin\B}\1\0
    &&\text{$|{\to}$,~no~$\2$}\\
    7.\quad&\implies\fm{\A\to((\B\to\rmin\A)\to\rmin\B)}\0\0
    &&\text{$|{\to}$,~no~$\1$}
  \end{align*}
  \endproof 

\begin{lemma}\label{T12.}
  $\rmin((\A\to\B)\to\rmin\A)\to\B$
\end{lemma}
\proof
\begin{align*}
  1.\quad&\fm\A\1\2\implies\fm\A\1\2
  &&\text{Axiom}\\
  2.\quad&\implies\fm\A\1\2,\,\fm{\rmin\A}\2\1
  &&\text{$\rmin|$}\\
  3.\quad&\fm\B\1\0\implies\fm{\rmin\A}\2\1,\,\fm\B\1\0
  &&\text{Axiom}\\
  4.\quad&\fm{\A\to\B}\2\0\implies\fm{\rmin\A}\2\1,\,\fm\B\1\0
  &&\text{${\to}|$}\\
  5.\quad&\implies\fm{(\A\to\B)\to\rmin\A}\0\1,\,\fm\B\1\0
  &&\text{$|{\to}$,~no~$\2$}\\
  6.\quad&\fm{\rmin((\A\to\B)\to\rmin\A)}\1\0\implies\fm\B\1\0
  &&\text{$|\rmin$}\\
  7.\quad&\implies\fm{\rmin((\A\to\B)\to\rmin\A)\to\B}\0\0
  &&\text{$|{\to}$,~no~$\1$}
\end{align*}
\endproof

\begin{lemma}\label{*T15.}
  $\rmin\A\to((\B\to\A)\to\rmin\B))$
\end{lemma}
\proof
\begin{align*}
  1.\quad&\fm\B\0\2\implies\fm\B\0\2
  &&\text{Axiom}\\
  2.\quad&\fm\A\0\1\implies\fm\A\0\1
  &&\text{Axiom}\\
  3.\quad&\fm{\B\to\A}\2\1,\,\fm\B\0\2\implies\fm\A\0\1
  &&\text{${\to}|$}\\
  4.\quad&\fm{\rmin\A}\1\0,\,\fm{\B\to\A}\2\1,\,\fm\B\0\2\implies
  &&\text{$\rmin|$}\\
  5.\quad&\fm{\rmin\A}\1\0,\,\fm{\B\to\A}\2\1\implies\fm{\rmin\B}\2\0
  &&\text{$|\rmin$}\\
  6.\quad&\fm{\rmin\A}\1\0\implies\fm{(\B\to\A)\to\rmin\B}\1\0
  &&\text{$|{\to}$,~no~$\2$}\\
  7.\quad&\implies\fm{\rmin\A\to((\B\to\A)\to\rmin\B)}\0\0
  &&\text{$|{\to}$,~no~$\1$}
\end{align*}
\endproof

\begin{lemma}\label{reflection}
  \begin{equation*}
    (\A\circ\B)\land\C\to((\A\land\rmin\D)\circ\B)\lor(\A\circ(\B\land(\D\circ\C)))
  \end{equation*}
  \begin{align*}
    \rmin(\A\to\rmin\B)\land\C&\to\rmin((\A\land\rmin\D)\to\rmin\B)
    \lor\rmin(\A\to\rmin(\B\land\rmin(\D\to\rmin\C)))
  \end{align*}
\end{lemma}
\proof
  \begin{align*}
    1.\quad &\fm\A\2\0\implies\fm\A\2\0
    &&\text{Axiom}\\
    2.\quad &\fm\D\0\2\implies\fm\D\0\2
    &&\text{Axiom}\\
    3.\quad &\implies\fm{\rmin\D}\2\0,\,\fm\D\0\2
    &&\text{$|\rmin$}\\
    4.\quad &\fm\A\2\0\implies\fm{\A\land\rmin\D}\2\0,\,\fm\D\0\2
    &&\text{1, 3, $|\land$}\\
    5.\quad &\fm\B\1\2\implies\fm\B\1\2
    &&\text{Axiom}\\
    6.\quad &\fm\B\1\2,\,\fm{\rmin\B}\2\1\implies
    &&\text{$\rmin|$}\\
    7.\quad&\fm\A\2\0,\,\fm\B\1\2,\,\fm{(\A\land\rmin\D)\to\rmin\B}\0\1
    \implies\fm\D\0\2
    &&\text{4, 6, ${\to}|$}\\
    8.\quad&\fm\C\1\0\implies\fm\C\1\0
    &&\text{Axiom}\\
    9.\quad&\fm\C\1\0,\,\fm{\rmin\C}\0\1\implies
    &&\text{$\rmin|$}\\
    10.\quad&\fm\C\1\0,\,\fm\A\2\0,\,\fm\B\1\2,\,\fm{(\A\land\rmin\D)\to\rmin\B}\0\1,\,
    \fm{\D\to\rmin\C}\2\1\implies
    &&\text{7, 9, ${\to}|$}\\
    11.\quad&\fm\C\1\0,\,\fm\A\2\0,\,\fm\B\1\2,\,\fm{(\A\land\rmin\D)\to\rmin\B}\0\1
    \implies\fm{\rmin(\D\to\rmin\C)}\1\2
    &&\text{$|\rmin$}\\
    12.\quad&\fm\C\1\0,\,\fm\A\2\0,\,\fm\B\1\2,\,\fm{(\A\land\rmin\D)\to\rmin\B}\0\1
    \implies\fm{\B\land\rmin(\D\to\rmin\C)}\1\2
    &&\text{5, 11, $|\land$}\\
    13.\quad&\fm\C\1\0,\,\fm\A\2\0,\,\fm\B\1\2,\,\fm{(\A\land\rmin\D)\to\rmin\B}\0\1,\,
    \fm{\rmin(\B\land\rmin(\D\to\rmin\C))}\2\1\implies
    &&\text{$|\rmin$}\\
    14.\quad&\fm\C\1\0,\,\fm\A\2\0,\,\fm\B\1\2,\,\fm{(\A\land\rmin\D)\to\rmin\B}\0\1,\,
    \fm{\A\to\rmin(\B\land\rmin(\D\to\rmin\C))}\0\1\implies
    &&\text{1, 13, ${\to}|$}\\
    15.\quad&\fm\C\1\0,\,\fm\A\2\0,\,\fm{(\A\land\rmin\D)\to\rmin\B}\0\1,\,
    \fm{\A\to\rmin(\B\land\rmin(\D\to\rmin\C))}\0\1
    \implies\fm{\rmin\B}\2\1
    &&\text{$|\rmin$}\\
    16.\quad&\fm\C\1\0,\, \fm{(\A\land\rmin\D)\to\rmin\B}\0\1,\,
    \fm{\A\to\rmin(\B\land\rmin(\D\to\rmin\C))}\0\1
    \implies\fm{\A\to\rmin\B}\0\1
    &&\text{$|{\to}$, no $\2$}\\
    17.\quad&\fm\C\1\0,\,\fm{\A\to\rmin(\B\land\rmin(\D\to\rmin\C))}\0\1
    \implies\fm{\A\to\rmin\B}\0\1,\,
    \fm{\rmin((\A\land\rmin\D)\to\rmin\B)}\1\0
    &&\text{$|\rmin$}\\
    18.\quad&\fm\C\1\0 \implies\fm{\A\to\rmin\B}\0\1,\,
    \fm{\rmin((\A\land\rmin\D)\to\rmin\B)}\1\0,\,
    \fm{\rmin(\A\to\rmin(\B\land\rmin(\D\to\rmin\C)))}\1\0
    &&\text{$|\rmin$}\\
    19.\quad&\fm{\rmin(\A\to\rmin\B)}\1\0,\,\fm\C\1\0
    \implies\fm{\rmin((\A\land\rmin\D)\to\rmin\B)}\1\0,\,
    \fm{\rmin(\A\to\rmin(\B\land\rmin(\D\to\rmin\C)))}\1\0
    &&\text{$\rmin|$}\\
    20.\quad&\fm{\rmin(\A\to\rmin\B)\land\C}\1\0
    \implies\fm{\rmin((\A\land\rmin\D)\to\rmin\B)}\1\0,\,
    \fm{\rmin(\A\to\rmin(\B\land\rmin(\D\to\rmin\C)))}\1\0
    &&\text{$\land|$}\\
    21.\quad&\fm{\rmin(\A\to\rmin\B)\land\C}\1\0
    \implies\fm{\rmin((\A\land\rmin\D)\to\rmin\B)
      \lor\rmin(\A\to\rmin(\B\land\rmin(\D\to\rmin\C)))}\1\0
    &&\text{$|\lor$}\\
    22.\quad&\implies\fm{\rmin(\A\to\rmin\B)\land\C\to\rmin((\A\land\rmin\D)\to\rmin\B)
      \lor\rmin(\A\to\rmin(\B\land\rmin(\D\to\rmin\C)))}\0\0
    &&\text{$|{\to}$, no $\1$}
  \end{align*}
  \endproof

\begin{lemma}\label{t??}
  $(\A\to\B)\land\(\C\circ\D\)\to ((\C\land\B)\circ\D) \lor
  (\C\circ(\D\land\rmin\A))$
  \begin{equation*}
    (\A\to\B)\land\rmin\(\C\to\rmin\D\)\to
    \rmin((\C\land\B)\to\rmin\D)
    \lor\rmin(\C\to\rmin(\D\land\rmin\A))
  \end{equation*}
\end{lemma}
\proof
  \begin{align*}
    1.\quad&\fm{\A}\2\1\implies\fm{\A}\2\1
    &&\text{Axiom}\\
    2.\quad&\fm{\B}\2\0\implies\fm{\B}\2\0
    &&\text{Axiom}\\
    3.\quad&\fm{\A\to\B}\1\0,\,\fm{\A}\2\1\implies\fm{\B}\2\0
    &&\text{${\to}|$}\\
    4.\quad&\fm{\A\to\B}\1\0\implies\fm{\rmin\A}\1\2,\,\ul{\fm{\B}\2\0}
    &&\text{$|\rmin$}\\
    5.\quad&\fm\C\2\0\implies\ul{\fm\C\2\0}
    &&\text{Axiom}\\
    6.\quad&\fm{\A\to\B}\1\0,\,\fm\C\2\0
    \implies\ull{\fm{\rmin\A}\1\2},\,\ul{\fm{\C\land\B}\2\0}
    &&\text{$|\land$}\\
    7.\quad&\fm\D\1\2\implies\ull{\fm{\D}\1\2}
    &&\text{Axiom}\\
    8.\quad&\fm{\A\to\B}\1\0,\,\fm\C\2\0,\,\ul{\fm\D\1\2}
    \implies\ull{\fm{\D\land\rmin\A}\1\2},\,\fm{\C\land\B}\2\0
    &&\text{$|\land$}\\
    9.\quad&\fm{\A\to\B}\1\0,\,\fm\C\2\0
    \implies\ull{\fm{\D\land\rmin\A}\1\2},\,\fm{\C\land\B}\2\0
    ,\,\ul{\fm{\rmin\D}\2\1}
    &&\text{$|\rmin$}\\
    10.\quad&\fm{\A\to\B}\1\0,\,
    \ull{\fm{\rmin(\D\land\rmin\A)}\2\1},\,\fm\C\2\0
    \implies\ul{\fm{\C\land\B}\2\0},\,\fm{\rmin\D}\2\1
    &&\text{$\rmin|$}\\
    11.\quad&\ul{\fm{\rmin\D}\2\1} \implies{\fm{\rmin\D}\2\1}
    &&\text{Axiom}\\
    12.\quad&\fm{\A\to\B}\1\0,\,\ul{\fm{(\C\land\B)\to\rmin\D}\0\1},\,
    \ull{\fm{\rmin(\D\land\rmin\A)}\2\1},\,\fm\C\2\0
    \implies\fm{\rmin\D}\2\1
    &&\text{${\to}|$}\\
    13.\quad&\fm\C\2\0 \implies\ull{\fm\C\2\0}
    &&\text{Axiom}\\
    14.\quad&\fm{\A\to\B}\1\0,\,\fm{(\C\land\B)\to\rmin\D}\0\1
    ,\,\ull{\fm{\C\to\rmin(\D\land\rmin\A)}\0\1},\,\ul{\fm\C\2\0}
    \implies\ul{\fm{\rmin\D}\2\1}
    &&\text{${\to}|$}\\
    15.\quad&\fm{\A\to\B}\1\0,\,\fm{(\C\land\B)\to\rmin\D}\0\1
    ,\,\fm{\C\to\rmin(\D\land\rmin\A)}\0\1
    \implies\ul{\fm{\C\to\rmin\D}\0\1}
    &&\text{$|{\to}$,~no~$\2$}\\
    16.\quad&\fm{\A\to\B}\1\0,\,\fm{\C\to\rmin(\D\land\rmin\A)}\0\1
    \implies\fm{\rmin((\C\land\B)\to\rmin\D)}\1\0,\,\fm{\C\to\rmin\D}\0\1
    &&\text{$|\rmin$}\\
    17.\quad&\fm{\A\to\B}\1\0
    \implies\fm{\rmin((\C\land\B)\to\rmin\D)}\1\0
    ,\,\fm{\rmin(\C\to\rmin(\D\land\rmin\A))}\1\0,\,\fm{\C\to\rmin\D}\0\1
    &&\text{$|\rmin$}\\
    18.\quad&\fm{\A\to\B}\1\0,\,\fm{\rmin\(\C\to\rmin\D\)}\1\0
    \implies\fm{\rmin((\C\land\B)\to\rmin\D)}\1\0,\,
    \fm{\rmin(\C\to\rmin(\D\land\rmin\A))}\1\0
    &&\text{$\rmin|$}\\
    19.\quad&\fm{(\A\to\B)\land\rmin\(\C\to\rmin\D\)}\1\0
    \implies\fm{\rmin((\C\land\B)\to\rmin\D)}\1\0,\,
    \fm{\rmin(\C\to\rmin(\D\land\rmin\A))}\1\0
    &&\text{$\land|$}\\
    20.\quad&\fm{(\A\to\B)\land\rmin\(\C\to\rmin\D\)}\1\0
    \implies\fm{\rmin((\C\land\B)\to\rmin\D)
      \lor\rmin(\C\to\rmin(\D\land\rmin\A))}\1\0
    &&\text{$|\lor$}\\
    21.\quad&\implies\fm{(\A\to\B)\land\rmin\(\C\to\rmin\D\)\to
      \rmin((\C\land\B)\to\rmin\D)
      \lor\rmin(\C\to\rmin(\D\land\rmin\A))}\0\0
    &&\text{$|{\to}$,~no~$\1$}
  \end{align*}
  \endproof

\begin{lemma}\label{prefixingA}
  $(\A\to\B)\to((\C\to\A)\to(\C\to\B))$
\end{lemma}
\proof
  \begin{align*}
    1.\quad&\fm\C\3\2\implies\fm\C\3\2&&\text{Axiom}\\
    2.\quad&\fm\A\3\1\implies\fm\A\3\1&&\text{Axiom}\\
    3.\quad&\fm{\C\to\A}\2\1,\fm\C\3\2\implies\fm\A\3\1&&\text{${\to}|$}\\
    4.\quad&\fm\B\3\0\implies\fm\B\3\0&&\text{Axiom}\\
    5.\quad&\fm{\A\to\B}\1\0,\fm{\C\to\A}\2\1,\fm\C\3\2\implies\fm\B\3\0
    &&\text{${\to}|$}\\
    6.\quad&\fm{\A\to\B}\1\0,\fm{\C\to\A}\2\1\implies\fm{\C\to\B}\2\0
    &&\text{$|{\to}$,~no~$\3$}\\
    7.\quad&\fm{\A\to\B}\1\0\implies\fm{(\C\to\A)\to(\C\to\B)}\1\0
    &&\text{$|{\to}$,~no~$\2$}\\
    8.\quad&\implies\fm{(\A\to\B)\to((\C\to\A)\to(\C\to\B))}\0\0
    &&\text{$|{\to}$,~no~$\1$}
  \end{align*}
  \endproof

\begin{lemma}\label{t10}
  $(\B\to(\C\to\A))\to(\rmin(\B\to\rmin\C)\to\A)$
\end{lemma}
\proof
  \begin{align*}
    1.\quad&\fm\B\3\1\implies\fm\B\3\1&&\text{Axiom}\\
    2.\quad&\fm\C\2\3\implies\fm\C\2\3&&\text{Axiom}\\
    3.\quad&\fm\A\2\0\implies\fm\A\2\0&&\text{Axiom}\\
    4.\quad&\fm{\C\to\A}\3\0,\,\fm\C\2\3\implies\fm\A\2\0&&\text{${\to}|$}\\
    5.\quad&\fm{\B\to(\C\to\A)}\1\0,\fm\B\3\1,\fm\C\2\3\implies\fm\A\2\0
    &&\text{1, 4, ${\to}|$}\\
    6.\quad&\fm{\B\to(\C\to\A)}\1\0,\fm\B\3\1\implies\fm\A\2\0,\fm{\rmin\C}\3\2
    &&\text{$|\rmin$}\\
    7.\quad&\fm{\B\to(\C\to\A)}\1\0\implies\fm\A\2\0,\,\fm{\B\to\rmin\C}\1\2
    &&\text{$|{\to}$,~no~$\3$}\\
    8.\quad&\fm{\B\to(\C\to\A)}\1\0,\,\fm{\rmin(\B\to\rmin\C)}\2\1\implies\fm\A\2\0
    &&\text{$\rmin|$}\\
    9.\quad&\fm{\B\to(\C\to\A)}\1\0\implies\fm{(\rmin(\B\to\rmin\C)\to\A)}\1\0
    &&\text{$|{\to}$,~no~$\2$}\\
    10.\quad&\implies\fm{(\B\to(\C\to\A))\to(\rmin(\B\to\rmin\C)\to\A)}\0\0
    &&\text{$|{\to}$,~no~$\1$}
  \end{align*}
  \endproof

\begin{lemma}\label{T19.}
  $(\rmin(\A\to\rmin\B)\to\C)\to(\A\to(\B\to\C))$
\end{lemma}
\proof
  \begin{align*}
    1.\quad&\fm\B\3\2\implies\fm\B\3\2
    &&\text{Axiom}\\
    2.\quad&\fm\A\2\1\implies\fm\A\2\1
    &&\text{Axiom}\\
    3.\quad&\fm\B\3\2,\fm{\rmin\B}\2\3\implies
    &&\text{$\rmin|$}\\
    4.\quad&\fm\A\2\1,\fm\B\3\2,\fm{\A\to\rmin\B}\1\3\implies
    &&\text{${\to}|$}\\
    5.\quad&\fm\A\2\1,\fm\B\3\2\implies\fm{\rmin(\A\to\rmin\B)}\3\1
    &&\text{$|\rmin$}\\
    6.\quad&\fm\C\3\0\implies\fm\C\3\0
    &&\text{Axiom}\\
    7.\quad&\fm{\rmin(\A\to\rmin\B)\to\C}\1\0,\,\fm\A\2\1,\fm\B\3\2\implies\fm\C\3\0
    &&\text{${\to}|$}\\
    8.\quad&\fm{\rmin(\A\to\rmin\B)\to\C}\1\0,\,\fm\A\2\1\implies\fm{\B\to\C}\2\0
    &&\text{$|{\to}$,~no~$\3$}\\
    9.\quad&\fm{\rmin(\A\to\rmin\B)\to\C}\1\0\implies\fm{\A\to(\B\to\C)}\1\0
    &&\text{$|{\to}$,~no~$\2$}\\
    10.\quad&\implies\fm{(\rmin(\A\to\rmin\B)\to\C)\to(\A\to(\B\to\C))}\0\0
    &&\text{$|{\to}$,~no~$\1$}
  \end{align*}
  \endproof

\begin{lemma}\label{t?}
  $(\A\to\B)\to(\rmin(\A\to\C)\to\rmin(\B\to\C))$
  \begin{equation*}
    (\A\to\B)\to((\A\circ\D)\to(\B\circ\D))
  \end{equation*}
\end{lemma}
\proof
  \begin{align*}
    1.\quad&\fm\C\3\2\implies\fm\C\3\2
    &&\text{Axiom}\\
    2.\quad&\fm\B\3\0\implies\fm\B\3\0
    &&\text{Axiom}\\
    3.\quad&\fm\B\3\0,\,\fm{\B\to\C}\0\2\implies\fm\C\3\2
    &&\text{${\to}|$}\\
    4.\quad&\fm\A\3\1\implies\fm\A\3\1
    &&\text{Axiom}\\
    5.\quad&\fm{\A\to\B}\1\0,\,\fm{\B\to\C}\0\2,\,\fm\A\3\1\implies\fm\C\3\2
    &&\text{${\to}|$}\\
    6.\quad&\fm{\A\to\B}\1\0,\,\fm{\B\to\C}\0\2\implies
    \fm{\A\to\C}\1\2
    &&\text{$|{\to}$, no $\3$}\\
    7.\quad&\fm{\A\to\B}\1\0\implies
    \fm{\rmin(\B\to\C)}\2\0,\,\fm{\A\to\C}\1\2
    &&\text{$|\rmin$}\\
    8.\quad&\fm{\A\to\B}\1\0,\,\fm{\rmin(\A\to\C)}\2\1\implies
    \fm{\rmin(\B\to\C)}\2\0
    &&\text{$\rmin|$}\\
    9.\quad&\fm{\A\to\B}\1\0\implies\fm{\rmin(\A\to\C)\to\rmin(\B\to\C)}\1\0
    &&\text{$|{\to}$, no $\2$}\\
    10.\quad&\implies\fm{(\A\to\B)\to(\rmin(\A\to\C)\to\rmin(\B\to\C))}\0\0
    &&\text{$|{\to}$, no $\1$}
  \end{align*}
  \endproof

\begin{lemma}\label{assocfusion}
  $(\A\circ\B)\circ\C\to\A\circ(\B\circ\C)$
\end{lemma}
\proof
  \begin{align*}
    1.\quad &\fm{\B\to\rmin\C}\3\1\implies\fm{\B\to\rmin\C}\3\1
    &&\text{Axiom}\\
    2.\quad&\implies\fm{\rmin(\B\to\rmin\C)}\1\3,\,\fm{\B\to\rmin\C}\3\1
    &&\text{$|\rmin$}\\
    3.\quad&\fm{\rmin\rmin(\B\to\rmin\C)}\3\1\implies\fm{\B\to\rmin\C}\3\1
    &&\text{$\rmin|$}\\
    4.\quad&\fm\A\3\0\implies\fm\A\3\0
    &&\text{Axiom}\\
    5.\quad&\fm{\A\to\rmin\rmin(\B\to\rmin\C)}\0\1,\,\fm\A\3\0\implies\fm{\B\to\rmin\C}\3\1
    &&\text{${\to}|$}\\
    6.\quad&\fm\B\2\3\implies\fm\B\2\3
    &&\text{Axiom}\\
    7.\quad&\fm\C\1\2\implies\fm\C\1\2
    &&\text{Axiom}\\
    8.\quad&\fm\C\1\2,\,\fm{\rmin\C}\2\1\implies
    &&\text{$\rmin|$}\\
    9.\quad&\fm{\B\to\rmin\C}\3\1,\,\fm\B\2\3,\,\fm\C\1\2\implies
    &&\text{6, 8, ${\to}|$}\\
    10.\quad&\fm{\A\to\rmin\rmin(\B\to\rmin\C)}\0\1,\,\fm\A\3\0,\,\fm\B\2\3,\,\fm\C\1\2\implies
    &&\text{5, 9, Cut}\\
    11.\quad&\fm{\A\to\rmin\rmin(\B\to\rmin\C)}\0\1,\,\fm\C\1\2,\,\fm\A\3\0\implies\fm{\rmin\B}\3\2
    &&\text{$|\rmin$}\\
    12.\quad&\fm{\A\to\rmin\rmin(\B\to\rmin\C)}\0\1,\,\fm\C\1\2\implies\fm{\A\to\rmin\B}\0\2
    &&\text{$|{\to}$, no $\3$}\\
    13.\quad&\fm{\A\to\rmin\rmin(\B\to\rmin\C)}\0\1,\,\fm{\rmin(\A\to\rmin\B)}\2\0\implies\fm{\rmin\C}\2\1
    &&\text{$|\rmin$}\\
    14.\quad&\fm{\A\to\rmin\rmin(\B\to\rmin\C)}\0\1\implies\fm{\rmin(\A\to\rmin\B)\to\rmin\C}\0\1
    &&\text{$|{\to}$, no $\2$}\\
    15.\quad&\fm{\A\to\rmin\rmin(\B\to\rmin\C)}\0\1,\,\fm{\rmin(\rmin(\A\to\rmin\B)\to\rmin\C)}\1\0\implies
    &&\text{$\rmin|$}\\
    16.\quad&\fm{\rmin(\rmin(\A\to\rmin\B)\to\rmin\C)}\1\0\implies\fm{\rmin(\A\to\rmin\rmin(\B\to\rmin\C))}\1\0
    &&\text{$|\rmin$}\\
    17.\quad&\implies\fm{\rmin(\rmin(\A\to\rmin\B)\to\rmin\C)\to\rmin(\A\to\rmin\rmin(\B\to\rmin\C))}\0\0
    &&\text{$|{\to}$, no $\1$}\\
    18.\quad&\implies\fm{\rmin(\A\circ\B\to\rmin\C)\to\rmin(\A\to\rmin(\B\circ\C))}\0\0
    &&\text{def $\circ$}\\
    19.\quad&\implies\fm{(\A\circ\B)\circ\C\to\A\circ(\B\circ\C)}\0\0
    &&\text{def $\circ$}
  \end{align*}
  \endproof
\section{Deductive rules of $\mathcal{L}_4$}
The next 13 lemmas establish the derived rules of inference for
Tarski's relevance listed Table~\ref{table3}.
\begin{lemma}[adjunction]\label{adjunction}
  If $\A,\B\in\FL$ then $\A\land\B\in\FL$.
\end{lemma}
\proof If $\A$ and $\B$ have 4-proofs we may concatenate them
  and add one more sequent to get a 4-proof of $\A\land\B$, as
  follows.
\begin{align*}
  1.\quad&\implies\fm\A\0\0 &&\text{by some 4-proof}\\
  2.\quad&\implies\fm\B\0\0 &&\text{by some 4-proof}\\
  3.\quad&\implies\fm{\A\land\B}\0\0 &&\text{$|\land$}
\end{align*}
\endproof

\begin{lemma}[modus ponens]\label{modusponens}
  If $\A\to\B\in\FL$ and $\A\in\FL$ then $\B\in\FL$.
\end{lemma}
\proof Assume $\A$ and $\A\to\B$ have 4-proofs.  Concatenate
  a 4-proof of $\A$ with a 4-proof of $\A\to\B$ and continue the
  sequence as follows, obtaining a 4-proof of $\B$, showing
  $\B\in\FL$.
\begin{align*}
  1.\quad&\implies\fm{\A\to\B}\0\0 &&\text{by some 4-proof}\\
  2.\quad&\implies\fm\A\0\0 &&\text{by some 4-proof}\\
  3.\quad&\fm\B\0\0\implies\fm\B\0\0 &&\text{Axiom}\\
  4.\quad&\fm{\A\to\B}\0\0\implies\fm\B\0\0 &&\text{${\to}|$}\\
  5.\quad&\implies\fm\B\0\0 &&\text{1, 4, Cut}
\end{align*}
\endproof

\begin{lemma}[disjunctive syllogism]\label{disjunctivesyllogism}
  If $\A\lor\B\in\FL$ and $\rmin\A\in\FL$ then $\B\in\FL$.
\end{lemma}
\proof If $\A\lor\B$ and $\rmin\A$ have 4-proofs, they may be
  continued to obtain a 4-proof of $\B$.
\begin{align*}
  1.\quad&\implies\fm{\A\lor\B}\0\0&&\text{by a 4-proof}\\
  2.\quad&\implies\fm{\rmin\A}\0\0&&\text{by a 4-proof}\\
  3.\quad&\fm\A\0\0\implies\fm\A\0\0&&\text{Axiom}\\
  4.\quad&\fm\B\0\0\implies\fm\B\0\0&&\text{Axiom}\\
  5.\quad&\fm{\A\lor\B}\0\0\implies\fm\A\0\0,\fm\B\0\0&&\text{$\lor|$}\\
  6.\quad&\implies\fm\A\0\0,\fm\B\0\0&&\text{1, 5, Cut}\\
  7.\quad&\fm{\rmin\A}\0\0\implies\fm\B\0\0&&\text{$\rmin|$}\\
  8.\quad&\implies\fm\B\0\0&&\text{2, 7, Cut}
\end{align*}
\endproof

\begin{lemma}[transitivity]\label{transitivity}
  If $\A\to\B\in\FL$ and $\B\to\C\in\FL$ then $\A\to\C\in\FL$.
\end{lemma}
\proof Assume $\A\to\B\in\FL$ and $\B\to\C\in\FL$.  Then the
  sequents $\implies\fm{\A\to\B}\0\0$ and $\implies\fm{\B\to\C}\0\0$
  have 4-proofs that can be concatenated with steps 2--5 inserted
  between them, followed by sequents 7--12, yielding a 4-proof of
  $\A\to\C$, hence $\A\to\C\in\FL$.
  \begin{align*}
    1.\quad&\implies\fm{\A\to\B}\0\0
    &&\text{by a 4-proof}\\
    2.\quad&\fm\A\1\0\implies\fm\A\1\0
    &&\text{Axiom}\\
    3.\quad&\fm\B\1\0\implies\fm\B\1\0
    &&\text{Axiom}\\
    4.\quad&\fm{\A\to\B}\0\0,\fm\A\1\0\implies\fm\B\1\0
    &&\text{${\to}|$}\\
    5.\quad& \fm\A\1\0\implies\fm\B\1\0
    &&\text{1, 4, CUT}\\
    6.\quad&\implies\fm{\B\to\C}\0\0
    &&\text{by a 4-proof}\\
    7.\quad&\fm\B\1\0\implies\fm\B\1\0
    &&\text{Axiom}\\
    8.\quad&\fm\C\1\0\implies\fm\C\1\0
    &&\text{Axiom}\\
    9.\quad&\fm{\B\to\C}\0\0,\fm\B\1\0\implies\fm\C\1\0
    &&\text{${\to}|$}\\
    10.\quad&\fm\B\1\0\implies\fm\C\1\0
    &&\text{6, 9, Cut}\\
    11.\quad&\fm\A\1\0\implies\fm\C\1\0
    &&\text{5, 10, Cut}\\
    12.\quad&\implies\fm{\A\to\C}\0\0 &&\text{$|{\to}$,~no~$\1$}
\end{align*}
\endproof 

\begin{lemma}[contraposition]\label{contraposition}
  If $\A\to\B\in\FL$ then $\rmin\B\to\rmin\A\in\FL$.
\end{lemma}
\proof Assume $\A\to\B\in\FL$. By interchanging $\0$ and $\1$
  throughout any 4-proof of $\implies\fm{\A\to\B}\0\0$, we obtain a
  4-proof of $\implies\fm{\A\to\B}\1\1$, which may be continued as
  follows to obtain a 4-proof of $\rmin\B\to\rmin\A$.
\begin{align*}
  1.\quad&\implies\fm{\A\to\B}\1\1&&\text{by a $(\0\1)$-permuted 4-proof}\\
  2.\quad&\fm\A\0\1\implies\fm\A\0\1&&\text{Axiom}\\
  3.\quad&\fm\B\0\1\implies\fm\B\0\1&&\text{Axiom}\\
  4.\quad&\fm\A\0\1,\fm{\A\to\B}\1\1\implies\fm\B\0\1&&\text{${\to}|$}\\
  5.\quad&\fm\A\0\1\implies\fm\B\0\1&&\text{1, 4, Cut}\\
  6.\quad&\implies\fm\B\0\1,\,\fm{\rmin\A}\1\0&&\text{$|\rmin$}\\
  7.\quad&\fm{\rmin\B}\1\0\implies\fm{\rmin\A}\1\0&&\text{$\rmin|$}\\
  8.\quad&\implies\fm{\rmin\B\to\rmin\A}\0\0&&\text{$|{\to}$,~no~$\1$}
\end{align*}
\endproof 

\begin{lemma}[contraposition.2]\label{contraposition.2}
  If $\A\to\rmin\B\in\FL$ then $\B\to\rmin\A\in\FL$.
\end{lemma}
\proof Assume $\A\to\rmin\B$ has a 4-proof.  Obtain a 4-proof
  of $\implies\fm{\A\to\rmin\B}\1\1$ by interchanging $\0$ and $\1$ in
  a 4-proof of $\implies\fm{\A\to\rmin\B}\0\0$.  Continue this 4-proof
  as follows to obtain a 4-proof of $\B\to\rmin\A$.
\begin{align*}
  1.\quad&\implies\fm{\A\to\rmin\B}\1\1 &&\text{by a $(\0\1)$-permuted 4-proof}\\
  2.\quad&\fm\A\0\1\implies\fm\A\0\1 &&\text{Axiom}\\
  3.\quad&\fm\B\1\0\implies\fm\B\1\0 &&\text{Axiom}\\
  4.\quad&\fm{\rmin\B}\0\1,\fm\B\1\0\implies &&\text{$\rmin|$}\\
  5.\quad&\fm{\A\to\rmin\B}\1\1,\fm\B\1\0,\fm\A\0\1\implies
  &&\text{2, 4, ${\to}|$}\\
  6.\quad&\fm\B\1\0,\fm\A\0\1\implies &&\text{1, 5, Cut}\\
  7.\quad&\fm\B\1\0\implies\fm{\rmin\A}\1\0 &&\text{$|\rmin$}\\
  8.\quad&\implies\fm{\B\to\rmin\A}\0\0 &&\text{$|{\to}$,~no~$\1$}
\end{align*}
\endproof

\begin{lemma}[cut]\label{cut} 
  If $\A\land\B\to\C\in\FL$ and $\B\to\C\lor\A\in\FL$ then
  $\B\to\C\in\FL$.
\end{lemma}
\proof The Cut Rule in relevance logic is a derived rule in
  Basic Logic, called DR2~\cite[p.\,291]{Routleyetal1982}.  To prove
  this simplified version of DR2, construct a 4-proof of $\B\to\C$
  from 4-proofs of $\A\land\B\to\C$ and $\B\to\C\lor\A$ as follows.
  It is interesting that Cut for sequents is used five times.
  \begin{align*}
    1.\quad&\fm\B\1\0\implies \fm\B\1\0&&\text{Axiom}\\
    2.\quad&\fm{\C\lor\A}\1\0\implies \fm{\C\lor\A}\1\0&&\text{Axiom}\\
    3.\quad&\fm{\B\to\C\lor\A}\0\0,\fm\B\1\0\implies \fm{\C\lor\A}\1\0&&\text{${\to}|$}\\
    4.\quad&\implies\fm{\B\to\C\lor\A}\0\0&&\text{by a 4-proof}\\
    5.\quad&\fm\B\1\0\implies\fm{\C\lor\A}\1\0&&\text{Cut}\\
    6.\quad&\fm\C\1\0\implies\fm\C\1\0&&\text{Axiom}\\
    7.\quad&\fm\A\1\0\implies\fm\A\1\0&&\text{Axiom}\\
    8.\quad&\fm{\C\lor\A}\1\0\implies\fm\C\1\0,\fm\A\1\0&&\text{${\lor}|$}\\
    9.\quad&\fm\B\1\0\implies\fm\C\1\0,\fm\A\1\0&&\text{5, 8, Cut}\\
    10.\quad&\fm{\A\land\B}\1\0\implies \fm{\A\land\B}\1\0&&\text{Axiom}\\
    11.\quad&\fm\C\1\0\implies \fm\C\1\0&&\text{Axiom}\\
    12.\quad&\fm{\A\land\B\to\C}\0\0,\fm{\A\land\B}\1\0\implies \fm\C\1\0&&\text{${\to}|$}\\
    13.\quad&\implies\fm{\A\land\B\to\C}\0\0&&\text{by a 4-proof}\\
    14.\quad&\fm{\A\land\B}\1\0\implies \fm\C\1\0	&&\text{Cut}\\
    15.\quad&\fm\A\1\0\implies \fm\A\1\0&&\text{Axiom}\\
    16.\quad&\fm\B\1\0\implies \fm\B\1\0&&\text{Axiom}\\
    17.\quad&\fm\A\1\0,\fm\B\1\0\implies \fm{\A\land\B}\1\0&&\text{$|{\land}$}\\
    18.\quad&\fm\A\1\0,\fm\B\1\0\implies \fm\C\1\0&&\text{14, 17, Cut}\\
    19.\quad&\fm{\B}\1\0\implies\fm\C\1\0&&\text{9, 18, Cut}\\
    20.\quad&\implies\fm{\B\to\C}\0\0&&\text{$|{\to}$,~no~$\1$}
  \end{align*}
  \endproof

\begin{lemma}[E-rule, BR1, R5]\label{E-rule}
  If $\A\in\FL$ then $(\A\to\B)\to\B\in\FL$.
\end{lemma}
\proof The E-rule~\cite[p.\,8]{Brady2003} is also called
  BR1~\cite[p.\,289]{Routleyetal1982} and
  R5~\cite[p.\,193]{Brady2003}.  If $\A$ has a 4-proof, then we obtain
  a 4-proof of $(\A\to\B)\to\B$ by appending sequents to a 4-proof of
  $\implies\fm\A\1\1$, as follows.
  \begin{align*}
    1.\quad&\implies\fm\A\1\1 &&\text{by a $(\0\1)$-permuted 4-proof}\\
    2.\quad&\fm\B\1\0\implies\fm\B\1\0 &&\text{Axiom}\\
    3.\quad&\fm{\A\to\B}\1\0\implies\fm\B\1\0 &&\text{${\to}|$}\\
    4.\quad&\implies\fm{(\A\to\B)\to\B}\0\0 &&\text{$|{\to}$,~no~$\1$}
  \end{align*}
  \endproof

\begin{lemma}[suffixing]\label{suffixing}
  If $\A\to\B\in\FL$ then $(\B\to\C)\to(\A\to\C)\in\FL$.
\end{lemma}
\proof Assume $\A\to\B$ is 4-provable.  Interchange $\0$ and
  $\1$ throughout a 4-proof of $\implies\fm{\A\to\B}\0\0$, obtaining a
  4-proof of $\implies\fm{\A\to\B}\1\1$, and continue it as follows to
  obtain a 4-proof of $(\B\to\C)\to(\A\to\C)$.
\begin{align*}
  1.\quad&\implies\fm{\A\to\B}\1\1
  &&\text{by a $(\0\1)$-permuted 4-proof}\\
  2.\quad&\fm\A\2\1\implies\fm\A\2\1&&\text{Axiom}\\
  3.\quad&\fm\B\2\1\implies\fm\B\2\1&&\text{Axiom}\\
  4.\quad&\fm{\A\to\B}\1\1,\fm\A\2\1\implies\fm\B\2\1&&\text{${\to}|$}\\
  5.\quad&\fm\C\2\0\implies\fm\C\2\0&&\text{Axiom}\\
  6.\quad&\fm{\B\to\C}\1\0,\fm\B\2\1\implies\fm\C\2\0
  &&\text{3, 5, ${\to}|$}\\
  7.\quad&\fm\A\2\1\implies\fm\B\2\1
  &&\text{1, 4, Cut}\\
  8.\quad&\fm{\B\to\C}\1\0,\fm\A\2\1\implies\fm\C\2\0
  &&\text{Cut}\\
  9.\quad&\fm{\B\to\C}\1\0\implies\fm{\A\to\C}\1\0
  &&\text{$|{\to}$,~no~$\2$}\\
  10.\quad&\implies\fm{(\B\to\C)\to(\A\to\C)}\0\0
  &&\text{$|{\to}$,~no~$\1$}
\end{align*}
\endproof

\begin{lemma}[cycling]\label{cycling}
  If $\A\to(\B\to\C)\in\FL$ then $\B\to(\rmin\C\to\rmin\A)\in\FL$
\end{lemma}
\proof If $\A\to(\B\to\C)$ is 4-provable, then there is a
  4-proof of $\implies\fm{\A\to(\B\to\C)}\2\2$, which may be
  incorporated into a 4-proof of $\B\to(\rmin\C\to\rmin\A)$ as
  follows.
  \begin{align*}
    1.\quad&\fm\A\0\2\implies\fm\A\0\2 &&\text{Axiom}\\
    2.\quad&\fm{\B\to\C}\0\2\implies\fm{\B\to\C}\0\2 &&\text{Axiom}\\
    3.\quad&\fm{\A\to(\B\to\C)}\2\2,\fm\A\0\2\implies\fm{\B\to\C}\0\2 &&\text{${\to}|$}\\
    4.\quad&\implies\fm{\A\to(\B\to\C)}\2\2 &&\text{by a $(\0\2)$-permuted 4-proof}\\
    5.\quad&\fm\A\0\2\implies\fm{\B\to\C}\0\2 &&\text{Cut}\\
    6.\quad&\fm\B\1\0\implies\fm\B\1\0 &&\text{Axiom}\\
    7.\quad&\fm\C\1\2\implies\fm\C\1\2 &&\text{Axiom}\\
    8.\quad&\fm{\B\to\C}\0\2,\fm\B\1\0\implies\fm\C\1\2 &&\text{${\to}|$}\\
    9.\quad&\fm\A\0\2,\fm\B\1\0\implies\fm\C\1\2 &&\text{5, 8, Cut}\\
    10.\quad&\fm\B\1\0\implies\fm\C\1\2,\fm{\rmin\A}\2\0&&\text{$|\rmin$}\\
    11.\quad&\fm\B\1\0,\fm{\rmin\C}\2\1\implies\fm{\rmin\A}\2\0&&\text{$\rmin|$}\\
    12.\quad&\fm\B\1\0\implies\fm{\rmin\C\to\rmin\A}\1\0&&\text{$|{\to}$,~no~$\2$}\\
    13.\quad&\implies\fm{\B\to(\rmin\C\to\rmin\A)}\0\0&&\text{$|{\to}$,~no~$\1$}
  \end{align*}
  \endproof 

\begin{lemma}[prefixing rule]\label{prefixingR}
  If $\A\to\B\in\FL$ then $(\C\to\A)\to(\C\to\B)\in\FL$.
\end{lemma}
\proof By Lemma~\ref{prefixingA},
  $(\A\to\B)\to((\C\to\A)\to(\C\to\B))\in\FL$, so if $\A\to\B\in\FL$
  then $((\C\to\A)\to(\C\to\B))\in\FL$ by Lemma~\ref{modusponens}.
\endproof

\begin{lemma}[affixing]\label{R3}
  If $\A\to\B\in\FL$ and $\C\to\D\in\FL$ then
  \begin{equation*}
    (\B\to\C)\to(\A\to\D)\in\FL.
  \end{equation*}
\end{lemma}
\proof By $\C\to\D\in\FL$ and Lemma~\ref{prefixingR},
  \begin{align*} 
    (\A\to\C)\to(\A\to\D)\in\FL.\\
    \intertext{ by $\A\to\B\in\FL$ and Lemma~\ref{suffixing},}
    (\B\to\C)\to(\A\to\C)\in\FL. \intertext{Hence, by
      Lemma~\ref{transitivity},} (\B\to\C)\to(\A\to\D)\in\FL.
\end{align*}
\endproof

\begin{lemma}[monotonic fusion]\label{monotonicfusion}
  If $\A\to\B\in\FL$ and $\C\to\D\in\FL$ then
  \begin{equation*}
    (\A\circ\C)\to(\B\circ\D)=
    \rmin(\A\to\rmin\C)\to\rmin(\B\to\rmin\D)\in\FL,
  \end{equation*}
\end{lemma}
\proof
\begin{align*}
  1.\quad&\A\to\B\in\FL
  &&\text{Assumption}\\
  2.\quad&\C\to\D\in\FL
  &&\text{Assumption}\\
  3.\quad&\rmin\D\to\rmin\C\in\FL
  &&\text{Lemma~\ref{contraposition.2}}\\
  4.\quad&(\B\to\rmin\D)\to(\A\to\rmin\D)\in\FL
  &&\text{1, Lemma~\ref{suffixing}}\\
  5.\quad&(\A\to\rmin\D)\to(\A\to\rmin\C)\in\FL
  &&\text{3, Lemma~\ref{prefixingR}}\\
  6.\quad&((\B\to\rmin\D)\to(\A\to\rmin\D))\to
  ((\B\to\rmin\D)\to(\A\to\rmin\C))\in\FL
  &&\text{Lemma~\ref{prefixingR}}\\
  7.\quad&(\B\to\rmin\D)\to(\A\to\rmin\C)\in\FL
  &&\text{4, 6, Lemma~\ref{modusponens}}\\
  8.\quad&\rmin(\A\to\rmin\C)\to\rmin(\B\to\rmin\D)\in\FL
  &&\text{Lemma~\ref{contraposition.2}}\\
  9.\quad&(\A\circ\C)\to(\B\circ\D)\in\FL
  &&\text{definition}
\end{align*}
\endproof


\end{document}